\documentclass[nosumlimits]{amsart}
\usepackage{amsmath}
\usepackage{amssymb}
\usepackage{layout}

\usepackage{graphicx}

\newtheorem{thm}{Theorem}[section]
\newtheorem{defin}[thm]{Definition}
\newtheorem{lemma}[thm]{Lemma}
\newtheorem{cor}[thm]{Corollary}

\newtheorem{exa}{Example}
\newtheorem{pro}[thm]{Proposition}
\numberwithin{equation}{section}

\newcommand{\N}{\mathcal{N}}  
\newcommand{\Power}{\mathcal{P}}  

\newcommand{\F}{\mathalpha{F}} 

 \newcommand{\into}{\rightarrow}
 \newcommand{\mto}{\mapsto}
 \newcommand{\impl}{\Longrightarrow}
 \newcommand{\less}{\setminus}
 \newcommand{\abs}[1]{\left\vert#1\right\vert}
 \newcommand{\set}[1]{\left\{#1\right\}}

 \newcommand{\norm}[1]{\left\Vert#1\right\Vert}
 
 \newcommand{\qtext}[1]{\quad\text{#1}\quad}
 \newcommand{\fa}{\qtext{for all}}
 \newcommand{\bb}{\begin{equation*}}
 \newcommand{\ee}{\end{equation*}}
 \newcommand{\bp}{\begin{proof}}
 \newcommand{\ep}{\end{proof}}

\begin{document}

\title[Generalized Lebesgue and Bochner integration theory]%
{Generalized vectorial Lebesgue and Bochner integration theory}

\author{ Victor M. Bogdan }

\address{Department of Mathematics, McMahon Hall 207, CUA, Washington DC 20064, USA}


\email{bogdan@cua.edu}


\begin{abstract}
This paper contains a development of the Theory
of Lebesgue and Bochner spaces of summable functions.
It represents a synthesis of the results due to
H.~Lebesgue, S.~Banach, S.~Bochner, G.~Fubini,
S.~Saks, F.~Riesz, N.~Dunford,
P.~Halmos,  and other contributors to this theory.

The construction of the theory is based on the notion of a measure
on a prering of sets in any abstract space X.
No topological structure of the space X
is required for the development of the theory.

Measures on prerings generalize the notion of  abstract Lebesgue
measures.
These measures are readily available and it is not
necessary to extend them beforehand onto
a sigma-ring for the development of the theory.

The basic tool in the development of the theory is the construction
and characterization of Lebesgue-Bochner spaces of summable functions
as in the paper of Bogdanowicz,
"A Generalization of the
Lebesgue-Bochner-Stieltjes Integral and a New Approach to the
Theory of Integration", Proc. of Nat. Acad. Sci. USA,
Vol.~53, No.~3, (1965), p.~492--498
\end{abstract}


\maketitle


\bigskip

\section{Introduction}

\bigskip

In this paper we shall present a development of the theory of
Lebesgue and Bochner spaces of summable functions and prove
the fundamental theorems of the theory.

The development of the integration theory beyond the classical
Riemann integral is important  for advancements
in  modern theory of differential equations,
theory of generalized functions, theory of
operators, probability, optimal control, and most of all
in theoretical physics.

Generalized functions introduced into mathematics by
P.~Dirac and put on precise mathematical footing
by L.~Schwartz \cite{schwartz}, I.~Gelfand
and G.~Shilov \cite{gelfand}, turned out to be essential in
analysis of flows of matter endowed with mass. For reference
see Bogdan \cite{bogdan72}.

This paper represents a partial synthesis of the results due to
H.~Lebesgue \cite{lebesgue}, S.~Banach \cite{banach}, S.~Bochner \cite{bochner},
S.~Saks \cite{saks}, F.~Riesz \cite{riesz}, N.~Dunford \cite{DS1},
and P.~Halmos \cite{halmos} and other contributors to this theory.

The construction of the theory is based on the notion of a measure
on a prering of any abstract space $X.$
No topological structure of the space $X$
is required for the development of the theory.
The notion of a prering represents an abstraction from the family of
intervals.

These measures are readily available and it is not
necessary to extend them beforehand for the development of the theory
to a Lebesgue type measure on a sigma-ring. They generalize Lebesgue
measures.

If $(X,V,v)$ represents such a measure space one can construct in a single
step the spaces $L(v,R)$ of Lebesgue summable functions and the spaces
$L(v,Y)$ of Bochner summable functions over the space $X$ and to develop
their properties as Banach spaces and to obtain the theory of Lebesgue
and Bochner integrals.

The main tools in developing of the theory are some elementary properties
of  Banach spaces concerning the norm $\norm{\ },$ some knowledge of calculus,
and familiarity with the geometric series
\begin{equation*}
 \sum_{n\ge0}x^n=1/(1-x)\fa x\in (0,1).
\end{equation*}

We shall follow the approach of Bogdanowicz \cite{bogdan10} and
\cite{bogdan14} with some modifications to construct a generalized
Lebesgue-Bochner-Stieltjes integral of the form $\int u(f,d\mu)$
where $u$ is a bilinear operator acting in the product of Banach
spaces, $f$ is a Bochner summable function, and $\mu$ a vector-valued
measure.

If the operator $u$ represents the multiplication by scalars that is
$u(\lambda,y)=\lambda\,y$ for $y$ in a Banach space $Y$ and $\mu$
represents a Lebesgue measure on a sigma ring, then the integral $\int u(f,d\mu)$
reduces to the classical Bochner integral $\int f\,d\mu$ and when
$Y$ represents the space of reals then the integral reduces to
the classical Lebesgue integral.

Every linear continuous operator from the space of Bochner summable
functions into any Banach space is representable by means of such
integral.

\bigskip


\bigskip

\section{Prerings and rings of sets}

\bigskip

Assume that $X$  is any abstract space and  $V$ some family of
subsets of $X$ that includes the empty set $\emptyset.$
Denote by $S(V)$ the family of all sets of $X$ that are disjoint unions
of finite collections of sets from the family $V.$
Clearly  empty set $\emptyset$ belongs to $S(V).$
Such a family $S(V)$ will be called the {\bf family of simple
sets} generated by the family $V.$
\bigskip

We shall say that such a family $V$ forms a {\bf prering} of the
space $X,$ if the following conditions are satisfied for any two sets
from $V$: if $A_{1},
A_{2}\in V$, then both the intersection $ A_{1}\cap A_{2}$ and set
difference $ A_{1}\less A_{2}$ belong to the family $S(V)$ of
simple sets.

The notion of a prering represents an abstraction
from the family of bounded intervals in the space of real numbers $R$
or  from the family of rectangles in the space $R^2.$

A family of sets $V$ of the space $X$ is called a {\bf ring} if
$V$ is a prering such that $V=S(V),$ which is equivalent to the
following conditions: if $A_{1}, A_{2}\in V$, then $ A_{1}\cup
A_{2}\in V$ and $ A_{1}\less A_{2}\in V$.

It is easy to prove that a family $V$ forms a prering if and only
if the family $S(V)$ of the simple sets forms a ring. Every ring
(prering) $V$ of a space $X$ containing the space $X$ itself is
called an {\bf algebra (pre-algebra)} of sets, respectively.

If the ring $V$ is  closed under countable unions it is called a
{\bf sigma ring} ($\sigma$-ring for short.) If the ring $V$ is
closed under countable intersections it is called a {\bf delta
ring} ($\delta$-ring for short.) It follows from de Morgan law
that $\delta$-algebra and $\sigma$-algebra represent the same
notion.

\bigskip


\bigskip

\section{Tensor product of prerings}

\bigskip

We shall need in the sequel the following notions.
Let $F=\set{A_i}$ be a nonempty finite family of sets of
the space $X.$ A finite family $G=\set{B_j}$ of disjoint sets is called a
{\bf refinement} of the family $F$ if every set of the family $F$ can be
written as disjoint union of some sets from the family $G.$

\bigskip

Given two abstract spaces $X_1$ and $X_2.$ Assume that $V_1$ is a collection
of subsets of the space $X_1$ and $V_2$ a collection of subsets of the space $X_2.$

By {\bf tensor product} $V_1\otimes V_2$ of the families $V_1$ and $V_2$  we shall
understand the family of subsets of the Cartesian product $X_1\times X_2$
defined by the formula
\begin{equation*}
 V_1\otimes V_2=\set{A_1\times A_2:\ A_1\in V_1,\ A_2\in V_2}.
\end{equation*}

\bigskip

\begin{thm}[Tensor product of prerings is a prering]
Assume that $V_j$ represents a prering of a space $X_j$ where  $j=1,2.$
Then the tensor product $V_1\otimes V_2$ represents a prering of the
space $X_1\times X_2.$
\end{thm}

\bigskip

\bp
    Notice that the following two properties of a family $V$ of subsets
    of a space $X$ are equivalent:
    \begin{itemize}
        \item The family $V$ forms a prering.
        \item The empty set belongs to $V$ and for every two sets $A,B\in V$
                there exists a finite disjoint refinement of the pair $\set{A,B}$
                in the family $V,$
                that is, there exists a finite collection $D=\{D_1,\ldots,D_k\}$
                of disjoint sets from $V$ such that each of the two sets $A,B$
                can be represented as a union of some sets from the collection $D.$
    \end{itemize}

    Clearly the tensor product $V_1\otimes V_2$ of the prerings contains the empty set.
    Now take any pair of sets $A,B\in V_1\otimes V_2.$ We have $A=A_1\times A_2$
    and $B=B_1\times B_2.$ If one of the sets  $A_1,A_2,B_1,B_2$ is empty then the
    pair $A,B$ forms its own refinement from $V.$ So consider the case when
    all the sets $A_1,A_2,B_1,B_2$ are nonempty.

    Let $C=\set{C_j\in V_1:\ j\in J}$ be a refinement of the
    pair $A_1,B_1$ and $$D=\set{D_k\in V_2:\ k\in K}$$ a refinement of $A_2,B_2.$
    We may assume that the refinements do not contain the empty set.

    The collection of sets $C\otimes D$
    forms a refinement of the pair $A,B.$
    Indeed each set of the pair $A_1,B_1$ can be uniquely
    represented as the union of sets
    from the refinement $C.$ Similarly each set of the pair $A_2,B_2$ can be represented
    in a unique way as union of sets from the refinement $D.$ Since
    the sets of the collection  $C\otimes D$ are disjoint and nonempty,
    each set of the pair
    $A_1\times A_2$ and $B_1\times B_2$ can be uniquely represented as the union
    of the sets from  $C\otimes D.$
    Thus $V=V_1\otimes V_2$ is a prering.
\ep

\bigskip


\bigskip

\section{Fundamental theorem on prerings}

\bigskip

The following theorem characterizes prerings of any abstract space $X.$ It shows that
prerings provide a natural sufficient structure to build the
theory of the integrals and of the spaces of summable functions.

\bigskip

\begin{thm}[Fundamental theorem on prerings] Let $V$ be a non-empty
family of subsets of an abstract space $X.$
Then the following statements are equivalent
\begin{enumerate}
    \item The family $V$ of sets forms a prering.
    \item Every finite family of sets $\set{A_1,A_2,\ldots,A_k}\subset V$
        has a finite refinement in the family $V.$
    \item For every finite collection of linear spaces $Y_1,Y_2,\ldots,Y_n, W$
        and any map
        $$u: Y_1\times Y_2\times\cdots Y_n\mto W$$
        preserving zero, that is such that $u(0,0,\ldots,0)=0,$ we have
        that the relations
        $$s_1\in S(V,Y_1),\,s_2\in S(V,Y_2),\,\ldots,\ s_n\in S(V,Y_n)$$
        imply $s\in S(V,W),$ where the function $s$ is defined by the formula
        \bb
            s(x)=u(s_1(x),s_2(x),\ldots,s_n(x))\fa x\in X.
        \ee
    \item The family $S(V)$ of simple sets generated by $V$ forms a ring
        of sets.
\end{enumerate}
\end{thm}

\bigskip

We will need the following lemma.
\begin{lemma}
Let $V$ be a prering of sets. If $B\in V$ and $C_1,C_2,\ldots,
C_k\in V$ then the set
$$B\less (C_1\cup C_2\cup \ldots \cup C_k)$$
is a simple set, that is it belongs to the family $S(V).$
\end{lemma}

\bigskip

    To prove the lemma use mathematical induction on $k$ and  the set
    identity
    $$B\less (C_1\cup C_2\cup \ldots \cup C_{k+1})=
    (B\less (C_1\cup C_2\cup \ldots \cup C_k))\less C_{k+1}.$$

\bigskip

\bp
    We shall prove the theorem using a circular argument.

\bigskip

    $1\impl 2.$ Assume that the family $V$ forms a prering. Consider
    any family
    \bb
        E=\set{B_1,\ldots,B_k}\subset V.
    \ee
    We shall apply induction with respect to $k.$ For $k=1$ the
    statement (2) is satisfied. Assume that for every family such that
    $k\le n$ there exists a finite refinement from the prering $V.$
    Consider a family
    \bb
        F=\set{B_1,\ldots,B_{n+1}}\subset V.
    \ee From the inductive assumption the family
    \bb
        \set{B_2,\ldots,B_{n+1}}
    \ee has a finite refinement
    \bb
        \set{C_1,\ldots,C_m} \subset V.
    \ee Notice that the following family
    \bb
        H=\set{C_j\cap B_1,\, C_j\less B_1,\, B_1\less
        (C_1\cup\ldots\cup C_m):\ j=1,2,\ldots,m}
    \ee forms a refinement of the family $F.$ Using the definition of
    a prering and the lemma one can prove that there exists a family
    $G\subset V$ forming a finite refinement of the family $H$ and
    consequently of the family $F.$ Thus we have proved the
    implication $1\impl 2.$
    \bigskip

\bigskip

    $2\impl 3.$ To prove  this implication take
    \bb
        s_j\in S(V,Y_j)\qtext{for} j=1,2,\ldots,n
    \ee
    We may assume that each of the functions is of the form
    \begin{equation}\label{1-rep}
        s_j=y_{j,1}c_{A_{j,1}}+
        \ldots+y_{j,m_j}c_{A_{j,m_j}}\qtext{for} j=1,2,\ldots,m_j
    \end{equation}
    where each of the sets $A_{j,i}$ is non-empty and
    $y_{j,i}\not=0.$ Let $\set{B_1,B_2,\ldots,B_t}\subset V$ be a
    finite refinement of the collection of sets
    \bb
        \set{ A_{j,i}:\ j=1,\ldots,n;\ i=1,\ldots,m_j}.
    \ee
    We may assume that each set $B_j$ is non-empty. It follows
    from formula (\ref{1-rep}) and from the definition of a  refinement that
    each of the functions $s_j$ is  constant on each set $B_i$ and is
    equal to $0$ outside the union $B_1\cup\ldots\cup B_t.$ Let
    \bb
        s_j(x)=z_{j,i}\qtext{when}x\in B_i;\ i=1,2,\ldots,t.
    \ee
    Let
    \bb
        w_i=u(z_{1,i},z_{2,i},\ldots,z_{n,i})\qtext{for} i=1,2,\ldots,t.
    \ee Then the composed functions $s(x)$
    \bb
        s(x)=u(s_1(x),s_2(x),\ldots,s_n(x))\fa x\in X
    \ee
    is of the form
    \bb
        s=w_1c_{B_1}+w_2c_{B_2}+\cdots+w_tc_{B_t}\in S(V,W).
    \ee
    This completes the proof of the implication $2\impl 3.$
    \bigskip

\bigskip

    $3\impl 4$ Notice the equivalence
    \bb
        A\in S(V)\Leftrightarrow c_A\in S(V,R).
    \ee Let $Y_1=Y_2=R$ be the space of reals. Define
    \bb
        u(r_1,r_2)=r_1+r_2-r_1r_2\fa r_1\in Y_1,\, r_2\in Y_2.
    \ee
    The function $u$ preserves zero $u(0,0)=0$ and we have
    \bb
        c_{A\cup B}(x)=u(c_A(x),c_B(x))=c_{A}(x)+c_{B}(x)-c_{A}(x)c_{B}(x)\fa x\in X.
    \ee
    Thus $c_{A\cup B}\in S(V,R)$ and so ${A\cup B}\in S(V).$ Thus
    the family $S(V)$ of simple sets is closed under the finite  union
    operation.

    To prove that it is closed under the difference of sets operation
    notice the identity
    $$c_{A\less B}=c_A-c_Ac_B.$$ Using the function $v$ defined by
    $$v(r_1,r_2)=r_1-r_1r_2\fa r_1\in Y_1,\, r_2\in Y_2$$
    similarly as before we can derive that the collection $S(V)$ of
    simple sets is closed under the difference of sets operation. Thus
    the implication $3\impl 4$ has been proved.
    \bigskip

\bigskip

    $4\impl 1.$ To prove this implication assume that $S(V),$ the
    collection of simple sets generated by the family $V,$ forms a
    ring of sets. Thus $S(V)$ is closed under finite union and set
    difference operations, that is
    \bb
        A,B\in S(V)\impl A\cup B,\, A\less B \in S(V).
    \ee
    Now from the set identity
    \bb
        A\cap B=(A\cup B)\setminus [(A\setminus B)\cup (B\setminus A)]
    \ee follows that the family $S(V)$ is closed under the operation
    of intersection of two sets.

    Since $V\subset S(V)$ the above properties imply
    \bb
        A,B\in V\impl A\cap B\in S(V)\qtext{and} A\less B\in S(V).
    \ee
    Thus the family $V$ forms a prering. This completes the proof
    of the implication $4\impl 1.$ Thus the theorem has been
    proved.
\ep

\bigskip

\begin{cor}[Case of a pre-algebra] Let $V$ be a non-empty
family of subsets of an abstract space $X.$
Then the following statements are equivalent
\begin{enumerate}
    \item The family $V$ of sets forms a pre-algebra.
    \item For every finite collection of linear spaces $Y_1,Y_2,\ldots,Y_n, W$
        and any map $u$ from the  Cartesian product $Y_1\times Y_2\times\cdots Y_n$
        into the space $W$  we have
        that the relations
        $$s_1\in S(V,Y_1),\,s_2\in S(V,Y_2),\,\ldots,\ s_n\in S(V,Y_n)$$
        imply $s\in S(V,W),$ where the function $s$ is defined by the formula
        \bb
            s(x)=u(s_1(x),s_2(x),\ldots,s_n(x))\fa x\in X.
        \ee
    \item The family $S(V)$ of simple sets generated by $V$ forms an algebra
        of sets.
\end{enumerate}
\end{cor}

\bp
    The proof is obvious and we leave it to the reader.
\ep


\bigskip


\bigskip

\section{Measure space}

\bigskip

A finite-valued function $v$ from a prering $V$ into $[0,\infty),$
 the non-negative reals, satisfying the following
implication
\begin{equation}\label{countable additivity}
    A=\bigcup_{t\in T}A_t\impl v(A)=\sum_{t\in T}v(A_t)
\end{equation}
for every set $A\in V,$ that can be decomposed into finite or
countable collection
\begin{equation*}
 A_t\in V\,(t\in T)
\end{equation*}
 of disjoint sets, will be called a
{\bf $\sigma$-additive positive measure.}

It was called a {\bf positive volume} in the earlier papers of
Bogdanowicz \cite{bogdan10}, \cite{bogdan14}, and \cite{bogdan23}.
Notice that since by definition every prering $V$
contains the empty set $\emptyset,$ from countable additivity
(\ref{countable additivity}) follows that $v(\emptyset)=0.$

By {\bf Lebesgue measure} over an abstract space $X$ we shall
understand any set function $v$ from a $\sigma$-ring $V$ of the
space $X$ into the extended non-negative reals $[0,\infty],$ that
satisfies the implication (\ref{countable additivity}) and has
value zero on the empty set $v(\emptyset)=0.$ We have to postulate
this explicitly to avoid the case of a trivial measure that is
identically equal to $\infty.$

\begin{defin}[Measure space]
\label{positive measure space}
 A triple $(X,V,v),$ where $X$ denotes an abstract space and
$V$  a prering of the space $X$ and $v$ a $\sigma$-additive
non-negative finite-valued measure on the prering $V,$ will be
called a {\bf measure space} or positive measure space.
\end{defin}

\bigskip

It is clear that every {\em finite Lebesgue measure} forms a positive
measure in our sense, and in the case when it has infinite values
by striping it of infinities we obtain a positive measure in our sense.

The infinite valued measures in the theory of integration
are like oversized tires on the wheels of a car.
They provide some convenience but they are not necessary.
This can be easily noticed in the theory of Bochner summable
functions: There no natural infinities in Banach spaces
and one has to discard all sets of infinite measure to construct the integral.
We shall discuss such measures in detail in a later section.

\bigskip


\bigskip

\section{Examples of measures on prerings}

\bigskip

It is good to have a few examples of the measure spaces at hand. The first
example corresponds to Dirac's $\delta$ function.

\bigskip

\begin{exa}[{\bf Dirac measure space}] Let $X$ be any abstract set and $V$
the family of all subsets of the space $X.$ Let $x_0$ be a fixed
point of $X.$ Let $v_{x_0}(A)=1$ if $x_0\in A$ and  $v_{x_0}(A)=0$
otherwise. Since $V$ forms a sigma ring the triple $(X,V,v)$ forms
in this case a Lebesgue measure space.
\end{exa}

\bigskip

\begin{exa}[{\bf Counting measure space}] Let $X$ be any abstract set and
$V$ consisting of the empty set and of all singletons
\begin{equation*}
 V=\set{\emptyset,\ \set{x}:\ x\in X}.
\end{equation*}
 Let $v(A)=1$ for all
singleton sets $A=\set{x}$ and $v(\emptyset)=0.$ The triple
$(X,V,v)$ forms a measure space that is not a Lebesgue measure
space.
\end{exa}

\bigskip

\begin{exa}[{\bf Striped Lebesgue measure space}]
Assume that $M$ is a sigma-ring of subsets of $X$ and $\mu$  is
any Lebesgue measure on $M$ finite valued or with infinite values.
Let $$V=\set{A\in
M:\,\mu(A)<\infty}.$$ Plainly $V$ forms  a prering.
Then restricting $\mu$ to $V$ yields a positive measure space
$(X,V,\mu.)$
\end{exa}

\bigskip

The most important measure space to the sequel is the following.

\begin{pro}[{\bf Riemann measure space}]
\label{Riemann measure space} Let $R$ denote the space of reals
and $V$ the collection of all bounded intervals $I$ open, closed,
or half-open. If $a\le b$ are the end points of an interval $I$
let $v(I)=b-a.$ Then the triple $(R,V,v)$ forms a measure space.
We shall call this space the {\bf Riemann measure space.}
\end{pro}

\bigskip

\bp
    The collection $V$ of intervals forms a prering. Indeed the intersection
    of any two intervals is an interval or an empty set. But empty set
    can be represented as an open interval $(a,a)=\emptyset.$
    The set difference of two intervals is either the union of two disjoint intervals
    or a single interval or an empty set. Thus we have that for any two intervals
    $I_1,\,I_2\in V$ we have $I_1\cap I_2\in S(V)$ and $I_1\less I_2\in S(V).$
    This proves that $V$ is a prering.

    To prove countable additivity
    assume that we have a decomposition of an interval $I$ with ends $a\le b$
    into disjoint countable collection $I_t(t\in T)$ of intervals with end
    points $a_t\le b_t,$ that is
    \begin{equation}\label{decomposition}
        I=\bigcup_{t\in T}I_t.
    \end{equation}
     The case when interval $I$ is empty or consists
    of a single point is obvious. So without loss of generality we may assume
    that the interval $I$ has a positive length and that
    our index set $T=\set{1,2,3,\ldots}.$ Take any $\varepsilon>0$ such that
    $2\varepsilon<v(I).$ Let $I^\varepsilon=[a+\varepsilon,b-\varepsilon]$ and
    $I^\varepsilon_t=(a_t-\varepsilon2^{-t},b_t+\varepsilon2^{-t})$ for all $t\in T.$

    The family $I^\varepsilon_t(t\in T)$ forms an open cover of the compact
    interval $ I^\varepsilon$ thus there exists a finite set $J\subset T$ of indexes
    such that
    \begin{equation*}
        I^\varepsilon\subset \bigcup_{t\in J} I^\varepsilon_t.
    \end{equation*}
    The above implies
    \begin{equation*}
        \begin{split}
        v(I)-2\varepsilon&=v(I^\varepsilon)\le \sum_{t\in J}v(I^\varepsilon_t)
        \le \sum_{t\in T}v(I^\varepsilon_t)\\
        &=\sum_{t\in T}(v(I_t)+2^{-t+1}\varepsilon)=
        \sum_{t\in T}v(I_t)+2\varepsilon.
        \end{split}
    \end{equation*}
    Passing to the limit in the above inequality when $\varepsilon\into 0$
    we get
    \begin{equation*}
        v(I)\le \sum_{t\in T}v(I_t)
    \end{equation*}
    On the other hand from the relation (\ref{decomposition})  follows that
    for any finite set $J$ of indexes we have
    \begin{equation*}
        I\supset \bigcup_{t\in J}I_t\impl v(I)\ge \sum_{t\in J}v(I_t).
    \end{equation*}
    Since $$\sup_J \sum_{t\in J}v(I_t)=\sum_{t\in T}v(I_t)$$ we get from
    the above relations that the set function $v$ is countably additive
    and thus it forms a measure.
\ep

\bigskip

As will follow from the development of this theory the Riemann
measure space generates the same space of summable functions and
the integral as the classical Lebesgue measure over the reals. It
is good to see a few more examples of measures related to this
one.

\bigskip

\begin{pro}[{\bf Riemann measure space over $R^n$}]
\label{Riemann measure space over R^n} Let $X=R^n$
and $V$ the collection of all n-dimensional cubes of the form
\begin{equation*}
 I_1\times\ldots\times I_n
\end{equation*}
where $I_j$  represent intervals in the space $R$ of reals.
If $a_j\le b_j$ are the end points of an interval $I_j$
let
\begin{equation*}
 v(I_1\times\ldots\times I_n)=(b_1-a_1)(b_2-a_1)\cdots(b_n-a_n).
\end{equation*}
 Then the triple $(R^n,V,v)$ forms a measure space.
We shall call this space the {\bf Riemann measure space over $R^n.$}
\end{pro}

\bigskip

\bp
    The proof is similar to the preceding one and we leave it to the reader.
\ep

\bigskip

\begin{pro}[{\bf Stieltjes measure space}]
Let $R$ denote the space of reals, and $g$  a nondecreasing
function from $R$ into $R,$ and $D$ the set of discontinuity
points of $g.$ Let $V$ denote the collection of all bounded
intervals $I$ open, closed, or half-open with end points
$a,b\not\in D.$ If $a\le b$ are the end points of an interval $I$
let $v(I)=f(b)-f(a).$ Then the triple $(R,V,v)$ forms a measure
space.
\end{pro}

\bigskip

\bp
    The proof is similar to the proof for Riemann measure and we leave it to the reader.
\ep

\bigskip

A nondecreasing left-side continuous function $F$ from the
extended closed interval $E=[\,-\infty,+\infty\,]$ such that
$F(-\infty)=0$ and $F(+\infty)=1$ is called a {\bf probability
distribution function.} Any measure space $(X,V,v)$ over a prering
$V$ such that $X\in V$ and $v(X)=1$ is called a {\bf probability
measure space.}

\bigskip

\begin{pro}[Probability distribution generates probability measure space]
Let $F$ be a probability distribution on the extended reals $E.$
Let $V$ consists of all intervals of the form $[\, a,b)$ or
$[\,a,\infty\,],$ where $a,b\in E.$ If $I\in V$ let $v(I)$ denote the
increment of the function on the interval $I$ similarly as in the
case of Stieltjes measure space.

Then the triple $(E,V,v)$ forms a probability measure space.
\end{pro}

\bigskip

\bp
    To prove this proposition notice that the space $E$ can be considered
    as compact space and the proof can proceed similarly as in the case of the
    Riemann measure space.
\ep

\bigskip

In the case of topological spaces there are two natural prerings
available to construct a measure space: The prering consisting of
differences $G_1\less G_2$ of open sets, and the prering
consisting of differences $Q_1\less Q_2$ of compact sets.

\bigskip


\bigskip

\section{Constructing new measure spaces from old}

\bigskip

Given some measure spaces we shall review here several construction allowing one to obtain
new measure spaces.
Assume that we have available measure spaces $(X_t,V_t,v_t)$ where the index $t$
runs through a set $T$ of any cardinality.
We can construct a direct sum of the sets $X_t\,(t\in T)$ by considering
the disjoint union $X$ of the sets $X_t\,(t\in T).$ Thus each space $X_t\subset X$ and
\begin{equation*}
 X_t\cap X_s=\emptyset\qtext{if}t\not=s;\, t,s\in T.
\end{equation*}
Therefore the union
\begin{equation*}
 V=\bigcup_{t\in T}V_t
\end{equation*}
consists of subsets of the space $X.$ It is evident that the family $V$
forms a prering
and the measure $v$ on it is uniquely defined by the condition
\begin{equation*}
 v(A)=v_t(A)\qtext{if}A\in V_t,\qtext{for some}t\in T.
\end{equation*}
This measure space $(X,V,v)$ will be called a {\bf direct sum of measure spaces}
$(X_t,V_t,v_t)$ over $t\in T.$

\bigskip

Now consider the case when all the underling spaces $X_t=X$ coincide.
Define
\begin{equation*}
 V=\set{A\subset X:\ A\in V_t\ \forall\ t\in T,\ v(A)=\sum_{t\in T}v_t(A)<\infty}.
\end{equation*}
It is easy to prove that the triple $(X,V,v)$ forms a measure space.

\bigskip

Another operation yielding a new measure space is the following.
Consider a fixed measure space $(X,V,v)$ and subfamily $V_0\subset V$
forming a prering. Then the restriction $(X,V_0,v)$ yields a measure space.
A common example is the following. Take any set $A$ in the prering
$V$ and let
\begin{equation*}
 V_A=\set{B\in V:\ B\subset A}.
\end{equation*}
 restriction $(A,V_A,v)$ yields a measure space.
Thus, for instance, the Riemann measure space restricted to
to any interval in $R$ yields a measure over a pre-algebra.

\bigskip

Another important case is the following. Consider a pair
of measure spaces $(X_i,V_i,v_i)$ for $i=1,2.$
Let $X=X_1\times X_2$ denote the Cartesian product and
$V=V_1\otimes V_2$ the tensor product of the prerings.
Let
\begin{equation*}
 v(A)=v_1(A_1)v_2(A_2)\fa A=A_1\times A_2\in V.
\end{equation*}
For shorthand we shall use the notation $v_1\otimes v_2$
to denote the set function $v.$ Then the triple
\begin{equation*}
 (X,V,v)=(X_1\times X_2,V_1\otimes V_2, v_1\otimes v_2)
\end{equation*}
forms a measure space called the {\bf tensor product} of the measure spaces.
This fact follows from the {\em Monotone Convergence Theorem,}
also known as Beppo Levi's theorem, and we shall establish it
in a later section.

\bigskip


\section{Mathematical Preliminaries}

This section contains review of the essential notions concerning construction
of Banach spaces \cite{banach} that are needed for understanding the development of the
theory presented here.

A reader familiar with these notions can skip this section at the first
reading and to return to it later as the need arises.

\subsection{Construction of Basic Banach spaces}

Let $\F$ denote either the Galois field $R$ of reals or the Galois field $C$
of complex numbers. We shall assume that the reader is familiar with
the notion a {\bf linear space} also called a {\bf vector space.}
In short in a vector space $X$ there are
given two operations: addition of vectors $x+y$ for any $x,y,$
and multiplication by scalars $\lambda x$ for any $\lambda\in \F$ and $x\in X.$
The basic examples of vector spaces are
\begin{itemize}
    \item Space $R$ of reals considered as a vector space over the field $R.$
    \item Space $C$ of complex numbers considered as a vector space over the field $R.$
    \item Space $C$ of complex numbers considered as a vector space over the field $C.$
    \item More generally the space $R^n$ or $C^n$ for $n=1,2,\ldots$
    \item The space of continuous functions from an interval $I$ into
            a space $R^n.$
          This space will be denoted by $C(I,R^n).$
\end{itemize}

\bigskip

We shall need in this paper methods of constructing Banach spaces from the basic
spaces.

\bigskip

\subsection{Topological space}

\bigskip

\begin{defin}
\label{Topological space}
Let $X$ be an abstract space. Assume that in the space $X$ there is given a collection
of subsets $G$ having the following properties
\begin{itemize}
    \item The empty set $\emptyset$ and the entire space $X$ belong to $G.$
    \item For any finite number of sets $A_i\in G,$
            $(i=1,\ldots,k)$ their intersection $\bigcap_{i\le k} A_i$ belongs to $G.$
    \item For any finite or infinite number of sets $A_i\in G,$ where
            $i\in T,$ their union $\bigcup_{i\in T} A_i$ belongs to $G.$
\end{itemize}
Then the pair $(X,G)$ is called a {\bf topological space.} The sets belonging to $G$ are
called {\bf open} sets and the sets, whose complements $X\less A$ are open, are called
{\bf closed} sets. By a {\bf neighborhood} of a point $x$ one understands any open set
containing the point $x.$
\end{defin}

\bigskip

From de Morgan laws one can obtain the following characteristic properties of closed sets.

\bigskip

\begin{cor}[Closed sets]
\label{Closed sets}
Denote by $F$ the collection of all closed sets of a topological space $X.$
Then the following is true
\begin{itemize}
    \item The empty set $\emptyset$ and the entire space $X$ belong to $F.$
    \item For any finite number of sets $A_i\in F,$ where
            $i=1,\ldots,k,$ their union $\bigcup_{i\le k} A_i$ belongs to $F.$
    \item For any finite or infinite number of sets $A_i\in F,$
            $i\in T,$ their intersection $\bigcap_{i\in T} A_i$ belongs to $F.$
\end{itemize}
\end{cor}

\bigskip

\subsection{Seminormed and normed spaces}

\bigskip

\begin{defin}[Seminorm, extended seminorm, norm] Let $R^+=[ 0,\infty)$
    denote the set of non-negative reals and $E^+=[ 0,\infty] $ the set
    of extended non-negative reals. A functional $p$ from a linear
    space $X$ into $R^+$ or into $E^+$  is called a
    {\bf seminorm} or an {\bf extended seminorm,}
    respectively, if
    \begin{equation*}
        p(x+y)\le p(x)+p(y),\quad p(a x)=|a|p(x)\fa x,y\in X,\ a\in \F.
    \end{equation*}
    The first of the above conditions is called the {\bf triangle inequality}
    and the second the {\bf homogeneity} condition.
    Such a functional is called  a {\bf norm}
    if in addition $p(x)=0$ if and only if $x=0.$
    Such a functional is called nontrivial if $0<p(x)<\infty$ at least for one
    point $x\in X.$
\end{defin}

\bigskip

\begin{pro}[Seminorm inequalities]
\label{Seminorm inequalities}
If $p$ is a seminorm defined on a vector space $X,$
then the following inequalities are true
\begin{equation*}
    |p(x)-p(y)|\le p(x-y)\fa x,y\in X
\end{equation*}
and
\begin{equation*}
    p(x-y)\geq p(x)-p(y)\fa x,y\in X.
\end{equation*}

\end{pro}

\bigskip

\bp
    The proof follows from the definition of a seminorm in particular from
    the triangle inequality $p(x+y)\le p(x)+p(y).$
\ep

\bigskip


A pair $(X,p)$\label{Seminormed and normed spaces}
where $X$ is a vector space and $p$ a finite-valued seminorm on $X$
forms a topological space if open sets are defined by the following condition.

A set $A$ is open if it has the following property,
for every element $x\in A$ there exists a number $r>0$ such that every point $y$ with
the property $p(y-x)<r$ belongs to $A.$

For the empty set this condition
is satisfied in the vacuum. It is easy to check that indeed such a family of
sets satisfies the axioms for topology.

We shall call such a pair $(X,p)$
a {\bf seminormed space.} If the functional $p$ forms a norm, then the pair $(X,p)$
is called a {\bf normed space.}

In this case we use a more convenient notation
\begin{equation*}
    \norm{x}=p(x)\fa x\in X.
\end{equation*}

\bigskip

The basic vector spaces $R$ and $C$ with the absolute value $\abs{\ \,}$
form a normed-space. The spaces $R^n$ and $C^n$ form normed spaces with
norm
defined by the usual formula
\begin{equation*}
    \norm{(x_1,x_2,\ldots,x_n)}=(\sum_{j\leq n}|x_j|^2)^{1/2}\fa
    x=(x_1,x_2,\ldots,x_n)\in R^n \text{ or }C^n.
\end{equation*}
In the space $C(I,R^n)$ we define the norm by
\begin{equation*}
    \norm{f}=\sup\set{|f(t)|:\ t\in I}\fa f\in C(I,R^n).
\end{equation*}
These spaces are all normed spaces. To see some example of
a seminormed space that is not a normed space consider take $R^2$
and let
\begin{equation*}
 p(x_1,x_2)=|x_1|\fa(x_1,x_2)\in R^2.
\end{equation*}

\bigskip

\subsection{Lipschitz condition}

\bigskip

We use function as synonym with functional, map, operator, or transformation.

\bigskip

\begin{defin}
\label{Lipschitz condition}
Assume that $(X,p)$ and $(Y,q)$ are two semi-normed spaces.
A function $f$ from a set $W\subset X$ into the space $Y$
is said to be {\bf Lipschitzian} on the set $W$ if for some constant $M$
we have
\begin{equation*}
    q(f(x_1)-f(x_2))\le M\,p(x_1-x_2)\fa x_1,x_2\in W.
\end{equation*}
\end{defin}

\bigskip

\begin{defin}[Continuity]
\label{Continuity}
A function $f,$
from a topological
space $X$ into a topological space $Y$ is said to be {\bf continuous} if for every
open set $A$ in the space $Y$ the set
\begin{equation}\label{continuous transformation}
    B=\set{x\in X:\ f(x)\in A}
\end{equation}
is open in the topological space $X.$ The relation between the sets $A$ and $B$ given in
(\ref{continuous transformation}) we write as $B=f^{-1}(A)$ and in such a
case the set $B$ is called the {\bf inverse image} of the set $A.$
\end{defin}

\bigskip

It is easy to prove that every Lipschitzian function is continuous. In particular
every seminorm $p$ is continuous when considered as a function from the seminormed
space $(X,p)$ into the space $R$ of reals.

Using the above definition one can prove that a function $f$  is continuous if and
only if the inverse image $f^{-1}(F)$ of every closed set $F$ is closed.

\bigskip

\subsection{Convergence of a sequence or a series}

\bigskip

In a semi-normed vector space $(X,p)$ we introduce the notion of a sequential
convergence.

\begin{defin}
\label{Convergence of a sequence or a series}
We say
that a sequence $x_n\in X$ {\bf converges} to a point $x\in X$ if and only if
\begin{equation*}
 p(x-x_n)\into 0\qtext{as}n\into \infty.
\end{equation*}

The limit point $x$ in a semi-normed space is unique if and only if
the functional $p$ forms a norm.

By a {\bf series} with the terms $x_n\in X$
we understand the sequence $$s_n=\sum_{j\le n}x_j.$$ A series $x_n$ is said to
be {\bf convergent absolutely} with respect to the seminorm $p$ if
\begin{equation*}
    \sum_{n=1}^\infty p(x_n)<\infty.
\end{equation*}
\end{defin}

\bigskip

In a semi-normed space we have a convenient characterization of closed sets.
A set $F$ in such a space is closed if and only if for every convergent
sequence $x_n$ with terms in $F$ the limit $x$ of the sequence also is in $F,$
that is the following implication is true
\begin{equation*}
 x_n\in F\fa n\qtext{and}x_n\to x\impl x\in F.
\end{equation*}

This gives us another characterization of continuous functions in such spaces:
A function $f$ from a set $W$ in a semi-normed space $X$ into a semi-normed
space $Y$ is continuous if and only if for every sequence $x_n\in W$ that converges
to a point $x\in W$ we must have that the values $f(x_n)$ converge to the value $f(x).$

Another important property of semi-normed spaces is the following.
For every absolutely convergent series with terms $x_n$ that converges to a point $x\in X$
we have the estimate
\begin{equation*}
    p(x-\sum_{j=1}^n x_j)\le \sum_{j=n+1}^\infty p(x_j)\fa n=1,2,\ldots
\end{equation*}

\bigskip

\subsection{Banach spaces and complete seminormed spaces}

\bigskip

\begin{defin}
\label{Banach spaces and complete seminormed spaces}
A normed space $(X,\norm{\ }),$ in which every absolutely convergent series
$x_n$ is convergent to
some point $x\in X,$ is called a {\bf Banach space.}
\end{defin}

In the theory of integration it is convenient to consider semi-normed spaces that are
not normed but still every absolutely convergent series is convergent to a point.

Such
spaces will be called in the sequel the {\bf complete semi-normed spaces} and
the corresponding seminorms, {\bf complete seminorms.}

\bigskip

Notice that the basic vector spaces that we introduced so far are Banach spaces.

\bigskip

\begin{defin}[Lower semi-continuity] If $f$ is a function from a topological
space $X$ into the extended reals $E=[-\infty,\infty] $ then we say
that it is lower semi-continuous if for every finite number $a\in R$ the set
\begin{equation*}
    \set{x\in X:\ f(x)>a}=\set{x\in X:\ f(x)\in (a,\infty] }
\end{equation*}
is open or equivalently, by taking the complements of the sets, that the set
\begin{equation*}
    \set{x\in X:\ f(x)\le a}=\set{x\in X:\ f(x)\in [-\infty,a] }
\end{equation*}
is closed.
\end{defin}

\bigskip

Notice that every continuous function is lower semi-continuous.

\bigskip


\begin{pro}
Assume  that $(X,\norm{\ })$ is a normed vector space and $p_{\,t}\,(t\in T)$ a family
of extended seminorms from $X$ into $E^+=[ 0,\infty] .$
If each seminorm $p_{\,t}$ is lower semi-continuous then
\begin{equation*}
    p(x)=\sup\set{p_{\,t}(x):\ t\in T}\fa x\in X
\end{equation*}
represents a lower semi-continuous extended seminorm on $X.$
\end{pro}

\bigskip

\bp
    It is plain that $p$ forms a seminorm.
    For any $a\in R$ notice the identity
    \begin{equation*}
        \set{x\in X:\ p(x)\le a}=\bigcap_{t\in T}\set{x\in X:\ p_{\,t}(x)\le a}.
    \end{equation*}
    Thus the set on the left is closed as an intersection of a family of closed sets.
    Hence the seminorm $p$ is lower semi-continuous.
\ep

\bigskip

\begin{thm}[Generating Banach space $X_p$]
\label{Generating Banach space $X_p$}
Assume that $(X,\| \ \|)$ represents a Banach space and $p:X\into E^+$
is a lower semi-continuous extended semi-norm.
Let
$$
        X_p=\set{x\in X:\, p(x)<\infty }\qtext{and}
        \norm{x}_p=\norm{x}+p(x)\text{ for }x\in X_p
$$
Then  the pair $(X_p,\norm{\ }_p)$ forms a Banach space.
\end{thm}

\bigskip

\bp
    It is easy to see that the space $X_p$ is linear and $\norm{\ }_p$ forms a norm
    on it. We need to prove completeness of the norm.
    Take any series $x_n\in X_p$ that is absolutely convergent
    $\sum_n\norm{x_n}_p<\infty$ and let
    $a$ be a number such that $$\sum_n\norm{x_n}_p<a<\infty.$$
    Since $\norm{x_n}\le \norm{x_n}_p$ and  $p(x_n)\le \norm{x_n}_p$ for all $n$
    we have the estimates
    \begin{equation*}
        \sum_n \norm{x_n}<a\qtext{and} \sum_n p_n(x_n)<a.
    \end{equation*}
    Since by assumption the norm $\norm{\ }$ is complete, there exists a point
    $x\in X$ such that
    \begin{equation}\label{norm limitofseries}
        \norm{x-\sum_{j\le n}x_j}\le\sum_{j>n}\norm{x_j}.
    \end{equation}
    Let us introduce notation
    \begin{equation*}
        B_r=\set{x\in X:\ p(x)\le r}\fa r\ge 0.
    \end{equation*}
    Notice that $B_r\subset X_p$ and that from lower semi-continuity of $p$
    follows that the sets $B_r$ are closed.
    Let
    \begin{equation*}
        r_n=\sum_{n<j}p(x_j)\fa n=1,2,\ldots
    \end{equation*}
    Since
    \begin{equation*}
        p({\sum_{n<j \le m}x_j})\le\sum_{n<j \le m}p(x_j)\le\sum_{n<j}p(x_j)= r_n,
    \end{equation*}
    we have
    \begin{equation*}
        {\sum_{n<j\le m}x_j}\in B_{r_n}\fa m>n\text{ and }n=1,2,\ldots
    \end{equation*}
    For fixed $n$ define the sequence $y_m$ by
    \begin{equation*}
        y_m=\sum_{n<j\le m}x_j\fa m>n\qtext{and}y_m=0\qtext{otherwise.}
    \end{equation*}
    Notice that when $m\into \infty,$ then
    \begin{equation*}
        y_m=\sum_{n<j\le m}x_j=\sum_{j\le m}x_j-\sum_{j\le n}x_j
        \into\sum_{j=1}^\infty x_j-\sum_{j\le n}x_j,
    \end{equation*}
    where the convergence is in the norm $\norm{\ }.$
    Since $y_m\in B_{r_n}$ for all $m,$ and
    by lower semi-continuity of $p$ the set $B_{r_n}$
    is closed, we must have
    \begin{equation*}
        x-\sum_{j\le n}x_j=\sum_{j=1}^\infty x_j-\sum_{j\le n}x_j\in B_{r_n}\fa n,
    \end{equation*}
    that is
    \begin{equation*}
        p(x-\sum_{j\le n}x_j)\le r_n=\sum_{n<j}p(x_j)\fa n=1,2,\ldots
    \end{equation*}
    Since all $x_j\in X_p$ the above implies that $x\in X_p$ as follows from the
    triangle inequality.
    Combining the above with (\ref{norm limitofseries}) we get
    \begin{equation*}
        \begin{split}
        \norm{x-\sum_{j\le n}x_j}_p&=\norm{x-\sum_{j\le n}x_j}+p(x-\sum_{j\le n}x_j)\\
           &\le \sum_{n<j}\norm{x_j}+\sum_{n<j}p(x_j)= \sum_{n<j}\norm{x_j}_p\fa n
        \end{split}
    \end{equation*}
    Thus we have proved that every absolutely convergent
    series in the normed space $(X_p,\norm{\ }_p)$ converges to some point
    $x\in X_p.$ Hence the pair $(X_p,\norm{\ }_p)$ forms a Banach space.
\ep

\bigskip

\begin{thm}[Space $C_\delta(J,Y)$]
\label{Space $C_delta(J,Y)$}
Assume that $J$ is any interval, bounded or unbounded, closed or open, or partially closed,
in the space $R$ of reals and $Y$ a Banach space.

Assume that $\delta$ is a non-negative continuous function from $J$ into $R.$
Assume that $C_\delta(J,Y)$ denotes the set of all continuous functions from
the interval $J$ into the Banach space $Y$ such that
\begin{equation*}
    |f(x)|\le M\delta(x)\fa x\in J\qtext{and some} M.
\end{equation*}
Let $\norm{f}$ denote the infimum, that is the greatest lower bound,
of all constants $M$ appearing in the above
condition.

Then the functional $\norm{\ }$ forms a norm and the set $C_\delta(J,Y)$ equipped with
the norm $\norm{\ }$ forms a Banach space.
\end{thm}
\bigskip

\bp
    Let us introduce a shorthand notation $C_\delta=C_\delta(J,Y).$
    It is clear that the set $C_\delta$ is non-empty since it contains
    for instance all functions of the form $\delta(x)y$ where $y$ is
    any element of $Y.$
    It follows from the definition of the infimum of a set that
    \begin{equation}\label{estimate of f in C_delta}
        |f(x)|\le \delta(x)\norm{f}\fa f\in C_\delta,\ x\in J.
    \end{equation}
    Let us prove that $C_\delta$ forms a linear space and $\norm{\ }$
    forms a norm.

    To this end take any scalars $\lambda_1,\lambda_2$ and functions $f_1,f_2\in C_\delta.$
    We have
    \begin{equation*}
        \begin{split}
        |\lambda_1f_1(x)+\lambda_2f_2(x)|&\le|\lambda_1|\,|f_1(x)|+|\lambda_2|\,|f_2(x)|\\
            &\le|\lambda_1|\,\delta(x)\norm{f_1}+|\lambda_2|\,\delta(x)\norm{f_2}\\
            &\le\delta(x)(|\lambda_1|\,\norm{f_1}+|\lambda_2|\,\norm{f_2}).\\
        \end{split}
    \end{equation*}
    Thus from above we get $\lambda_1f_1+\lambda_2f_2\in C_\delta,$ so $C_\delta$
    is linear, and
    \begin{equation*}
        \norm{\lambda_1f_1+\lambda_2f_2}\le
        |\lambda_1|\norm{f_1}+|\lambda_2|\norm{f_2}\fa
        \lambda_1,\lambda_2\text{ and }f_1,f_2\in C_\delta.
    \end{equation*}
    The above inequality implies that $\norm{\ }$ forms a seminorm. Thus from the
    estimate (\ref{estimate of f in C_delta}) follows that the seminorm is
    zero if and only if the function $f$ is identically zero. Thus $\norm{\ }$ is
    a norm.

    To prove completeness take any absolutely convergent
    series $f_n\in C_\delta.$
    We have
    \begin{equation*}
        \sum_n\norm{f_n}<\infty
    \end{equation*}
    and
    \begin{equation*}
        |f_n(x)|\le \delta(x)\norm{f_n}\fa x\in J,\ n=1,2,\ldots
    \end{equation*}
    Since the function $\delta$ as a continuous function is bounded on
    every  bounded subinterval of $J,$ the series
    \begin{equation*}
        \sum_nf_n(x)=f(x)
    \end{equation*}
    converges uniformly and absolutely to some continuous function $f(x)$
    on such subinterval. Thus the limit  function $f$ is continuous.

    Since
    \begin{equation*}
        |f(x)|\le \sum_n|f_n(x)|\le \delta(x) \sum_n \norm{f_n}\fa x\in J,
    \end{equation*}
    the function $f$ belongs to the set $C_\delta.$
    Thus we have
    \begin{equation*}
        \begin{split}
        |f(x)-\sum_{j\le n}f_j(x)|
            &\le \delta(x)\norm{f-\sum_{j\le n}f_j}\\
            &\le\delta(x)\sum_{j>n}\norm{f_j}
                \fa x\in J,\ n=1,2,\ldots
        \end{split}
    \end{equation*}
    Therefore by definition of the norm $\norm{\ }$ we have
    \begin{equation*}
        \norm{f-\sum_{j\le n}f_j}\le\sum_{j>n}\norm{f_j}.
    \end{equation*}
    Thus the space $C_\delta$ forms a Banach space.
\ep

\bigskip

In the above theorem we did not assume about the function $\delta$ that it is bounded
or that it does not take on the value zero at some points. When
$\delta(x)=1$ for all $x\in J,$ we get that the space $C_\delta$ coincides with
the set of all continuous bounded functions on the interval $J$  and $C_\delta$
forms a Banach space. We will need a more general theorem for the case
of continuous bounded functions.

\bigskip

\begin{thm}[Space $C(W,Y)$]
\label{Space $C(W,Y)$}
Assume that $W$ is a topological space and $Y$ a Banach space.
Assume that $C(W,Y)$ denotes the set of all continuous bounded functions from
the space $W$ into the Banach space $Y.$
Let
\begin{equation*}
    \norm{f}=\sup\set{|f(x)|:\ x\in W}\fa f\in C(W,Y).
\end{equation*}
Then the set $C(W,Y)$ equipped with the above norm forms a Banach space.
\end{thm}

\bigskip

\bp
    The proof is similar to the previous one an we leave it to the reader.
\ep

\bigskip



\bigskip

\section{Linearity and bilinearity in Banach spaces}

\bigskip

Assume that $Y, Z, W$ represent some Banach spaces either over the
field $R$ of reals or over the field $C$ of complex numbers. By
linear operator $T,$ say from the space $Y$ into the space $W,$ we
shall understand any {\bf additive operator,} that is such that
\begin{equation*}
    T(y_1+y_2)= T(y_1)+T(y_2)\fa y_1,y_2\in Y,
\end{equation*}
and moreover it is either {\bf homogeneous}
\begin{equation*}
    T(\lambda y)= \lambda T(y)\fa y\in Y\text{ and }\lambda\in C
\end{equation*}
or {\bf conjugate homogeneous}
\begin{equation*}
    T(\lambda y)= \overline{\lambda}\, T(y)\fa y\in Y\text{ and }\lambda\in C.
\end{equation*}

In the case of Banach spaces over the field $C$ of complex numbers
every linear operator with respect to complex scalars is at the
same time a linear operators with respect to the real scalars.

By a bilinear operator $u$ from the product $Y\times Z$ into the
space $W$ we shall understand an operator such that the operators
\begin{equation*}
    y\mto u(y,z)\qtext{and}z\mto u(y,z)
\end{equation*}
are linear when the other variable is fixed.

As an example of such an operator consider Dirac's bra-ket
$\langle y|z\rangle$ operator. In this case $Y$ is a Hilbert space
over the field of complex numbers and the operator is given by the
formula
\begin{equation*}
    u(y,z)= \langle y|z\rangle\fa y,z\in Y,
\end{equation*}
where $u$ is from $Y\times Y$ into $C,$ and represents an inner
product operator. Such operator is additive in each variable and
homogeneous in variable $z$ and conjugate homogeneous in the
variable $y.$

We shall require the operator $u$ to be {\bf bounded} that is
to satisfy the following condition
\begin{equation*}
 |u(y,z)|\le m |y|\,|z|\fa y\in Y,\ z\in Z\text{ and some }m\in R.
\end{equation*}

The set of all such bilinear operators with the same type of
homogeneity with respect to complex numbers with usual operations
of addition of functions and multiplication of a function by a
number forms a Banach space with norm defined as $|u|=\inf\set{m}$
where $m$ runs through all constants satisfying the above condition.

Denote by $U$ the space of all bilinear bounded operators $u$
from the space $Y\times Z$ into $W$. Norms of elements in the
spaces $Y, Z, W, U$ will be denoted by $|\ \,|$.

\bigskip


\bigskip

\section{Classical approach to Lebesgue-Bochner integration}

\bigskip

The development of the classical Lebesgue-Bochner theory of the
integral goes through the following main stages as in Halmos
\cite{halmos} and Dunford and Schwartz \cite{DS1}:
\begin{itemize}
\item The construction and development of the Caratheodory theory
of outer measure $v^*$
        over an abstract space $X.$
\item The construction of the Lebesgue measure $v$ on the sigma
ring $V$ of measurable sets
        of the space $X$ induced by the outer measure $v^*.$
\item The development of the theory of real-valued measurable
functions $M(v,R).$ \item The construction of the Lebesgue
integral $\int f\,dv.$ \item The construction and development of
the theory of the space $L(v,R)$
        of Lebesgue summable functions.
\item The construction and development of the theory of the space
$M(v,Y)$ of
        Bochner measurable functions.
\item The construction of the Bochner integral $\int f\,dv$ and of
the space
        $L(v,Y)$ of Bochner summable functions $f$ from the space $X$ into
        any Banach space $Y.$
\end{itemize}

The construction of the classical Lebesgue integral is an abstraction from
the area under the graph of the function similar to the ideas of
Riemann though different in execution.

From the point of view of Functional Analysis both the Lebesgue
and Bochner integrals are particular linear continuous operators
from the space $L(v,Y)$ of Bochner summable functions into the
Banach space $Y.$ Moreover from the theory of the space $L(v,Y)$
one can easily derive the theory of the spaces $M(v,Y),$ $M(v,R),$
and $L(v,R)$ and of the Lebesgue and Bochner integrals and also
the theory of Lebesgue measure. For details see Bogdanowicz
\cite{bogdan14} and \cite{bogdan23}.

We shall show in brief how one can develop the theory of the space
$L(v,Y)$ and to construct an integral of the form $\int u(f,d\mu),$
where $u$ is any bilinear bounded operator from the product
$Y\times Z$ of Banach spaces into a Banach space $W$ and $\mu$
represents a vector measure. This integral for the case, when the
spaces $Y,Z,W$ are equal to the space $R$ of reals and the
bilinear operator $u$ represents multiplication $u(y,z)=yz$,
coincides with the Lebesgue integral
\begin{equation}\label{Lebesgue integral}
    \int f\,dv=\int u(f,dv)\fa f\in L(v,R).
\end{equation}
In the case, when $Y=W$ and $Z=R$ and $u(y,z)=zy$ represents the
scalar multiplication, the integral coincides with the Bochner
integral
\begin{equation}\label{Bochner integral}
    \int f\,dv=\int u(f,dv)\fa f\in L(v,Y).
\end{equation}

\bigskip


\bigskip

\section{Vector measures}

\bigskip

\begin{defin}[Vector measure space]
Assume that $(X,V,v)$ represents a positive measure space over the
prering $V.$
A set function $\mu$ from a prering $V$ into a Banach space $Z$ is
called a {\bf vector measure} if for every finite family of
disjoint sets $A_{t}\in V(t\in T)$ the following implication is
true
\begin{equation}\label{additivity}
    A=\bigcup_T A_{t}\in V\impl \mu(A)=\sum_T\mu(A_{t}).
\end{equation}
Denote by $K(v,Z)$ the space of all vector measures $\mu$ from the
prering $V$ into the space $Z$, such that
\begin{center}
$|\mu(A)|\leq m v(A)\fa A\in V$ and some $m$.
\end{center}
The least constant $m$ satisfying the above inequality is denoted
by $\norm{\mu}$. It is easy to see that the pair $(K(v,Z),
\norm{\mu})$ forms a Banach space.
\end{defin}

\bigskip


\bigskip

\section{Construction of the elementary integral spaces}

\bigskip

Assume that $c_A$ denotes the characteristic function of the set
$A$ that is $c_A(x)=1$ on $A$ and takes value zero elsewhere. Let
$S(V,Y)$ denote the space of all functions of the form
\begin{equation}\label{B}
    h=y_{1}c_{A_1}+\ldots+y_{k}c_{A_{k}},\text{ where }\
    y_{i}\in Y,\ A_{i}\in V
\end{equation}
and the sets $A_i$ in above formula are disjoint.

Notice that we extended the multiplication by scalars by
agreement $y\lambda=\lambda y$ for all vectors $y$ and scalars
$\lambda.$ The family $S(V,Y)$ of functions will be called the
family of {\bf simple functions} generated by the prering $V.$ For
fixed $u\in U$ and $\mu\in K(v,Z)$ and any simple function
$h\in S(v,Y)$ define the operator
\begin{center}
    $\int u(h,d\mu)=u(y_{1},\mu(A_{1}))+\ldots+u(y_{k},\mu(A_{k}))$
\end{center}
and the integral operator
\begin{center}
    $\int h\,dv=y_{1}v(A_{1})+\ldots+y_{k}v(A_{k})$.
\end{center}
The operators $\int h\,dv$ and $\int u(h,d\mu)$ are well defined,
that is, they do not depend on the representation of the function
$h$ in the form (\ref{B}). To prove this use the fact that
any finite collection of sets from a prering has a finite
refinement.

Let $|h|$ denote the function defined by the formula
$|h|(x)=|h(x)|$ for $x\in X$. We see that if $h\in S(V,Y)$, then
$|h|\in S(V,R)$. Therefore the following functional $\norm{h}=
\int|h|d\,v$ is well defined for all $h\in S(V,Y)$.

\bigskip

The following development of the theory of Lebesgue and Bochner
summable  functions and of the integrals  are  from Bogdanowicz
\cite{bogdan10}.

\begin{lemma}[Elementary integrals on simple functions]
\label{Elementary integrals} \label{Elementary lemma}
\label{Lemma1}
The following statements describe the basic
relations between the notions that we have just introduced.

\begin{enumerate}
\item   The space $S(V,Y)$  is linear, $||h||$
 is a seminorm on it, and $\int h\,dv$  is a
 linear operator on $S(V,Y),$  and moreover
\begin{equation*}
  \left|\int h\,dv\right|\leq\norm{h}\fa h\in S(V,Y).
\end{equation*}

\item   If $g\in S(V,R)$  and $f\in S(V,Y)$,  then $gf\in S(V,Y).$

\item  $\int h\, dv\geq 0$  if $h\in S(V,R)$  and $h(x)\geq 0$  for
all $x\in X$.

\item   $\int g\, dv\geq\int fdv$  if $g, f\in S(V,R)$ and $g(x)\geq
f(x)$  for all $x\in X$. \item  The operator $\int u(h,d\mu)$  is
trilinear from the product space $$U\times S(V,Y)\times K(v,Z)$$
into the space $W$  and
\begin{equation*}
    \left|\int u(h,d\mu)\right|\leq|u|\,\norm{h}\,\norm{\mu}
\end{equation*}
for all
\begin{equation*}
 u\in U,\ h\in S(V,Y),\ \mu\in K(v,Z).
\end{equation*}
\end{enumerate}
\end{lemma}

\bigskip

\bp
    The proof of the lemma is obvious and we leave it to the reader.
\ep

\bigskip


\bigskip

\section{Null sets}

\bigskip

Assume that $(X,V,v)$ represents a positive measure space over the
prering $V$ of an abstract space $X.$

\bigskip

\begin{defin}[Null sets]
Let $N$ be the family of all sets $A\subset X$ such that for every
$\varepsilon>0$ there exists a countable family $A_{t}\in V(t\in
T)$ such that
\begin{equation*}
 A\subset\bigcup_{T}A_{t}\qtext{ and} \sum_{T}v(A_{t})<\varepsilon.
\end{equation*}
 Sets of the family $N$ will be
called {\bf null-sets} with respect to the measure $v.$
\end{defin}

\bigskip

The sets of the family $N$ play the role of the sets of Lebesgue
measure zero and in fact this is how Lebesgue himself defined
this family in the case of sets in the space $R$ of reals.

This family represents a {\bf sigma-ideal} of sets
in the power set $\Power(X),$ that is, it has
the following properties: if $A\in N$, then $B\cap A\in N$ for any
set $B\subset X,$ and the union of any countable family of
null-sets $A_{t}\in N(t\in T)$ is also a null-set $\bigcup_T
A_{t}\in N$.

To see this fact take any $\varepsilon>0$ and assume
that $T=\set{1,2,\ldots}$ represents the set of natural numbers.
It is sufficient to cover each set $A_t$ with a countable
collection $Q_t$ of sets from the prering $V$ with total sum of measures
that does not exceed
\begin{equation*}
 2^{-t}\varepsilon\fa t=1,2,\ldots
\end{equation*}
Clearly the collection
\begin{equation*}
 Q=\bigcup_TQ_t
\end{equation*}
is countable, it covers the set $A=\bigcup_TA_t,$
and the total sum of the measures of sets in this collection
does note exceed $\varepsilon.$

A condition $C(x)$ depending on a parameter $x\in X$ is said to be
{\bf satisfied almost everywhere} if there exists a  null-set $A\in N$
such that the condition is satisfied at every
point of the set $X\less A$. 

\bigskip


\bigskip

\section{Null sets in case of Riemann measure space}

\bigskip

To see some examples of null sets consider the Riemann measure space $(R,V,v)$ over
the reals $R.$
Clearly the empty set $\emptyset=(a,a),$ and any singleton $[b,b]$
is a null. Moreover any countable set of points
forms a null set.

There exist null sets that are uncountable.
A typical example of such a set  is {\bf Cantor's set.}

To construct Cantor's set take the closed interval $[0,1]$
and divide it into three equal intervals. From the middle remove
the open interval $(1/3,2/3).$ The remaining two closed intervals
have total length $2/3$ and they form a closed set $F_1$. Repeat
this process with each of the remaining intervals.

After
n-steps the remaining set $F_n$ will consist of the union of $2^n$
disjoint closed intervals of total length of $(2/3)^n.$
The sets $F_n$ are nested and their intersection $F=\bigcup_nF_n$ will
represent a nonempty closed set called the Cantor set.

Cantor's set is a null set of cardinality equal to cardinality of
the interval $[0,1].$
Indeed, the set $F$ can be covered by a countable number
of intervals of total length as small as we please. Notice
that a finite cover by intervals we can always augment by
a sequence of intervals of the form $(a,a),$ that is by empty sets
to get a countable cover.

To prove that the set $F$ is of the same cardinality as the interval $[0,1]$
consider expansions into infinite fractions at the base $3$ of points belonging
to $F.$ Notice that the expansion must be of the form
\begin{equation*}
 x=0.\,d_1,d_2,d_3,\ldots
\end{equation*}
where the digits $d_i\in\set{0,2}.$
Ignore the set of points which have  periodic expansions since
they represent some rational numbers that form a countable set.

Similarly consider the binary expansions, that is at base $2,$  into non-periodic
sequence of digits of points $y$ of the set $[0,1].$ We have
\begin{equation*}
 y=0.\,a_1,a_2,a_3,\ldots
\end{equation*}
where $a_i\in\set{0,1}.$

Clearly the map $x\mapsto y$
given by the formula
\begin{equation*}
 a_i=d_i/2\fa i=1,2,\ldots
\end{equation*}
is one-to-one and onto. Thus the cardinalities of $F$
and $[0,1]$ are equal.
Hence Cantor's set $F$ is a null set of cardinality of a continuum.


\bigskip


\bigskip

\section{Fundamental lemmas}

\bigskip

\begin{defin}[Basic sequence]
By a {\bf basic sequence} we shall understand a sequence $s_{n}\in
S(V,Y)$ of simple functions for which there exists a series with terms $
h_{n}\in S(V,Y)$
 and a constant $M>0$ such that $s_{n}=h_{1}+h_{2}+\ldots+h_{n},$
 where  $||h_{n}||\leq M4^{-n}$ for all
$n=1,2, \ldots$
\end{defin}

\bigskip

The idea of the basic sequence as a series of simple functions with geometric rate
of convergence can be traced to the work of Riesz. See \cite{riesz}, page~59, the proof
of Riesz-Fisher theorem. The proof of Egoroff's theorem \cite{egoroff}
provided the idea for the rest of the needed structure. Thus it seems appropriate
to name the next lemma as Riesz-Egoroff property of a basic sequence.

\bigskip

\begin{lemma}[Riesz-Egoroff property of a basic sequence]
\label{Riesz-Egoroff property of basic sequence} \label{Basic
sequence lemma} \label{Lemma 2} Assume that $(X,V,v)$ is a
positive measure space on a prering $V$ and $Y$ is a Banach space.
Then the following is true.
\begin{enumerate}
\item {{\rm [Riesz]} \it If} $s_{n}\in S(V,Y)$ {\it is a basic
sequence, then there exists a function} $f$ {\it from} {\it the
set} $X$ {\it into the Banach space} $Y$ {\it and a null-set} $A$
{\it such that} $s_{n}(x)\rightarrow f(x)$ {\it for all} $ x\in$
$X\backslash A$. \item {{\rm [Egoroff]} \it Moreover, for every}
$\varepsilon>0$ {\it and} $\eta>0$, {\it there exists an index}
$k$ {\it and a countable} {\it family of sets} $A_{t}\in V(t\in
T)$ {\it such that}
\begin{center}
$A  \subset\bigcup_{T}A_{t}$ \quad{\it and} \quad $
\sum_{T}v(A_{t})<\eta$
\end{center}
{\it and for every} $n\geq k$
\begin{center}
$|s_{n}(x)-f(x)|<\varepsilon$  \quad{\it if} \quad $x\not\in
\bigcup_{T}A_{t}$.
\end{center}
\end{enumerate}
\end{lemma}

\bp
    For proof of this lemma see Lemma 2 of Bogdanowicz \cite[page 493]{bogdan10}.
\ep

\bigskip

The idea for the following lemma came from Dunford-Schwartz \cite{DS1},
Part I, p.~111, Lemma~16.

\bigskip

\begin{lemma}[Dunford's Lemma]
\label{Dunford's Lemma} \label{Lemma 3}
\label{Basic sequenceconverging a.e. to 0 converges in seminorm}
Assume that $(X,V,v)$
is a positive measure space on a prering $V$ and $Y$ is a Banach
space. Then the following is true.

If $s_{n}\in S(V,Y)$  is a basic sequence converging almost
everywhere to zero $0$, then the sequence of seminorms $\norm{s_n}$
converges to zero.
\end{lemma}
\bp
    For proof of this lemma see Lemma 3 of Bogdanowicz \cite[pages 493-495]{bogdan10}.
\ep

\bigskip


\bigskip

\section{The Spaces of Lebesgue and Bochner Summable Functions}

\bigskip

\begin{defin}[Lebesgue and Bochner spaces]
Assume that $(X,V,v)$ is a measure space over a prering $V$ of an
abstract space $X$ and $Y$ is a Banach space.

Let $L(v,Y)$ denote the set of all functions $f:X\into Y,$ such
that there exists basic sequence $s_{n}\in S(V,Y)$ that converges
almost everywhere to the function $f.$

The space $L(v,Y)$ is called the space of {\bf Bochner summable}
functions and, for the case when $Y$ is equal to the space $R$ of
reals, $L(v,R)$ represents the space of {\bf Lebesgue summable}
functions.
\end{defin}

Define
\begin{center}
    $\norm{f}=  \lim_n\norm{s_n},\ \int u(f,d\mu)=
    \lim_n\int u(s_{n},d\mu),\ \int fdv=\lim_n\int s_{n}dv$.
\end{center}
Since the difference of two basic sequences is again a basic
sequence, therefore it follows from the Elementary Lemma
\ref{Elementary lemma} and Dunford's Lemma \ref{Lemma 3} that the operators
are well defined, that is, their values do not depend on the
choice of the particular basic sequence convergent to the function
$f$.

\bigskip

\begin{thm}[Basic properties of the space $L(v,Y)$]
\label{Basic properties of the space $L(v,Y)$} \label{Theorem 1}
Assume that $(X,V,v)$ is a positive measure space on a prering $V$
and $Y$ is a Banach space. Then the following is true.

\begin{enumerate}
\item The space $L(v,Y)$ is linear and $\norm{f}$ represents a
seminorm
    being an extension of the seminorm from the space $S(V,Y)$ of simple functions.
\item We have $\norm{f}=0$ if and only if $f(x)=0$ almost everywhere.
\item The functional $\norm{f}$ is $a$
    complete seminorm on $L(v,Y)$ that is given a sequence of functions $f_n\in L(v,Y)$
    such that $\norm{f_n-f_m}\stackrel{nm}\into 0$ there exists a function $f\in L(v,Y)$
    such that $\norm{f_n-f}\stackrel{n}\into 0.$
\item If $f_{1}(x)=f_{2}(x)$ almost everywhere and $f_{2}\in
L(v,Y)$, then $f_{1}\in L(v,Y)$ and
\begin{center}
$\norm{f_{1}}=\norm{f_{2}},\ \int f_{1}dv=\int f_{2}dv,\ \int
u(f_{1},d\mu)=\int u(f_{2},d\mu)$.
\end{center}
\item The operator $  \int fdv$ is linear and represents an
extension onto $L(v,Y)$ of the operator from $S(V,Y)$. It
satisfies the condition $|\int fdv|\leq\norm{f}$ for all $f\in
L(v,Y)$. \item The operator $  \int u(f,d\mu)$ is trilinear on
$U\times L(v,Y)\times  K(v,Z)$ and represents an extension of the
operator from the space $ U\times S(V,Y)\times K(v,Z).$
 It satisfies the condition:
\begin{center}
$|\int u(f,d\mu)|\leq|u|\,\norm{f}\,\norm{\mu}\fa u\in U,\ f\in L(v,Y),\
\mu\in  K(v,Z).$
\end{center}
\end{enumerate}
\end{thm}

\bigskip

\bp
    For proof of this theorem see Theorem 1 of Bogdanowicz \cite[page 495]{bogdan10}.
\ep

\bigskip

\begin{lemma}[Density of simple functions in $L(v,Y)$]
\label{Lemma 4} \label{Density of simple functions in $L(v,Y)$}
Assume that $(X,V,v)$ is a positive measure space on a prering $V$
and $Y$ is a Banach space.
Let  $s_{n}\in S(V,Y)$ be a basic sequence convergent almost
everywhere to a function $f$. Then $\norm{s_{n}-f}\into 0.$
\end{lemma}

\bigskip

\bp
    For proof of this lemma see Lemma 4 of Bogdanowicz \cite[page 495]{bogdan10}.
\ep

\bigskip

From Theorem \ref{Theorem 1} we see that the obtained integrals
are continuous under the convergence with respect to the seminorm
$\norm{\ },$ that is, if $\norm{f_{n}-f}\rightarrow 0$, then
\begin{equation*}
 \int f_{n}dv\rightarrow\int fdv\qtext{and}\int
u(f_{n},d\mu)\rightarrow\int u(f,d\mu).
\end{equation*}
The following theorem characterizes convergence with respect to
this seminorm.

\bigskip

\begin{thm}[Characterization of the seminorm convergence]
\label{Characterization of the seminorm convergence}
\label{Theorem 2} Assume that $(X,V,v)$ is a positive measure
space on a prering $V$ and $Y$ is a Banach space.

Assume that we have a sequence of summable functions $f_{n}\in
L(v,Y)$ and some function $f$ from the set $X$ into the Banach
space $Y.$

Then the following conditions are equivalent
\begin{enumerate}
\item The sequence $f_n$ is Cauchy, that is
$\norm{f_{n}-f_{m}}\stackrel{nm}\into 0,$
 and there exists a subsequence $f_{k_{n}}$  convergent
 almost everywhere to the function $f$.
\item The function $f$ belongs to the space $L(v,Y)$  and
$\norm{f_{n}-f}\rightarrow 0.$
\end{enumerate}
\end{thm}

\bigskip

\bp
    For proof of this theorem see Theorem 2 of Bogdanowicz \cite[page 496]{bogdan10}.
\ep

\bigskip

When the space $Y=R,$ the space $L(v,R)$ represents the space of
Lebesgue summable functions. We have the following relation
between Bochner summable functions and Lebesgue summable
functions.

\bigskip

\begin{thm}[Norm of Bochner summable function is Lebesgue summable]
\label{Norm of Bochner summable function is Lebesgue summable}
\label{Theorem 3} Let $(X,V,v)$ be a positive measure space on a
prering $V$ and assume that $Y$ is a Banach space.

If $f$ belongs  space $L(v,Y)$ of Bochner summable
functions, then the function $|f|$
defined by the formula
\begin{equation*}
 |f|(x)=|f(x)|\fa x\in X
\end{equation*}
 belongs to the space
$L(v,R)$ of Lebesgue summable functions and we have the identity
\begin{center}
    $\norm{f}=\int|f|\,dv\fa f\in L(v,Y).$
\end{center}
\end{thm}

\bigskip

\bp
    For proof of this theorem see Theorem 3 of Bogdanowicz \cite[page 496]{bogdan10}.
\ep

\bigskip

\begin{thm}[Properties of Lebesgue summable functions]
\label{Properties of Lebesgue summable functions} \label{Theorem
4} Let $(X,V,v)$ be a positive measure space on a prering $V$ and
$L(v,R)$ the Lebesgue space of $v$-summable functions.
\begin{description}
    \item[(a)]  If $f\in L(v,R)$  and
        $f(x)\geq 0$ almost everywhere on $X$ then $\int fdv\geq 0$.
    \item[(b)]  If $\ f,g\in L(v,R)$
        and $f(x)\geq g(x)$ almost everywhere on $X$ then $\int f\,dv\geq\int g\,dv.$
    \item[(c)]  If $f,g\in L(v,R)$ and $h(x)=\sup\{f(x),g(x)\}$ for all $x\in X$
        then $h\in L(v,R).$
    \item[(d)]  Let $f_{n}\in L(v,R)$ be a monotone sequence with respect
        to the relation less or equal almost everywhere.
        Then there exists a function $f\in L(v,R)$
        such that $f_{n}(x)\into f(x)$ almost everywhere on $X$
        and $\norm{f_{n}-f}\rightarrow 0$ if and only if the sequence of numbers
        $\int f_{n}\,dv$ is bounded.
    \item[(e)] Let $g,f_{n}\in L(v,R)$ and $f_{n}(x)\leq g(x)$ almost
        everywhere on $X$ for $n=1,2, \ldots$.
        Then the function $h(x)=\sup\{f_{n}(x):\ n=1,2,\ldots\}$
        is well defined almost everywhere on $X$
        and is summable, that is,  $h\in L(v,R)$.
        A function defined almost everywhere is said to be
        summable if it has a summable extension onto the space $X.$
\end{description}
\end{thm}

\bigskip

\bp
    For proof of this theorem see Theorem 4 of Bogdanowicz \cite[page 496-497]{bogdan10}.
\ep

\bigskip

From part (d) of the above theorem we can get the classical
theorem due to Beppo Levi \cite{levi}.

\bigskip

\begin{thm}[Beppo Levi's Monotone Convergence Theorem]
\label{Beppo Levi's Monotone Convergence Theorem}
\label{Monotone Convergence Theorem}
Assume that $(X,V,v)$ is a positive measure
space on a prering $V$ and $L(v,R)$ the Lebesgue space of
$v$-summable functions.

Let $f_{n}\in L(v,R)$ be a monotone sequence with respect to the
relation less or equal almost everywhere. Then there exists a
function $f\in L(v,R)$ such that $f_{n}(x)\into f(x)$ almost
everywhere on $X$ and $$\int f_{n}dv\rightarrow \int f\,dv$$ if
and only if the sequence of numbers $\int f_{n}\,dv$ is bounded.
\end{thm}

\bigskip

Beppo Levi's theorem has  an equivalent formulation in terms of a series.

\bigskip

\begin{thm}[Beppo Levi's theorem for a series]
\label{levi for series}
Assume that $(X,V,v)$ is a positive measure
space on a prering $V$ and $L(v,R)$ the Lebesgue space of
$v$-summable functions.
Let $f_{n}$ be a sequence of nonnegative Lebesgue summable functions.

 Then there exists a
function $f\in L(v,R)$ such that
\begin{equation*}
 \sum_nf_{n}(x)= f(x)\qtext{almost everywhere on} X
\end{equation*}
 and
\begin{equation*}
 \sum_n\int f_{n}\,dv=\int\sum_n f_{n}\,dv= \int f\,dv
\end{equation*}
if
and only if the  sum of the series
\begin{equation*}
 \sum_{n=1}^\infty \int f_{n}\,dv<\infty
\end{equation*}
 is finite.
\end{thm}

\bigskip

\begin{defin}[Linear lattice]
Assume that $L$ represents a linear space of functions from an abstract space
$X$ into the reals $R.$
If the space $L$ is closed under the map $f\mto |f|$ then the operations
$f\vee g$ and $f\wedge g$ given by the formula
    \begin{equation*}
    \begin{split}
    (f\vee g)(x)=&\sup\set{f(x),g(x)}=\frac{1}{2}(f(x)+g(x)+|f(x)-g(x)|)\qtext{and}\\
    (f\wedge g)(x)=&\inf\set{f(x),g(x)}=\frac{1}{2}(f(x)+g(x)-|f(x)-g(x)|)\fa x\in X
    \end{split}
    \end{equation*}
are well defined. Such a space $L$ will be called
a {\bf linear lattice} of functions.

The operation $f\vee g$ is called the {\bf meet} of functions $f,g$
and the operation $f\wedge g$ the joint. Notice that meet operation
is commutative: $f\wedge g= g\wedge f,$ and associative
\begin{equation*}
 (f\wedge g)\wedge h=f\wedge (g\wedge h),
\end{equation*}
and so is the {\bf joint} operation.
\end{defin}

\bigskip

   Notice that from Theorem (\ref{Theorem 3}) follows that
    absolute value of a Lebesgue summable function is itself summable, that is,
    \begin{equation}\label{summability of absolute value}
    |f|\in L(v,R)\fa f\in L(v,R).
    \end{equation}
    Thus the space $L(v,R)$ of Lebesgue summable functions
    forms a linear lattice.

\bigskip

The following theorem is due to P.~Fatou \cite{fatou} and is known in the literature as
Fatou's Lemma.

\bigskip

\begin{thm}[Fatou's Lemma]
Assume that $(X,V,v)$ is a measure space over the prering $V.$
Given a sequence $f_n\in L(v,R)$ consisting of nonnegative Lebesgue summable functions
such that the sequence of integrals
\begin{equation*}
 \int f_ndv
\end{equation*}
is bounded,
then the function $f=\liminf f_n$ belongs to the Lebesgue space $L(v,R)$ and we have
the inequality
\begin{equation*}
 \int f\,dv=\int \liminf f_n\,dv\le \liminf \int f_n dv.
\end{equation*}
\end{thm}

\bigskip

\bp
    For fixed $n$ define the sequence
\begin{equation*}
     g_{nk}=f_n\wedge f_{n+1}\wedge\cdots\wedge f_{n+k}\fa k=0,1,2,\ldots
\end{equation*}
    Since the functions $g_{nk}$ belong to the Lebesgue space $L(v,R)$
    and are nonnegative
    and for fixed $n$ form a decreasing sequence
    $g_{nk}\,(k=1,2,\ldots)$ with respect to relation less or equal $\le\,,$
    from Monotone Convergence Theorem there exists a function $g_n\in L(v,R)$
    such that
\begin{equation*}
     g_n=\lim_k g_{nk}\quad a.e.\qtext{and} \int g_n\,dv=\lim_k\int g_{nk}\,dv.
\end{equation*}

    Notice that
\begin{equation*}
     g_n(x)=\inf\set{f_k(x):\ k\ge n}\quad a.e.\fa n.
\end{equation*}
    Since
\begin{equation}\label{estimate on g sub n}
     g_n\le g_{nk}\le f_{n+k}\quad a.e.\fa n\text{ and }k
\end{equation}
    and the sequence of integrals $\int f_n\, dv$ is bounded,
    the sequence of integrals $\int g_n\, dv$ is bounded.
    Since the sequence $g_n$ monotonically converges almost everywhere
    to the function
\begin{equation*}
     f=\liminf f_n\quad a.e.
\end{equation*}
    by Monotone Convergence Theorem we get that $f\in L(v,R)$
    and
\begin{equation*}
     \lim_n \int g_n\, dv=\int f\, dv.
\end{equation*}

    On the other hand we have for fixed $n$ from the estimate (\ref{estimate on g sub n})
\begin{equation*}
     \int g_n\,dv\le\liminf_k\int f_{n+k}\,dv=\liminf_k\int f_{k}\,dv.
\end{equation*}
    Thus passing to the limit $n\into \infty$ in the above relation
    we get
\begin{equation*}
     \int f\,dv=\int \liminf_k f_k\,dv\le\liminf_k\int f_{k}\,dv.
\end{equation*}
\ep

\bigskip

\begin{thm}[Lebesgue's Dominated Convergence Theorem]
\label{Dominated Convergence Theorem} \label{Theorem 5}
Let
$(X,V,v)$ be a positive measure space on a prering $V$ and $Y$ a
Banach space.

Assume that we are given a sequence $f_{n}\in L(v,Y)$ of Bochner
summable functions that can be majorized by a Lebesgue summable
function $g\in L(v,R),$ that is, for some null set $A\in N$ we have
the estimate
$$|f_{n}(x)|\leq g(x)\fa x\not\in A\text{ and } n=1,2, \ldots$$
Then the condition
   $$f_{n}(x)\rightarrow f(x)\qtext{a.e. on X}$$
implies the relations
    $$f\in L(v,Y)\qtext{and} \norm{f_{n}-f}\rightarrow 0$$
and, therefore, also the relations
\begin{equation*}
 \int f_n\,dv\into \int f\,dv\qtext{and}\int u(f_n,d\mu)\into \int u(f,d\mu)
\end{equation*}
for every bilinear continuous operator $\ u\ $ from the product
$Y\times Z$ into the Banach space $W$ and any vector measure
$\mu\in K(v,Z).$
\end{thm}

\bigskip

\bp
    For proof of this theorem see Theorem 5 of Bogdanowicz \cite[page 497]{bogdan10}.
\ep

\bigskip

\begin{thm}[On absolutely summable series in $L(v,Y)$]
\label{absolutely summable series}
Let
$(X,V,v)$ be a positive measure space on a prering $V$ and $Y$ a
Banach space.
Assume that we are given a sequence $f_{n}\in L(v,Y)$ of Bochner
summable functions.

If $\sum_n\norm{f_n}<\infty$ then there exist a Bochner summable function
$f\in L(v,Y)$ and a Lebesgue summable
function $g\in L(v,R)$ and a null set $A\in N$ such that we have
\begin{equation*}
 \sum_{n=1}^\infty f_{n}(x)= f(x)
\qtext{and}\sum_{n=1}^\infty |f_{n}(x)|\leq g(x)\fa x\not\in A
\end{equation*}
and moreover when $k\to \infty$
we have the relations
\begin{equation*}
\begin{split}
    \norm{\sum_{n\le k}f_{n}-f}     &\into 0,\\
    \int \sum_{n\le k}f_n\,dv       &\into \int f\,dv,\\
    \int u( \sum_{n\le k}f_n,d\mu)  &\into \int u(f,d\mu),
\end{split}
\end{equation*}
or equivalently
\begin{equation*}
         \sum_{n=1}^\infty f_{n}=f,
\end{equation*}
in the sense of convergence in the space $L(v,Y)$ of Bochner summable functions,
and we have the commutativity of the operations
\begin{equation*}
\begin{split}
     \sum_{n=1}^\infty \int f_n\,dv &=\int \sum_{n=1}^\infty f_n\,dv \\
    \sum_{n=1}^\infty \int u( f_n,d\mu)&=\int \sum_{n=1}^\infty  u(f_n,d\mu)=\int u( \sum_{n=1}^\infty f_n,d\mu),
\end{split}
\end{equation*}
for every bilinear continuous operator $\ u\ $ from the product
$Y\times Z$ into the Banach space $W$ and any vector measure
$\mu\in K(v,Z).$
\end{thm}

\bigskip


\bigskip

\section{A characterization of Bochner summable functions}

\bigskip

In this section we shall present several theorems that will be
utilized in the following sections.


\begin{pro}[A characterization of Bochner summable functions]
\label{A characterization of Bochner summable functions} Assume
that $(X,V,v)$ is  a measure space on a prering $V$ of subsets of
an abstract space $X $ and $Y$ a Banach space.

Denote by $G$ the set of non-negative Lebesgue summable functions $g\in L(v,R)$
being a limit of an increasing sequence of simple functions converging
almost everywhere to the function $g.$

A function $f$ mapping $X$ into $Y$ belongs to the space $L(v,Y)$
of Bochner summable functions, if
and only if, there exist a sequence $s_n\in S(V,Y)$ of simple
functions and a function $g\in G,$ such
that $s_n(x)\into f(x)$ almost everywhere on $X$ and
\begin{equation*}
    |s_n(x)|\le g(x)\fa n=1,2,\ldots \qtext{and almost all} x\in X.
\end{equation*}
\end{pro}

\bigskip

\bp
    If $f\in L(v,Y)$ then there exists a basic sequence of the form
    \bb
        s_n=h_1+h_2+\cdots+h_n
    \ee
    converging almost everywhere to the function $f.$
    Notice that the sequence
    \bb
        S_n=|h_1|+|h_2|+\cdots+|h_n|
    \ee
     is nondecreasing and is basic. Thus it converges almost
     everywhere to some summable function $g\in G.$
    Since
    \bb
        |s_n(x)|\le S_n(x)\le g(x)\qtext{for almost all}x\in X,
    \ee
    we get the necessity of the condition.

    The sufficiency of the condition follows from the Dominated Convergence Theorem.
\ep

\bigskip

\begin{thm}[Summability of $u(f,g)$]
\label{Summability of $u(f,g)$} Assume that
$(X,V,v)$ is  a measure space on a prering $V$ of subsets of an
abstract space $X.$ Let $Y,\,Z,\,W$  be Banach spaces and $\ u\ $
a bilinear bounded operator from the product $Y\times Z$ into $W.$

Assume that $f\in L(v,Y)$ and $g\in L(v,Z)$
and either $f$ or $g$ is bounded almost everywhere on $X.$

Then the composed function $u(f,g)=h$ belongs to the space $L(v,W),$
where
\begin{equation*}
 h(x)=u(f(x),g(x))\fa x\in X.
\end{equation*}
\end{thm}

\bigskip
\bp
    Assume that $f_n\in S(V,Y)$ is a basic sequence converging almost
    everywhere to the function $f$ and $g_n\in S(V,Z)$ a basic sequence
    converging to $g.$ Let $\tilde{f}$ and $\tilde{g}$ be the majorants
    of $f$ and $g,$ respectively, from the family $G$ as in the preceding theorem.

    Then the sequence $h_n=u(f_n,g_n)$ represents a sequence of simple
    functions from $S(V,W)$ and it can be majorized by the function
    $m(\tilde{f}+\tilde{g})\in G,$ where either
\begin{equation*}
     |f(x)|\le m\qtext{or}|g(x)|\le m\qtext{for almost all}x\in X.
\end{equation*}
    By continuity of the operator $u$ we get
\begin{equation*}
     h_n(x)\to h(x)\qtext{for almost all}x\in X.
\end{equation*}
    Thus by Theorem \ref{A characterization of Bochner summable functions}
    we get $h\in L(v,W).$
\ep

\bigskip

\begin{pro}[Some product summability]
\label{Summability of a product of functions} Assume that
$(X,V,v)$ is  a measure space on a prering $V$ of subsets of an
abstract space $X $ and $\ Y$ a Banach space. Let
\begin{equation*}
    u(y)=\frac{1}{|y|}\,y\text{ if }|y|>0\qtext{and}u(y)=0\text{ if }|y|=0.
\end{equation*}
Assume that $f\in L(v,Y)$ and $g\in L(v,R).$ Then the product
function $u\circ f\cdot g$ is summable that is $u\circ f\cdot g\in
L(v,Y),$ where $u\circ f$ denotes the composition $(u\circ
f)(x)=u(f(x))$ for all $x\in X.$
\end{pro}

\bigskip

\bp
    Take any natural number $k$ and define a function
    \bb
        u_k(y)=(k|y|\wedge 1)\frac{1}{|y|}y\fa y\in Y,\,|y|>0\qtext{and}u_k(0)=0.
    \ee
    Notice that the function $u_k$ is continuous and
    \bb
        \lim_k u_k(y)=u(y) \qtext{ and } |u_k(y)|\le 1\fa y\in Y.
    \ee
    Let $s_n$ be a basic sequence converging almost everywhere to $f$ and
    $S_n$ a basic sequence converging almost everywhere to $g.$
    Let $S\in G$ be a majorant for the sequence $S_n.$
    Then we have that the sequence
    \bb
        h_{k\,n}=u_k\circ s_n\cdot S_n\in S(V,Y)
    \ee
    consists of simple functions and when $n\into \infty$
    it converges almost everywhere to the function
    \bb
        h_k=u_k\circ f\cdot g.
    \ee
    Since $S$ majorizes the sequence $h_{k\,n},$ from the Dominated Convergence
    Theorem we get $h_k\in L(v,Y)$ and moreover
    \bb
        |h_k(x)|\le S(x)\qtext{for almost all}x\in X.
    \ee
    Passing to the limit $k\into\infty$ and applying the Dominated Convergence
    Theorem yields that $$u\circ f\cdot g\in L(v,Y).$$
\ep

\bigskip

\bigskip

\section{Summable sets}

\bigskip

Assume now again that we have a measure space $(X,V,v)$ on a
prering $V$ of subsets of an abstract space $X.$ Following
Bogdanowicz \cite{bogdan15} and \cite{bogdan20} denote by $V_c$
the family of all sets $A\subset X$ whose characteristic function
$c_A$ is $v$-summable that is $c_A\in L(v,R).$ Put $v_c(A)=\int
c_Adv$ for all sets $A\in V_c.$

\bigskip

\begin{defin}[Summable sets and completion of a measure]
The family $V_c$ will be called the family of {\bf summable sets}
and the set function $v_c$ from $V_c$ into $R$ will be called the {\bf completion}
of the measure $v.$
\end{defin}

\bigskip

From Theorem
\ref{Properties of Lebesgue summable functions}
concerning properties of Lebesgue summable functions
we can deduce the following proposition.


\bigskip

\begin{pro}[Summable sets form a delta ring]
Assume that $(X,V,v)$ is  a measure space on a prering $V$ of
subsets of an abstract space $X.$

Then the family $V_c$ of summable sets forms a $\delta$-ring and
$v_c$ forms a measure. If in addition $X\in V_c$ then $V_c$ forms
a $\sigma$-algebra.
\end{pro}


\bp
    The space $L(v,R)$ forms a linear lattice.
    Thus from the identities
    \begin{equation*}
    c_{A\cup B}=c_A\vee c_B\qtext{and}c_{A\cap B}=c_A\wedge c_B\qtext{and}c_{A\less B}=
    c_A-c_A\wedge c_B
    \end{equation*}
    we can conclude that the family $V_c$ of summable sets forms a ring.

    Now if $A_n\in V_c$ is a sequence of summable sets
    and $B_n=\bigcap_{j\le n} A_j$ and $B=\bigcap_{j\ge 1} A_j,$
    from the Dominated Convergence Theorem \ref{Dominated Convergence Theorem}
    and from the relations
    \begin{equation*}
    |c_{B_n}(x)|=c_{B_n}(x)\le c_{A_1}(x)\qtext{and}c_{B_n}(x)\into c_B(x)\fa x\in X
    \end{equation*}
    we get that $B\in V_c.$ Thus the family $V_c$ of summable sets forms a
    $\delta$-ring.

    In the case, when $X\in V,$ we get from the de Morgan law
    and the fact that $V_c$ forms a $\delta$-ring that
    \begin{equation*}
    \bigcup_{n\ge 1}A_n=X\less \bigcap_{n\ge 1}(X\less A_n)\in V_c.
    \end{equation*}
    Hence in this case $V_c$ forms a $\sigma$-algebra.

    To show that the triple $(X,V_c,v_c)$ forms a positive
    measure space assume that $A\in V_c$
    and a sequence of disjoint sets $A_n\in V_c$ forms a decomposition of the set $A.$
    So $A=\bigcup_{j\ge 1} A_j.$ Let $B_n=\bigcup_{j\le n} A_j.$
    From the Dominated Convergence Theorem \ref{Dominated Convergence Theorem}
    and from the relations
    \begin{equation*}
    |c_{B_n}(x)|=c_{B_n}(x)\le c_{A}(x)\qtext{and}c_{B_n}(x)\into c_A(x)\fa x\in X
    \end{equation*}
    and linearity of the integral,
    we get that
    \begin{equation*}
    \begin{split}
    v_c(A)&=\lim_n v_c(B_n)=\lim_n \int c_{B_n}dv\\&=\lim_n \sum_{j\le n}\int c_{A_j}dv=
    \lim_n \sum_{j\le n}v_c(A_j)=\sum_{j}v_c(A_j).
    \end{split}
    \end{equation*}
     Thus $v_c$ is countably additive on the delta ring $V_c.$
\ep

\bigskip



\bigskip

\section{Summability on sets}

\bigskip


Assume that $(X,V,v)$ is  a measure space on a prering $V$ of
subsets of an abstract space $X $ and $Y$ a Banach space.

A vector measure $\mu$ from a prering $V$ into a Banach space $Y$
is said to be of {\bf finite variation} on $V$
if 
\bb
    |\mu|(A)=\sup\set{\sum_{t\in T}|\mu(A_t)|:\ A=\bigcup_{t\in T}A_t}<\infty\fa A\in V
\ee
where the supremum is taken over all finite disjoint
decompositions $A_t\in V\,(t\in T)$ of the set
$A=\bigcup_{t\in T}A_t.$ The set function $|\mu|$ is called the {\bf variation} of
the vector measure $\mu.$

We shall say that a function $f:X\into Y$ is {\bf summable on a
set} $A\subset X$ if the product function $c_Af$  belongs to the space
$L(v,Y)$ of Bochner summable functions and
we shall write
$$\int_Af\,dv=\int c_Af\,dv$$
to denote the integral of the function $f$ on the set $A.$

Denote by $V_f$ the family of all sets $A\subset X$ on which the function
$f$  is summable.
Notice that the family $V_c$ of summable sets can be thought of as
family of sets on which the characteristic function $c_X$ of the entire space $X$
is summable.

\bigskip


\begin{pro}[Sets on which a function is summable form a $\delta$-ring]
Assume that $f$ is a function from the space $X$ into the Banach space $Y.$

Then the family of sets $V_f$ forms a
$\delta$-ring and the set function
\begin{equation*}
 \mu(A)=\int_Af\,dv\fa A\in V_f
\end{equation*}
forms a $\sigma$-additive
vector measure of finite variation on $V_f.$
\end{pro}

\bigskip

\bp
    Assume that $A,B\in V_f.$ Then $c_Af\in L(v,Y)$ and $|c_Bf|\in L(v,R).$
    From Proposition \ref{Summability of a product of functions} we get
    \bb
        c_{A\cap B}f=c_Ac_B\,u\circ{f}\cdot|f|=u\circ (c_Af)\cdot |c_Bf|\in L(v,Y).
    \ee
    Thus $A\cap B\in V_f.$ It follow from linearity of the space $L(v,Y)$
    and the identities
    \bb
        c_{A\less B}=c_A-c_{A\cap B}\qtext{and}
        c_{A\cup B}=c_{A\less B}+c_{A\cap B}+c_{B\less A}
    \ee
    that $V_f$ forms a ring.

    Now using the Dominated Convergence Theorem we can easily prove that
    $V_f$ forms a $\delta$-ring and the set function
    \bb
        \mu(A)=\int_Af\,dv\fa A\in V_f
    \ee
    forms a $\sigma$-additive vector measure and
    \bb
        |\mu|(A)\le \int_A|f|\,dv<\infty\fa A\in V_f.
    \ee
\ep

\bigskip

The family of sets $V_f$ on which a function $f$ is summable may
consist only of the empty set. However in the case of a summable
function this family is rich as follows from the following
corollary.

\bigskip


\begin{cor}[Sets on which a summable function is summable form $\sigma$-algebra]
Assume that $(X,V,v)$ is  a measure space on a prering $V$ of
subsets of an abstract space $X $ and $Y$ is a Banach space.

If $f\in L(v,Y)$ is a summable function then the family $V_f$ of
sets, on which $f$ is summable, forms a $\sigma$-algebra
containing all summable sets, that is, we have the inclusion
$$V\subset V_c\subset V_f$$
and the set function $\mu(A)=\int_Af\,dv$ is $\sigma$-additive of
finite total variation $|\mu|(X)\le \int|f|\,dv.$
\end{cor}

\bigskip

\bp
    To prove this corollary notice that similarly as before we can
    prove that product $gf$ of a summable bounded function $g\in L(v,R)$
    with a summable function $f\in L(v,Y)$ is summable $gf\in L(v,Y).$
    This implies that $V\subset V_c\subset V_f.$
\ep

\bigskip

For further studies of vector measures we recommend Dunford and
Schwartz \cite{DS1}, and for extensive survey of the state of the
art in the theory of vector measures see the monograph of Diestel
and Uhl \cite{diestel-uhl}.
Compare also Bogdanowicz and Oberle \cite{bogdan55}.

\bigskip



\section{Measures generating the same integration}

\bigskip

In view of the existence of a variety of extensions of a measure
from a prering onto $\delta$-rings and the multiplicity of
extensions to Lebesgue measures it is important to be able to
identify measures that generate the same class of Lebesgue-Bochner
summable functions $L(v,Y)$ and the same trilinear integral $\int
u(f,\,d\mu)$ and thus the ordinary Bochner integral $\int f\,dv.$
In this regard we have the following theorems.

Assume that $(X,V_j,v_j),\,(j=1,2)$ are two measure spaces over
the same abstract space $X$ and $Y,Z,W$ are any Banach spaces and
$U$ is the Banach space of bilinear bounded operators from the
product $Y\times Z$ into $W.$

\bigskip

\begin{thm}[When $L(v_2,Y)$ extends $L(v_1,Y)$?]
\label{When $L(v_2,Y)$ extends $L(v_1,Y)$?} For every Banach space
$Y$ we have $L(v_1,Y)\subset L(v_2,Y)$ and \bb
    \int f\,dv_1=\int f\,dv_2\fa f\in L(v_1,Y)
\ee if and only if $V_{1c}\subset V_{2c}$ and \bb
    v_{1c}(A)=v_{2c}(A)\fa A\in V_{1c}
\ee that is the measure $v_{2c}$ represents an extension of the
measure $v_{1c}.$
\end{thm}

\bigskip

Consequently we have the following theorem.

\begin{thm}
For any Banach space $Y$ and any bilinear bounded transformation
$u\in U$ we have $L(v_1,Y)= L(v_2,Y)$ and \bb
    \int f\,dv_1=\int f\,dv_2\fa f\in L(v_1,Y)
\ee and the spaces of vector measures
$K(v_1,Z),K(v_2,Z),K(v_{1c},Z),K(v_{2c},Z)$ are isometric and
isomorphic and \bb
    \int u(f,d\mu_1)=\int u(f,d\mu_2)=\int u(f,d\mu_{1c})
    =\int u(f,d\mu_{2c})\fa f\in L(v_1,Y),
\ee where $\mu_1,\mu_2,\mu_{1c},\mu_{2c}$ are vector measures that
correspond to each other through the isomorphism, if and only if,
the completions of the measures $v_1,v_2$ coincide
$v_{1c}=v_{2c}.$
\end{thm}

For proofs of the above theorems see Bogdanowicz \cite{bogdan20}.
It is important to relate the above theorems to the classical
spaces of Lebesgue and Bochner summable functions and the
integrals generated by Lebesgue measures.

Since there is a great
variety of approaches to construct these spaces we shall
understand by a classical construction of the Lebesgue space
$L(\mu,R)$ the construction developed in Halmos \cite{halmos} and
by classical approach to the theory of the space $L(\mu,Y)$ of
Bochner summable functions as presented in Dunford and Schwartz
\cite{DS1}.

Now if $(X,V,v)$ is a measure space on a prering $V$ and
$(X,M,\mu)$ represents a Lebesgue measure space where $\mu$ is the
smallest extension of the measure $v$ to a Lebesgue complete
measure on the $\sigma$-ring $M,$ then we have the following
theorem.

\begin{thm}
For every Banach space $Y$ the spaces $L(v,Y)$ and $L(\mu,Y)$
coincide and we have \bb
    \int_A f\,dv=\int_A f\,d\mu \fa f\in L(\mu,Y) \qtext{and} A\in M.
\ee
\end{thm}

The above theorem is a consequence of the theorems developed in
Bogdanowicz~\cite{bogdan14}, page~267, Section~7.

\bigskip



\section{Continuity and summability on locally compact spaces}

\bigskip

Now let the measure space $(X,V,v)$ be the {\em Riemann measures space}
over the Euclidean space $R^{m}=X$ on the prering $V$ of all cubes of the form
\begin{equation*}
  A=J_{1}\times\ldots\times J_{m},
\end{equation*}
where $J_{i}$ denotes an interval with end points
 $a_{i}\leq b_{i}$

or let $X$ be a locally compact Hausdorff
space, the family $V$ consist of all sets $A=F\cap G$, where $F$
is a compact set and $G$ is an open set, and let $v$ be any
positive measure on $V$.

\begin{thm}[Continuous functions on compact sets are summable]
\label{Continuous functions on compact sets are summable}
\label{Theorem 8}
Assume that the triple $(X,V,v)$ represents either the Riemann
measure space or the measure space as defined above over
a locally compact Hausdorff space $X.$

If a function $f$ is continuous
from a compact set $Q\subset X$ into the Banach space $Y$, then the function
is Bochner summable on the set $Q$ that is $c_Qf\in L(v,Y)$ and thus
the integrals
\begin{equation*}
 \int_Qf\,dv\qtext{and}\int_Q u(f,d\mu)=\int u(c_Qf,d\mu)
\end{equation*}
exist and moreover we have the estimates
\begin{equation*}
 \norm{\int_Qf\,dv}\le\int_Q\norm{f}\,dv
 \qtext{and}\norm{\int_Q u(f,d\mu)}\le |u|\norm{\mu}\int_Q\norm{f}\,dv
\end{equation*}
for all $\mu\in K(v,Z),$ and $u\in U.$
\end{thm}

\bigskip

\bp
    Consider first the case of the Euclidean space $R^m.$
    In this case there exists a cube $A$
    such that $Q\subset A$. There exist disjoint cubes of diameter less than
    $1/n,$  $A_{i}^{n}(i=1, \ldots , k_{n}),$
    such that the intersections $A_{i}^{n}\cap Q$  are nonempty
    and $Q\subset\bigcup_iA^{n}_i\subset A$.
    Pick a point $x_{i}^{n}$ from each of the sets $A_{i}^{n}\cap Q$.

    In the case
    of a locally compact space consider the compact set $f(Q)=K\subset Y$.
    There exist nonempty disjoint sets $B_{i}^{n}(i=1, \ldots,k_{n})$
    of diameter less than $1/n$ such that $\bigcup_iB^{n}_i=K$ and each of
    the sets $B_{i}^{n}$ is the intersection of an open set
    with a closed set.
    Therefore,
    $A_{i}^{n} =f^{-1}(B_{i}^{n})\bigcap Q\in V$.
    Choose a point $x_{i}^{n}$ in each of the sets $A_{i}^{n}$.

    Define $f_{n}(x) =  \sum_{i}f(x_{i}^{n})c_{A_{i}^{n}}$. We
    have $f_{n}(x)\rightarrow c_{Q}(x)f(x)$ for every $x\in X,$ and the sequence
    $|f_{n}(x)|$ is dominated either by $Mc_{A}$ or by
    $Mc_{Q}$.
    Notice that $f_n\in S(V,Y)$ by construction of the sequence.
    Using the Dominated Convergence Theorem \ref{Theorem 5} we conclude the proof.
\ep

\bigskip

In the case when the space $X$ is topological, we say that the {\bf measure is regular} if
$$v(A)=  \inf\{v(E):\ A\subset \mathrm{int}(E),\ E\in V\}=
\sup\{v(E):\ \mathrm{clo}(E)\subset A,\ E\in V\},$$ where
$\mathrm{int}(E)$ denotes the interior and $\mathrm{clo}(E)$ the
closure of the set $E.$ The measure defined in the Euclidean space $R^m$
by the formula $$v(J_{1}\times\ldots\times
J_{m})=(b_{1}-a_{1})\cdots(b_{m}-a_{m})$$ is regular.

Denote by $C(Y)$ the family of
all functions $f$ from the set $X$ into the space $Y$
such that there exists an increasing
sequence of compact sets $Q_{n}$ with the following property
$$\{x\in X:\ f(x)\neq 0\}\less \bigcup_n Q_{n}\in N,$$
the function $f$ restricted to each of the
sets $Q_{n}$ is continuous on  $Q_{n},$
and $$\sup\int_{Q_n}|f|\,dv<\infty.$$
Such functions will be called {\bf almost $\sigma$-continuous}
(sigma continuous.)

\bigskip

The following theorem is equivalent to a theorem of Lusin~\cite{lusin}.
One can prove it using the fact that continuous functions on compact sets are summable
\ref{Continuous functions on compact sets are summable},
and the Dominated Convergence Theorem
 \ref{Dominated Convergence Theorem},
 and Riesz-Egoroff property of a basic sequence
 \ref{Riesz-Egoroff property of basic sequence}.

\bigskip

\begin{thm}[Lusin's Theorem]
\label{Lusin's Theorem}
\label{Theorem 9}  If $(X,V,v)$ is the Riemann measure space
over $X=R^n,$ that is, the prering $V$ consists of cubes
as above, and the measure $v$ represents the volume of the $n$-dimensional
cube, then the space of Bochner summable functions
coincides with the space of almost $\sigma$-continuous functions, that is,
we have the identity $$C(Y)=L(v,Y).$$
\end{thm}

\bigskip

To prove the following generalization of Lusin's theorem one may use in addition
to previously mentioned properties the fact that in a locally compact topological Hausdorff
space one can separate any two disjoint closed sets by means of a continuous function.
This fact is known as Urysohn's Theorem, see Yosida \cite{yosida}, page~7.

For a complete proof of the following theorem see
Bogdanowicz~\cite{bogdan23}, page~230, Theorem~9, Part~3,
and the definition of the family $C_\sigma(Y)$ on page~221, preceding Theorem~1.

\begin{thm}[Generalized Lusin's Theorem]
\label{Generalized Lusin's Theorem}
If $X$ is a locally compact Hausdorff topological space and
prering $V$ consists of sets of the form $Q\cap G,$ where $Q$ is compact
and $G$ is an open set, then
every almost $\sigma$-continuous function is Bochner summable $C(Y)\subset L(v,Y)$.

If  measure is regular, then the space of Bochner summable functions
coincides with the space of almost $\sigma$-continuous functions $$C(Y)=L(v,Y).$$
\end{thm}


\bigskip


\bigskip

\section{Extensions to Lebesgue measures}

\bigskip

If $V$ is any nonempty collection of subsets of an abstract space
$X$ denote by $V^{\sigma}$ the collection of sets that are
countable unions of sets from $V$ and denote by $V^{r}$ the
collection \bb
    V^r=\set{A\subset X:\ A\cap B\in V\fa B\in V}.
\ee Now assume that $(X,V,v)$ is a measure space and $V_c$ the
$\delta$-ring of summable sets and $v_c(A)=\int c_Adv.$ It is easy
to prove that $V_c^\sigma$ forms the smallest $\sigma$-ring
containing the prering $V$ and the family $N$ of $v$-null sets
(see Bogdanowicz \cite{bogdan14}, section 3, page 258).

\bigskip

Moreover the set function defined by
\begin{equation}\label{completing v_c to Lebesgue measure}
    \mu(A)=\sup\set{v_c(B):\ B\subset A,\,B\in V_c}\fa A\in V_c^\sigma
\end{equation}
forms a Lebesgue measure on $V_c^\sigma.$

\bigskip

A Lebesgue measure is called {\bf complete} if all subsets of sets
of measure zero are in the domain of the measure and thus have
measure zero. The above Lebesgue measure $\mu$ can be
characterized as the smallest extension of the measure $v$ to a
complete Lebesgue measure.
This fact follows from Part 7, of Theorem 4, of \cite{bogdan14}, page 259.
Hence this measure is unique.

\bigskip

For any fixed function $f$ from $X$ into a Banach space $Y$ denote by
$V_f$ the collection of all sets $A\subset X$ such that $c_Af\in L(v,Y).$
The family $V_f$ may consist of just the empty set but in the case of
a summable function $f\in L(v,Y)$ it represents a sigma algebra extending
the family $V_c$ of summable sets.

\bigskip

One can prove the following identity
\bb
    V_c^r=\bigcap_{f\in L(v,Y)}V_f=\bigcap_{f\in L(v,R)}V_f.
\ee
The family $V_c^r$ as an intersection of $\sigma$-algebras
forms itself a $\sigma$-algebra containing the $\delta$-ring
$V_c.$ The smallest $\sigma$-algebra $V^a$ containing $V_c$ is
given by the formula
\bb
    V^a=\set{A\subset X:\ A\in V_c^\sigma\text{ or }X\less A\in V_c^\sigma}.
\ee

If $X\in V_c^\sigma$ then the sigma algebras coincide
$V_c^\sigma=V^a=V_c^r.$ If  $X\not\in V_c^\sigma,$ one can always
extend the measure $v$ to a Lebesgue measure on  $V^a$ or $V_c^r$
by the formula
\begin{equation}\label{trivial extension}
    \mu(A)=v_c(A)\text{ if } A\in V_c
    \qtext{and}\mu(A)=\infty\text{ if }A\not\in V_c.
\end{equation}
However if $\sup\set{v_c(A):\ A\in V_c}=a<\infty$ and $X\not\in
V_c^\sigma$ the extensions are not unique. Indeed one can take in
this case $\mu(X)=b,$ where $b$ is any number from the interval
$[a,\infty),$ and put \bb
    \mu(A)=v_c(A)\text{ if }A\in V_c\qtext{and}
    \mu(A)=b-v_c(X\less A)\text{ if }X\less A\in V_c
\ee to extend  the measure $v_c$ onto the $\sigma$-algebra $V^a$
preserving sigma additivity.

Consider an example. Let $(X,V,v)$ be the following measure space:
\bb
    X=R,\quad V=\set{\emptyset,\{n\}:\ n=1,2,\ldots},\quad v(\emptyset)=0,\,v(\{n\})=2^{-n}.
\ee In this case the family $N$ of null sets contains only the
empty set $\emptyset,$ the family $S$ of simple sets consists of
finite subsets of the set of natural numbers $\N,$ the family
$V_c$ of summable sets consists of all subsets of $\N,$ we have
$V_c^\sigma=V_c,$ the smallest $\sigma$-algebra extending $V_c$
consists of sets that either are subsets of $\N$ or their
complements are subsets of $\N,$ finally $V_c^r=P(R)$ consists of
all subsets of $R.$ Since $v_c(\N)=1$ is the supremum of $v_c,$
the measure $v_c$ has many extensions onto the $\sigma$-algebras
$V^a$ and $V_c^\sigma.$ An infinite valued extension onto the power set
$P(R)$ is given by the formula
(\ref{trivial extension}) and a totaly finite valued extension, for instance, by \bb
    \mu(A)=\sum_{n\in A\cap \N}2^{-n}\fa A\subset R.
\ee

\bigskip


\bigskip

\section{Extension to outer measure}

\bigskip

In the proofs of some theorems it is useful to use
the notion of the outer measure also called the exterior measure in the literature.

We shall say that a nonempty family $V$ of subsets
of a space $X$ is {\bf hereditary} if it satisfies the implication
\begin{equation*}
 A\subset B\in V\impl A\in V.
\end{equation*}

\bigskip

\begin{defin}[Outer measure]
A set function $\eta:V\mto [0,\infty]$ is called an {\bf outer measure}
if its domain $V$ forms a hereditary sigma ring of subsets of an abstract space $X$
and it takes value zero on the empty set, is monotone, and countably subadditive,
that is it has the following properties
\begin{enumerate}
\item $\eta(\emptyset)=0,$
\item if $A\subset B\in V$ then $\eta(A)\le \eta(B),$
\item for any sequence $A_n\in V$ we have
$
 \eta(\bigcup_n A_n)\le \sum_n\eta(A_n).
$
\end{enumerate}
\end{defin}

\bigskip

\begin{defin}[Set function $v^*$]
Assume that $(X,V,v)$ is a positive measure space over the prering $V.$

For any subset $A$ of the space $X$ define the set function $v^*$
as follows if there exists a summable set $B$ containing the set $A$
\begin{equation*}
 v^*(A)=\inf\set{v_c(B):\ A\subset B\in V_c}
\end{equation*}
else
let $v^*(A)=\infty.$
\end{defin}

\bigskip

\begin{thm}[Function $v^*$ forms an outer measure extending $v_c$]
Assume that  $(X,V,v)$ forms a positive measure space over the prering $V$
and $v^*$ represents the set function as defined above.

Then
\begin{enumerate}
\item if for some set $v^*(A)<\infty$ then there exists a summable
set $B\in V_c$ such that $A\subset B$ and $v^*(A)=v_c(B),$

\item the set function $v^*$ forms an outer measure,

\item it extends the measure $v_c$ from the
family $V_c$ of the summable sets onto the family of all subsets of the space $X,$
and therefore it extends also the measure $v.$
\end{enumerate}

\end{thm}

\bigskip

\bp
    To prove part (1) take any positive integer $n.$ Since the number
    $$v^*(A)+1/n$$ is greater than the greatest lower bound $v^*(A)$ of the
    set of numbers
    \begin{equation*}
        D=\set{v_c(B):\ A\subset B\in V_c},
    \end{equation*}
    the number $v^*(A)+1/n$ is not a lower bound of the number set $D.$
    Therefore, for each $n,$ there exists a summable set $B_n\in V_c,$ such that
    \begin{equation*}
        A\subset B_n\qtext{and}v^*(A)\le v_c(B_n)<v^*(A)+1/n.
    \end{equation*}
    Put $B=\bigcap_n B_n.$ We have $A\subset B\in V_c$ and from monotonicity
    of the measure $v_c$ we get
    \begin{equation*}
        v^*(A)\le v_c(B)\le v_c(B_n)<v^*(A)+1/n\fa n=1,2,\ldots
    \end{equation*}
    Hence $v^*(A)=v_c(B).$

    As a consequence if $A\in V_c$ then $v^*(A)=v_c(A)$ thus the set function
    $v^*$ extends the measure $v_c$ from the family $V_c$ of summable sets to
    all subsets of the space $X.$

    Notice that since the empty set $\emptyset$ is in $V_c$
    we must have $v^*(\emptyset)=v_c(\emptyset)=0.$

    To prove the monotonicity of the set function $v^*,$ that is, the validity
    of the implication
\begin{equation*}
     A\subset B\impl v^*(A)\le v^*(B),
\end{equation*}
    notice that the above implication is always true when $v^*(B)=\infty.$
    To establish the validity of the implication when $v^*(B)<\infty$
     apply part (1) of the theorem.

    Finally to prove countable subadditivity
\begin{equation*}
     v^*(\bigcup_n A_n)\le \sum_n v^*(A_n),
\end{equation*}
    notice that the inequality is always true when the right side is infinite.
    So consider the case when it is finite. In this case each term of the sum
    must be finite and so there exist summable sets $B_n\in V_c$
    such that
\begin{equation*}
     A_n\subset B_n\qtext{and}v^*(A_n)=v_c(B_n)\fa n=1,2,\ldots
\end{equation*}
    Thus we have
\begin{equation*}
     A=\bigcup_n A_n\subset \bigcup_n B_n=B.
\end{equation*}
    Consider the set $B.$ Introduce the sets
\begin{equation*}
     D_1=B_1,\quad D_n=B_n\less (B_1\cup\cdots\cup B_{n-1})\fa n>2
\end{equation*}
    Notice that the sets $D_n\in V_c$ and they are disjoint. Moreover
    $D_n\subset B_n$ and for all $m$ we have the estimate
\begin{equation*}
     \sum_{n\le m}v_c(D_n)\le \sum_{n\le m}v_c(B_n)\le \sum_{n=1}^\infty v_c(B_n)<\infty.
\end{equation*}
    Since $B=\bigcup_n D_n$ the set $B$ is in the ring $V_c$ and
    from countable additivity of the measure $v_c$ we have
\begin{equation*}
     v^*(A)\le v_c(B)=\sum_n v_c(D_n)\le\sum_n v_c(B_n)=\sum_n v^*(A_n).
\end{equation*}
    The above inequality completes the proof.
\ep

\bigskip

The outer measure $v^*$ has other convenient representations and in the case of
topological spaces it can be defined by or related to  topology.
To do that we will need some theorems on approximation of sets.

We shall use the following notation.
If $V$ is any nonempty family of sets then $V^\sigma$
will denote the family of all sets that are
countable unions of sets from the
family $V,$ and by $V^\delta$ the family of all sets that are
countable intersections of sets from the
family $V.$ We shall use the following shortcut notation
\begin{equation*}
 V^{\sigma\delta}=(V^\sigma)^\delta.
\end{equation*}
We remind the reader that if $V$ is any collection of sets containing at least the
empty set $\emptyset$ then $S(V)$ denotes the collection of simple sets
consisting of finite disjoint unions of sets from the collection $V.$

The operation $\div$ will denote the symmetric difference of sets
\begin{equation*}
 A\div B=(A\less B)\cup(B\less A).
\end{equation*} The collection of all subsets of a space $X$ with this operation
forms an Abelian group.

\bigskip

The following theorem is a consequence of
Theorem~4, Part~6, page~259, of Bogdanowicz~\cite{bogdan14}.
We shall present its proof for the sake of completeness, since
it is essential in the sequel development of the theory.

\bigskip

\begin{thm}[A characterization of summable sets]
\label{sigma-delta sets}
Assume that $(X,V,v)$ is a positive measure space. Assume that the prering $V$
contains the space $X.$ Let $S=S(V)$ and $N$ denote the collection of null sets.

Then a set $A\subset X$ is summable, that is $A\in V_c,$
if and only if one of the following conditions is satisfied
\begin{enumerate}

\item There is a set $B\in S^{\delta\sigma}$
and  $D\in N$ such that $A=B\div D.$
\item There is a set $B\in V^{\sigma\delta}$
and  $D\in N$ such that $A=B\div D.$
\end{enumerate}
\end{thm}

\bigskip

\bp
    Assume that the set $A$ is summable.
Denote by $S$ the collection
$S(V)$ of simple sets. It follows
from the definition of a summable function
that there exists a sequence of simple functions $s_n \in S(V,R)$
and a null set $D\in N$  such that $s_n(x)\to c_A(x)$ if $x \not\in D.$
Put
\begin{equation*}
 f_n(x)= \inf\set{s_m(x):\ m > n}.
\end{equation*}
The sequence of values $f_n(x)$ increasingly converges to the value $c_A(x)$ for
every point $x \not\in D.$

Denote by $B$ the set of all points $x \in X$ for which there exists an index $n$
such that $f_n (x) > 0.$ We see that $B = \bigcup_{np} B_{np}$ where
$B_{np} = \set{x \in X :\  f_n(x) \ge 1/p}$
and further $B_{np} = \bigcap_{\,m\ge p} E_{mp}$
where $E_{mp} = \set{x \in X :\  s_m(x) \ge 1/p}.$ Since $E_{mp}\in S$
therefore $B \in S^{\delta\sigma}.$

We notice that $A \div B \subset D.$ Hence the set $C = A \div B$ as a subset of the
null set D is itself a null set. From the properties of the symmetric difference we
get $A = (B\div B ) \div A = B\div (B\div A) = B \div C.$ Thus the necessity of the
condition is proved.

To prove the sufficiency of the condition, take any set $B \in S^{\delta\sigma}.$
The set $B$ obviously belongs to the sigma algebra $V_c.$
Let $A = B \div C$ where $C \in N.$ The characteristic functions
$c_A$ and $c_B$ are equal almost everywhere.

Since $c_B\in V_c$ therefore
there exists a basic sequence of simple functions
$s_n\in S (V,R)$ converging almost everywhere to the function $c_B$ and therefore also
converging almost everywhere to the function $c_A.$ Thus we have $c_A \in L(v,R),$
that is, $A \in V_c.$ Thus the condition (1) is equivalent to summability of the
set $A.$
\bigskip

Now assume that the condition (1) is satisfied. Thus
\begin{equation*}
 A=B\div D,\quad B\in S^{\delta\sigma},\quad D\in N.
\end{equation*}

Since $V_c$ forms a sigma algebra we have $A_1=X\div B\in V_c.$
Thus $A_1$ must satisfy the condition~(1).
Therefore there exists sets $B_1$ and $D_1$
such that
\begin{equation*}
 A_1=B_1\div D_1,\quad B_1\in S^{\delta\sigma},\quad D_1\in N.
\end{equation*}

We have
\begin{equation*}
\begin{split}
 A&=X\div X\div B\div D\\
 &=X\div (X\div B)\div D\\
 &=X\div A_1\div D\\
 &=X\div B_1\div D_1\div D.
\end{split}
\end{equation*}
Since $X\div B_1=X\less B_1,$ it follows from De Morgan law on complements of sets
that the set $B_2$ satisfies the relations
\begin{equation*}
 B_2=X\less B_1\in S^{\sigma\delta}=V^{\sigma\delta}.
\end{equation*}
Since the set $D_2=D_1\div D\subset D_1\cup D$ as a subset of a null set
it is itself a null set. Thus we have
the required representation $A=B_2\div D_2.$ Thus from condition (1)
we derived the condition (2). It is easy to see from the symmetry
of the argument that assuming that condition (2) holds for any
summable set then the condition (1) must be also true.
\ep

\bigskip

\begin{pro}[For a prering $V$ family $V^\sigma$ is closed under finite intersection]
\label{pro-V sigma closed under intersection}
Let $X$ be an abstract space and $V$ a prering of subsets of $X.$
If $A,B\in V^\sigma$ then $A\cap B\in V^\sigma.$
\end{pro}

\bigskip
\bp
    Assume $A=\cup_m A_m$ and $B=\cup_n B_n$ where the sequences $A_m$ and $B_n$
    are from the prering $V.$ Then
\begin{equation*}
     A\cap B=\bigcup_{m,n}(A_m\cap B_n).
\end{equation*}
    By definition of a prering each set $A_m\cap B_n$ is a union of a finite number
    of sets from the family $V.$ Thus $A\cap B\in V^\sigma.$
\ep

\bigskip

\begin{pro}
\label{pro-approx by V sigma}
Assume that $(X,V,v)$ is a positive measure space.
Assume that the prering $V$ forms a pre-algebra,
that is $X\in V.$

Then for every summable $A\in V_c$ and every $\varepsilon>0$ there exists
a set $B\in V^\sigma$ such that
\begin{equation*}
 A\subset B\qtext{and}v_c(B)\le v_c(A)+\varepsilon.
\end{equation*}
\end{pro}

\bigskip

\bp
    Take any summable set $A$ and any $\varepsilon>0.$ By Theorem \ref{sigma-delta sets}
    there exists a sequence $B_n\in V^\sigma$ such that $A=B\div D$
    where $B=\bigcap_n B_n$ and $D$ is null set.

    By Proposition \ref{pro-V sigma closed under intersection}
    the family of sets $V^\sigma$ is closed under finite intersections.
    Thus we have
\begin{equation*}
     B'_n=\bigcap_{j\le n}B_j\in V^\sigma\fa n.
\end{equation*}
    The sequence of characteristic functions $c_{B'_n}$ of sets $B'_n$ decreasingly converges
    to the characteristic function $c_B$ and therefore it converges
    almost everywhere to the characteristic function $c_A.$
    The sequence is bounded by summable
    function, namely by the characteristic function $c_X.$ By Lebesgue's
    Dominated Convergence Theorem
\begin{equation*}
     v_c(B'_n)\to v_c(A).
\end{equation*}

    Since $D$ is a null set, there exists a sequence $C_j\in V$ such that
\begin{equation*}
     \sum_{j=1}^\infty v(C_j)<\infty\qtext{and}D\subset D_n=\bigcup_{j>n}C_j\fa n.
\end{equation*}
    The sequence $D_n$ is decreasing and since
\begin{equation*}
     v_c(D)\le \sum_{j>n}v(C_j)\to 0.
\end{equation*}

    Thus we have
\begin{equation*}
     A\subset B'_n\cup D_n\in V^\sigma\fa n,
\end{equation*}
    and
\begin{equation*}
     \quad v_c(A)\le v_c(B'_n\cup D_n)\le v_c(B'_n)+v_c(D_n)\fa n.
\end{equation*}
    Thus for sufficiently large index $n$ we will have
\begin{equation*}
     \quad v_c(A)\le v_c(B'_n\cup D_n)\le v_c(A)+\varepsilon.
\end{equation*}
\ep

\bigskip

\begin{pro}
\label{pro-decompose into disjoint sets}
Let $X$ be an abstract space and $V$ a prering of subsets of $X.$
If $A\in V^\sigma$ then there exists a sequence $I_n$ of disjoint sets
from the prering $V$ such that
\begin{equation*}
 A=\bigcup_n I_n.
\end{equation*}
\end{pro}

\bigskip

\bp
    By definition of the family $V^\sigma$ there exists a sequence $A_n\in V$
    such that $A=\bigcup_n A_n.$
    Consider the sequence defined recursively by
\begin{equation*}
     B_1=A_1,\quad B_n=A_n\less \bigcup_{j<n} A_j\fa n>1.
\end{equation*}
    Since the family $S(V)$ of simple sets generated by the prering $V$
    forms a ring, each of the sets $B_n$ is a union of a finite number
    of disjoint sets from $V.$ Since the sequence $B_n$ consists of disjoint
    sets, combining all the sets from the decompositions of sets $B_n$ into
    a sequence $I_m$ we get
\begin{equation*}
     A=\bigcup_n B_n=\bigcup_m I_m.
\end{equation*}
\ep

\bigskip

The following theorem represents a generalization of Proposition
\ref{pro-approx by V sigma} to the case of arbitrary positive measure
space over a prering.

\bigskip

\begin{thm}[Approximation by summable $V^\sigma$ sets]
\label{thm-approx by V sigma}
Assume that $(X,V,v)$ is a positive measure space.

Then for every summable $B\in V_c$ and every $\varepsilon>0$ there exists
a set $C\in V^\sigma$ such that
\begin{equation*}
 B\subset C\qtext{and}v_c(C)\le v_c(B)+\varepsilon.
\end{equation*}
\end{thm}

\bigskip

\bp

    It follows from the definition of the space $L(v,R)$ of Lebesgue summable
    functions that the support of the characteristic function of the set $B$ can
    be covered by a countable number of the sets from the prering $V.$

    Indeed assume that $s_n$ is a basic sequence converging to the function
    $c_B$ at every point $x\not\in D$ where $D$ is a null set. By definition of
    a null set $D$ can be covered by a countable family $F_1$ of sets
    from the prering $V$ and the support of each simple function $s_n$
    consists of finite number of sets from $V$ thus there exists a countable
    family $F_2$ of set from $V$ covering the supports of all the functions $s_n.$

    Denote by $G$ the set covered by the family $F_1\cup F_2$ of sets.
    Plainly $G\in V^\sigma$ and at every point $x\not\in G$ all functions $s_n(x)$
    have value zero, so their limit must be also $0$ and this means that $c_B(x)=0,$
    and this means that is $x\not\in B.$ Therefore we must have the inclusion $B\subset G.$

        By Proposition \ref{pro-decompose into disjoint sets},
    there exists a sequence of disjoint
    sets
\begin{equation*}
     X_m\in V
\end{equation*}
    such that $G=\bigcup_m X_m.$ Let
\begin{equation*}
     V_m=\set{I\in V:\ I\subset X_m}\qtext{and}v_m(I)=v(I)\fa I\in V_m.
\end{equation*}

    Each triple $(X_m,V_m,v_m)$ forms a measure space and, it is easy to see that,
    any function $f$ whose
    support is in the set $X_m$ belongs to the Lebesgue space $L(v_m,R)$ if and only
    if it belongs to the space $L(v,R).$ Thus for each $m$ the set
\begin{equation*}
     B_m=X_m\cap B
\end{equation*}
    is summable with respect to the measure space $(X_m,V_m,v_m).$

    Now take any $\varepsilon>0.$ Since $X_m\in V_m$ by
    Proposition \ref{pro-approx by V sigma}
    there exist sets
\begin{equation*}
     C_m\in V_m^\sigma\subset V^\sigma
\end{equation*}
    such that
\begin{equation*}
     B_m\subset C_m\qtext{and}v_c(C_m)\le v_c(B_m)+2^{-m}\varepsilon\fa m.
\end{equation*}
    Put $C=\bigcup_m C_m.$ Clearly $C\in V^\sigma$ and
\begin{equation*}
     B\subset C\qtext{and}v_c(C)\le v_c(B)+\varepsilon,
\end{equation*}
    which completes the proof.
\ep

\bigskip

\begin{thm}[Generating outer measure by $V^\sigma$ sets]
Assume that $(X,V,v)$ is a positive measure space over the prering $V.$

For any set $A\subset X$ let $\eta(A)$ be defined by the following formula
if the set of sums
\begin{equation*}
 H=\inf\set{\sum_{n=1}^\infty v(B_n)<\infty:\
B_n\in V,\  A\subset \bigcup_n B_n}
\end{equation*}
 is nonempty then $\eta(A)=\inf H$ else let $\eta(A)=\infty.$

Then the set function $\eta$ coincides with the outer measure $v^*$ on
all subsets $A$ of the space $X.$
\end{thm}

\bigskip

\bp
    We have the following inequality
\begin{equation*}
     v^*(A)\le \eta(A)\fa A\subset X.
\end{equation*}
    Indeed when $\eta(A)=\infty$ the inequality clearly is satisfied. When  $\eta(A)<\infty$
    consider any sequence $B_n\in V$ covering the set $A$ and such that
\begin{equation*}
     \sum_{n=1}^\infty v(B_n)<\infty.
\end{equation*}
    Since $B_n\in V\subset V_c$ by the Monotone Convergence Theorem we get that the set
\begin{equation*}
     B= \bigcup_n B_n
\end{equation*}
    is summable $B\in V_c.$ Since $A\subset B$ thus by definition of the set function
    $v^*$ and by additivity and monotonicity of the measure $v_c$ we have
\begin{equation*}
     v^*(A)\le v_c(B)\le \sum_{n=1}^\infty v_c(B_n)= \sum_{n=1}^\infty v(B_n).
\end{equation*}
    Thus
\begin{equation*}
     v^*(A)\le \eta(A)\fa A\subset X.
\end{equation*}
\bigskip

    Now let us establish the opposite inequality
\begin{equation*}
     \eta(A)\le v^*(A)\fa A\subset X.
\end{equation*}
    This inequality is true when $v^*(A)=\infty,$ so consider the case $v^*(A)<\infty.$
    In this case there exists a summable set $B\in V_c$ such that
\begin{equation*}
     A\subset B\qtext{and}v^*(A)=v_c(B).
\end{equation*}

    Take any $\varepsilon>0.$
    By Theorem \ref{thm-approx by V sigma} there exists a set $C\in V^\sigma$
    such that
\begin{equation*}
     B\subset C\qtext{and}v_c(C)\le v_c(B)+\varepsilon.
\end{equation*}
    By Proposition \ref{pro-decompose into disjoint sets}
    there exists a sequence of disjoint sets $I_n\in V$ such that $C=\bigcup_n I_n.$
    Since $A\subset B\subset C$ we get by definition of the set function $\eta$ that
\begin{equation*}
     \eta(A)\le \sum_n v(I_n)=v_c(C)\le v_c(B)+\varepsilon=v^*(A)+\varepsilon.
\end{equation*}
    Since $\varepsilon$ was fixed but arbitrary the above relation implies
\begin{equation*}
     \eta(A)\le v^*(A)\fa A\subset X.
\end{equation*}

    Hence we have established that for any measure space $(X,V,v)$ the following identity
    is true $v^*(A)=\eta(A)$ for all subsets $A$ of the space $X.$
\ep

\bigskip

\begin{defin}[Outer regularity] Assume that $X$ is a topological space and $(X,V,v)$
a measure space. If for every set $A$ in the prering $V$ we have
\begin{equation*}
 v(A)=\inf\set{v(B):\ A\subset int(B)},
\end{equation*}
where $int(B)$ denotes the interior of the set $B,$
then such a measure space is called {\bf outer regular.}
\end{defin}

\bigskip

Observe that the Riemann measure space $(R,V,v)$ defined over the prering $V$ of
all bounded intervals of $R$ is outer regular. The same is true for the Riemann
measure  space $(R^n,V,v)$ where the prering $V$ consists of all $n$-dimensional cubes.

\bigskip

\begin{thm}
Assume that $X$ is a topological space and $(X,V,v)$ is an outer regular measure
space over the prering $V.$

Define a set function $\lambda$ as follows.
For any set $A\subset X$ if the set of sums
\begin{equation*}
H= \set{\sum_{n=1}^\infty v(B_n)<\infty:\ B_n\in V,\  A\subset \bigcup_n int(B_n)}
\end{equation*}
is nonempty then $\lambda(A)=\inf H$ else let $\lambda(A)=\infty.$

Then the set function $\lambda$ coincides with the outer measure $v^*$ on
all subsets $A$ of the space $X.$
\end{thm}

\bigskip

\bp
    It is obvious that
    \begin{equation*}
         v^*(A)\le \lambda(A)\fa A\subset X.
    \end{equation*}
    To prove the opposite inequality take any $\varepsilon>0$ and
     any sequence of sets $A_n\in V$ such that
    \begin{equation*}
         \sum_n v(A_n)<\infty\qtext{and}A\subset \bigcup_n A_n.
    \end{equation*}
    By outer regularity of the measure space there exist sets $B_n\in V$
    such that
    \begin{equation*}
        A_n\subset int(B_n)\qtext{and}v(B_n)<v(A_n)+2^{-n}\varepsilon\fa n.
    \end{equation*}
    We have
    \begin{equation*}
         A\subset \bigcup_n int(B_n)\qtext{and}\sum_n v(B_n)
         \le \sum_n v(A_n)+\varepsilon<\infty.
    \end{equation*}
    Thus
    \begin{equation*}
         \lambda(A)\le \sum_n v(B_n)\le \sum_n v(A_n)+\varepsilon
    \end{equation*}
    that is
    \begin{equation*}
         \lambda(A)-\varepsilon\le \sum_n v(A_n).
    \end{equation*}
    Since $A_n\in V$ was an arbitrary sequence covering the set $A$ and having
    finite sum $\sum_n v(A_n),$ the above inequality yields
    \begin{equation*}
         \lambda(A)-\varepsilon\le v^*(A).
    \end{equation*}
    Passing to the limit $\varepsilon\to0$ in the above we get
    \begin{equation*}
         \lambda(A)\le v^*(A)\fa A\subset X.
    \end{equation*}

    Hence we have $\lambda(A)=v^*(A)$ for all $A\subset X.$
\ep

\bigskip

\begin{thm}[Outer regularity of outer measure]
\label{Outer regularity of outer measure}
Assume that $X$ is a topological space and $(X,V,v)$ is an outer regular measure
space over the prering $V.$

Then for every $\varepsilon>0$ and every set $A\subset X$
such that $v^*(A)<\infty$ there exists an open set
$G$ such that
\begin{equation*}
 A\subset G\qtext{and}v^*(G)< v^*(A)+\varepsilon.
\end{equation*}
\end{thm}

\bigskip

\bp
    The proof follows from the monotonicity and subadditivity of the outer measure $v^*$
    and from the preceding theorem.
\ep

\bigskip


\bigskip

\section{Properties of Vitali's coverings}

\bigskip

In this section $X$ will denote a close bounded interval $[a,b].$
Let $(X,V,v)$ denote the Riemann measure space over the prering $V$ of
all subintervals $I$ of the interval $X$ and let $v(I)$ denote the length of the interval $I.$
Plainly since $X\in V$ we have $v(I)\le v(X)<\infty$ for all $I\in V.$
\bigskip

\begin{defin}[Vitali covering]
We shall say that a family $H\subset V$
forms a {\bf Vitali covering} of a set $E\subset X$
if $H$ consists of closed sets including the empty set $\emptyset$ and for every point $x\in E$
and every $\varepsilon>0$ the family $H$ contains a set $I$
such that
$x\in I$ and $0<v(I)<\varepsilon.$
\end{defin}

\bigskip

\begin{thm}[Property of Vitali's covering]
\label{Property of Vitali's cover}
Assume that the family $H$  forms a Vitali covering of a set $E\subset X.$
If $H$ does not contain a finite disjoint family of closed sets $I_j$ covering
the set $E$ then one can find a sequence of disjoint closed sets $I_j\in H$ such that
\begin{enumerate}
\item for each $j$ the intersection $E\cap I_j$ is nonempty,
\item for each $n$ the family $H_n$ of sets defined by
\begin{equation*}
 H_n=\set{I\in H:\ E\cap I\not=\emptyset,\ 0<v(I),\ I\cap I_j=\emptyset\ \forall\ j\le n}
\end{equation*}
is nonempty,
\item for each $n$ we have $I_{n+1}\in H_n$ and
$
  \varepsilon_n<2v(I_{n+1}),
$
where
$$
\varepsilon_n=\sup\set{v(I):\ I\in H_n}.
$$
\end{enumerate}
\end{thm}

\bigskip

\bp
    Take any point $x\in E.$ Let $I_1$ denote any interval from $H$ containing the point~$x.$
    Since $I_1$ by assumption cannot cover $E$ the set $E\less I_1$ is nonempty.
    Thus there exists a point $x_1$ in $E$ that does not belong to $I_1$ that is
    $x_1$ belongs to the set $X\less I_1.$ Since complement of a closed set is open
    there is an open interval $J$ containing $x_1$ such that $J\subset X\less I_1.$
    By definition of Vitali's covering there is a closed set $I\in H$
    of sufficiently small length such that $x_1\in I$ and
    $v(I)>0$ and $I\subset J.$ Hence the family of sets
\begin{equation*}
     H_1=\set{I\in H:\ E\cap I\not=\emptyset,\ 0<v(I),\ I\cap I_1=\emptyset}
\end{equation*}
    is nonempty.
    Thus the set $R_1=\set{v(I):\ I\in H_1}$ is nonempty and bounded from above
    by $v(X).$ Thus by axiom of completeness of the space $R$ of reals the least upper
    bound $\varepsilon_1$ of the set $R_1$ is well defined and we have
\begin{equation*}
     0<\varepsilon_1=\sup\set{v(I):\ I\in H_1}.
\end{equation*}
    The number $\frac{1}{2}\varepsilon_1$ is smaller than the least upper bound
    $\varepsilon_1$ of the set $R_1.$ Therefore there exists a set $I_2\in H_1$
    such that $\frac{1}{2}\varepsilon_1<v(I_2),$ that is
\begin{equation*}
     \varepsilon_1<2v(I_2).
\end{equation*}

    Clearly the pair $I_1,\,I_2$ forms disjoint sets and the conditions (2) and
    (3) of the theorem are satisfied for $n=1.$

    It is also clear that we can continue by induction to obtain a sequence $I_j\in H$
    of disjoint sets, and a sequence $H_n$ of subfamilies of $H,$ and a sequence
    $\varepsilon_n$ of positive numbers, satisfying
    the conditions (1), (2), and (3) as required in the theorem.
\ep

\bigskip

\begin{defin}[Sparse sequence]
Assume that $H$ forms a Vitali covering of a set $E\subset X.$ Any sequence of
sets $I_j\in H$ satisfying the conditions (1), (2), and (3), of
Theorem \ref{Property of Vitali's cover}
will be called a {\bf sparse sequence} with respect to the set $E.$
\end{defin}

\bigskip

\begin{pro}[For sparse sequence $\lim_n \varepsilon_n=0$]
\label{pro lim epsilon sub n = 0}
If $I_j\in H$ is a sparse sequence relatively to the set $E$ and $\varepsilon_j$
is the associated sequence satisfying the condition (3) of Theorem
\ref{Property of Vitali's cover}, then $\varepsilon_j\to 0$ as $j\to \infty.$
\end{pro}

\bigskip

\bp
    Since $X\in V$ the family $V_c$ of summable sets forms a sigma algebra. Thus the set
    $A=\bigcup_n I_n$ belongs to $V_c.$ Since the sets $I_n$ are disjoint, we have
\begin{equation*}
     \sum_n v(I_n)=\sum_n v_c(I_n)=v_c(A)\le v_c(X)=v(X)<\infty.
\end{equation*}

    From the convergence of the above series follows that $v(I_n)\to 0.$
    Since from condition (3) of Theorem \ref{Property of Vitali's cover}
    follows that $\varepsilon_n<2v(I_{n+1})$ for all $n$
    we must have $\varepsilon_n\to 0.$
\ep

\bigskip
\begin{defin}[Set operator $S$]
For any interval $I\in V$ define the set operator $S:V\mto V$ as follows:
If $I$ is the empty set $\emptyset,$ let $S(I)=\emptyset$ else let
\begin{equation*}
S(I)=[c-2v(I),d+2v(I)]\cap X
\end{equation*}
where $c\le d$ are the end points of the interval $I.$
\end{defin}
\bigskip

Notice that we have a convenient estimate $v(S(I))\le 5v(I)$ for every $I\in V.$

\bigskip

\begin{lemma}\label{lem on 2 points}
Given an interval $I\in V$ and two points $x,y\in X.$
If
\begin{equation*}
 |x-y|\le 2v(I)\qtext{ and }y\in I
\end{equation*}
then the point $x$ belongs to the interval $S(I).$
\end{lemma}

\bigskip

\bp
    The proof is obvious.
\ep

\bigskip

\begin{defin}[Interlaced cover]
We shall say that the sequence
$I_j\in V$ forms an {\bf interlaced cover} of a set $E\subset X$ if the
sets $I_j$ are closed and disjoint and
\begin{equation*}
 E\subset (\bigcup_{j\le n}I_j)\cup (\bigcup_{j> n}S(I_j))\fa n=1,2,\ldots
\end{equation*}
\end{defin}

\bigskip

\begin{thm}[Vitali Covering Theorem]
\label{Vitali Covering Theorem}
Assume that $E\subset X$ and $H$ forms a Vitali covering of the set $E.$
Then there exists a sequence $I_j\in H$ of intervals forming an
interlaced covering of the set $E.$
\end{thm}

\bigskip

\bp
    If the set $E$ can be covered by a finite number of disjoint sets
     $I_1,\,I_2,\ldots,I_n$ from the family $H,$ by setting $I_j=\emptyset$
     for $j>n$ we get an interlaced cover from the family $H$ of the set $E.$

    If a finite cover from the family $H$ does not exist then there exists
    a sparse  sequence $I_j\in H$ relatively to the set $E.$
    Let the sequence of families $H_n\subset H$
    and the sequence $\varepsilon_n$ be defined as in
    Theorem \ref{Property of Vitali's cover}. Since the sequence $H_n$ of families
    is decreasing the sequence $\varepsilon_n$ is non-increasing.

    We shall prove that the sequence of sets $I_j$ forms
    an interlaced cover of the set $E.$
    To this end take any index $n$ and consider the set $$B_n=E\less \bigcup_{j\le n}I_j.$$
    Since the set $B_n$ is nonempty there exist a point $x\in E$ and an interval $I\in H$
    such that $x\in I$ and $0<v(I)$ and $I\cap I_j=\emptyset$ for all $j\le n.$
    Thus $I\in H_n$ and so $v(I)\le \varepsilon_n.$
    By Proposition~\ref{pro lim epsilon sub n = 0} we must have $\varepsilon_m\to 0.$
    Thus for all sufficiently large indexes $m$ we must have
\begin{equation*}
     \varepsilon_m<v(I).
\end{equation*}
    Thus the set of indexes
\begin{equation*}
     K=\set{m:\ v(I)\le \varepsilon_m}
\end{equation*}
    contains the index $n$ and is at most finite. Let $k$ be the largest index
    belonging to $K.$ Thus we must have $k\ge n$ and
\begin{equation*}
     \varepsilon_{k+1}<v(I)\le\varepsilon_k.
\end{equation*}

    The relation $v(I)\le \varepsilon_k$ implies $I\in H_k$ and thus
    $I\cap I_j=\emptyset$ for all $j\le k.$ The relation $\varepsilon_{k+1}<v(I)$
    implies that it is not true that $I\cap I_j=\emptyset$ for all $j\le (k+1).$
    Hence we must have $I\cap I_{k+1}\not=\emptyset.$ So there exists a point $y\in X$
     such that $y\in I$ and $y\in I_{k+1}.$ Thus since $x\in I$ we have
\begin{equation*}
      |x-y|\le v(I)\le\varepsilon_k\le 2v(I_{k+1})\qtext{and}y\in I_{k+1}.
\end{equation*}
    By Lemma~\ref{lem on 2 points} we must have $x\in S(I_{k+1}),$ which implies
    the inclusion
\begin{equation*}
     E\subset (\bigcup_{j\le n}I_j)\cup (\bigcup_{j> n}S(I_j)).
\end{equation*}
    Since the index $n$ was fixed but arbitrary the theorem is established.
\ep

\bigskip


\bigskip

\section{Lebesgue Theorem on differentiability of monotone functions}

\bigskip

In this section we will prove Lebesgue's theorem asserting that
every real-valued monotone function defined on an interval $I$ has a
finite derivative almost everywhere on $I.$ To this end
it is convenient to introduce derivatives
known in the literature as Dini's derivatives.

\bigskip

\subsection{Dini's derivatives}

\bigskip
Consider a real-valued function $f$ defined on an open interval $(a,b).$

\begin{defin}

By a right-sided upper {\bf Dini derivative} of $f$ at a point $x\in (a,b)$
we shall understand the  finite or infinite limit
\begin{equation*}
     D_r^+ f(x)=\limsup_{h>0,\,h\into 0}\frac{1}{h}(f(x+h)-f(x)).
\end{equation*}

Similarly we define the right-sided lower Dini derivative by
\begin{equation*}
     D_r^- f(x)=\liminf_{h>0,\,h\into 0}\frac{1}{h}(f(x+h)-f(x)),
\end{equation*}
and left-sided upper derivative
\begin{equation*}
     D_l^+ f(x)=\limsup_{h>0,\,h\into 0}\frac{1}{h}(f(x)-f(x-h)),
\end{equation*}
and left-sided lower derivative
\begin{equation*}
     D_l^- f(x)=\liminf_{h>0,\,h\into 0}\frac{1}{h}(f(x)-f(x-h)).
\end{equation*}
\end{defin}

\bigskip

The function $f$ has derivative $f'(x)$ at a point $x$
if and only if all four Dini's derivatives are finite and have equal values.
We have the following relations between Dini's derivatives.
\begin{equation*}
 D_r^- f(x)\le D_r^+ f(x)\qtext{and}D_l^- f(x)\le D_l^+ f(x)\fa x\in (a,b).
\end{equation*}

\bigskip

By an {\bf increasing} function $f$ on an interval $I$ we shall understand a function
satisfying the following condition
\begin{equation*}
 f(x)\le f(y)\fa x<y,\ x,y\in I.
\end{equation*}

Such a function is also called in the literature nondecreasing.

By an {\bf decreasing} function $f$ on an interval $I$ we shall understand a function
satisfying the condition
\begin{equation*}
 f(x)\ge f(y)\fa x<y,\ x,y\in I.
\end{equation*}

A function that is either increasing or decreasing on an interval is called
{\bf monotone} on that interval.

\bigskip

\subsection{Differentiability of monotone functions}
\bigskip

\begin{thm}[Lebesgue]
Let $X$ denote a closed bounded interval $[a,b]$ and assume that
$(X,V,v)$ denote the Riemann measure space on the prering $V$ of all subintervals of $X.$

If $f$ is a monotone real-valued function on $X,$
then the derivative $f'(x)$ exists in the open interval $(a,b)$ and is finite for
almost all $x.$
\end{thm}

\bigskip

\bp
    We may assume without loss of generality that the function $f$ is increasing
    otherwise we would consider the function $g(x)=-f(x)$ for all $x\in I.$

    Introduce a set function $w:V\mto R$ by the formula
\begin{equation*}
     w(I)=f(\beta)-f(\alpha),
\end{equation*}
    where $\alpha\le\beta$ are the end points of the interval $I.$
    The set function $w$ is finitely additive on the prering $V$
    and as such it admits a unique extension to a finitely additive set
    function on the ring $S(V)$ of simple sets generated by $V.$
    We shall use the same symbol $w$ to denote that extension.
    Notice that $w(A)\ge 0$ for all sets $A\in S(V)$ and as a consequence
\begin{equation*}
     w(A)\le w(B)\le w(X)\qtext{if} A\subset B;\ A,B\in S(V).
\end{equation*}

    Let $c(x)$ and $d(x)$  denote any two Dini derivatives of the function $f(x)$
    on the interval $(a,b).$
    We will prove that the functions $c$ and $d$ are equal almost
    everywhere on $X.$

    The following property follows from the definition
    of $\limsup$ and $\liminf\!.$ If $c(x)$ denotes either lower or upper right-sided
    Dini derivative and $c(x)<u,$ then there exists a sequence $h_n>0$ such that
    $h_n\to 0$ and
\begin{equation*}
     f(x+h_n)-f(x)<u\,h_n\fa  n.
\end{equation*}
    Thus denoting by $I_n$ the interval with the end points $x$ and $x+h_n$
    we get
\begin{equation}\label{estimate one}
     w(I_n)<u\,v(I_n)\fa n.
\end{equation}

    Similar argument can be applied to the left-sided derivatives yielding the
    existence of a sequence of intervals $I_n$ with the end points at $x-h_n$
    and $x$ such that the estimate (\ref{estimate one}) holds.

    Proceeding in a similar manner we can prove that if $d(x)$ denotes any Dini derivative,
    left or right-sided, upper or lower, and if $s<d(x)$ then there exists
    a sequence of intervals $I_n$ each of positive lengths
    and such that $x\in I_n$ and $v(I_n)\to 0,$ when $n\to \infty,$
    and
\begin{equation*}
     s\,v(I_n)<w(I_n)\fa n.
\end{equation*}

    Now consider the set
\begin{equation*}
     E=\set{x\in (a,b):\ c(x)<d(x)}.
\end{equation*}
    Since between any two real numbers one can find a pair of rational numbers
\begin{equation}\label{u less than s}
    u<s,
\end{equation}
the set $E$ is equal to the countable union of the sets
\begin{equation*}
     E_{u,s}=\set{x\in (a,b):\ c(x)<u<s<d(x)}.
\end{equation*}

    Assuming that the set $E$ is not a null set, will yield a contradiction.
     Indeed from countable subadditivity
    of the outer measure $v^*$ we get that at least for one of the sets $E_{u,s} $
    we must have
\begin{equation*}
     0<t=v^*(E_{u,s}).
\end{equation*}
    Take any $\varepsilon$ such that $0<\varepsilon<t/2$
    and let $G$ denote an open set containing the
    set $E_{u,s}$ and such that
\begin{equation}\label{outer meas of G bound}
    v^*(G)<t+\varepsilon.
\end{equation}

    Since Riemann measure space is regular such an
    open set exists in accord with Theorem
    \ref{Outer regularity of outer measure}
    on outer regularity of outer measure~$v^*.$

    Consider the collection of closed sets
\begin{equation}\label{collection H}
     H=\set{I\subset G:\ I=\emptyset\text{ or } I\cap E_{u,s}\not=\emptyset,\ 0<v(I),
     \ w(I)<u\,v(I)}.
\end{equation}
    It follows from the previous considerations that the collection $H$
    forms a Vitali covering of the set $E_{u,s}$

    By Vitali covering theorem there exists a sequence of disjoint intervals $I_k$
    forming an interlaced covering of the set $E_{u,s}$ that is
\begin{equation*}
     E_{u,s}\subset [\bigcup_{k\le n}I_k]\cup [\bigcup_{k>n}S(I_k)]\fa n.
\end{equation*}
    Let $J_k$ denote the interior of the interval $I_k.$
    Notice the estimate
\begin{equation*}
     E_{u,s}\less \bigcup_{k\le n}J_k\subset
    [\bigcup_{k\le n}(I_k\less J_k)]\cup [\bigcup_{k>n}S(I_k)]\fa n.
\end{equation*}
    Since each set $I_k\less J_k$ consists of just two points,
    their union forms a null set.
    Thus from the preceding set estimate we get
\begin{equation*}
     v^*(E_{u,s}\less \bigcup_{k\le n}J_k)\le 5\sum_{k>n}v(I_k)<\varepsilon
\end{equation*}
    for sufficiently large $n.$

    Denote by $G_u$ the open set $\bigcup_{k\le n}J_k$ and notice the inequality
\begin{equation*}
     t=v^*(E_{u,s})\le v^*(E_{u,s}\cap G_u)+v^*(E_{u,s}\less G_u)
\end{equation*}
    yielding the estimate
\begin{equation}\label{lower est of v*(E sub us cap G sub u}
     0<t-\varepsilon\le v^*(E_{u,s}\cap G_u).
\end{equation}
    From the definition of the collection $H$ of sets (\ref{collection H})
    and the fact that $G_u$ forms a simple set from the ring $S(V),$
    and the estimate (\ref{outer meas of G bound}) on $v^*(G),$ follows that
\begin{equation}\label{upper est on v* of G sub u}
 \begin{split}
    w(G_u)&=\sum_{k\le n}w(J_k)=\sum_{k\le n}w(I_k)\\
    &\le u\,\sum_{k\le n}v(I_k)\le u\,v(G)\le u\,t+u\,\varepsilon.
\end{split}
\end{equation}
\bigskip

    Now consider the collection of closed sets
\begin{equation}\label{collection H subs s}
     H_s=\set{I\subset G_u:\ I=\emptyset\text{ or }
     I\cap E_{u,s}\cap G_u\not=\emptyset,\ 0<v(I),
     \ s\,v(I)<w(I)}.
\end{equation}
    This collection forms a Vitali covering of the set $E_{u,s}\cap G_u.$
    Let $I'_k$ denote the sequences of disjoint sets forming an interlaced
    covering of the set. Thus we have
\begin{equation*}
     E_{u,s}\cap G_u\subset [\bigcup_{k\le n}I'_k]\cup [\bigcup_{k>n}S(I'_k)]\fa n.
\end{equation*}
    Again for sufficiently large index $n$ we get the estimate
\begin{equation}\label{upper estimate of v* of E sub us cap G sub u}
     v^*(E_{u,s}\cap G_u)\le \sum_{k\le n}v^*(I'_k)
     +\varepsilon= \sum_{k\le n}v(I'_k)+\varepsilon.
\end{equation}
    By the lower estimate
    (\ref{lower est of v*(E sub us cap G sub u}) of $v^*(E_{u,s}\cap G_u)$, and
    the upper estimate (\ref{upper estimate of v* of E sub us cap G sub u}),
    and the upper estimate (\ref{upper est on v* of G sub u}) on $w(G_u)$,
    we get
\begin{equation*}
\begin{split}
     s\,(t-\varepsilon)
     &\le s\,v^*(E_{u,s}\cap G_u)\le s\sum_{k\le n}v(I'_k)+s\,\varepsilon\\
    &\le \sum_{k\le n}w(I'_k)+s\,\varepsilon
    \le w(\bigcup_{k\le n}I'_k)+s\,\varepsilon\\
    &\le w(G_u)+s\,\varepsilon
    \le u\,t +u\,\varepsilon+s\,\varepsilon
\end{split}
\end{equation*}
    that is
\begin{equation*}
     s\,(t-\varepsilon)\le ut +u\varepsilon+s\varepsilon\fa \varepsilon\in (0,t/2).
\end{equation*}
    Thus passing to the limit $\varepsilon\to 0$ in the above inequality
    and dividing by $t$ we get $s\le u,$ which contradicts the assumption
    (\ref{u less than s}) that $u<s.$
\bigskip

    Switching the roles of the functions $c$ and $d$ we get that the set
    where $c(x)>d(x)$ also forms a null set.
    Thus we have proved that any increasing function $f$ on an interval
    has a derivative almost everywhere inside of the interval.
\bigskip

    Now let us prove that the derivative is finite almost everywhere.
    Since for an increasing function all the difference quotients
    are nonnegative, the derivative $f'(x)$ may take on only nonnegative values.
    So we need to consider only the set
\begin{equation*}
     E=\set{x\in (a,b):\ f'(x)=\infty}.
\end{equation*}
    Take any $\varepsilon>0$ and
    chose a sufficiently large number $r>0$ so that
\begin{equation*}
     w(X)/r\le \varepsilon/2.
\end{equation*}
    Introduce a collection of closed sets
\begin{equation*}
     H=\set{I\in V:\ I=\emptyset\text{ or }I\cap E\not=\emptyset,\ 0<v(I),\ r\,v(I)<w(I)}.
\end{equation*}
    The collection $H$
    forms a Vitali covering of the set $E.$
    Let $I_k$ be a sequence of disjoint intervals forming an interlaced
    cover of the set $E$ that is
\begin{equation*}
     E\subset [\bigcup_{k\le n}I_k]\cup [\bigcup_{k>n}S(I_k)]\fa n.
\end{equation*}
    Since
\begin{equation*}
     \sum_{k>0} v(I_k)=\sum_{k>0} v_c(I_k)= v_c(\bigcup_{k>0}I_k)\le v_c(X)=f(b)-f(a)<\infty,
\end{equation*}
    the above series converges and thus its remainder tends to zero.
    Choosing sufficiently large index $n$ we can get
\begin{equation*}
     5\sum_{k>n} v(I_k)\le \varepsilon/2.
\end{equation*}
    We have the estimate
\begin{equation*}
\begin{split}
    v^*(E)&\le \sum_{k\le n}v(I_k)+\sum_{k>n}v(S(I_k))\\
    &\le  \frac{1}{r}\sum_{k\le n}w(I_k)+5\sum_{k>n}v(I_k)\\
    &\le \frac{1}{r}w(\bigcup_{k\le n}I_k)+5\sum_{k>n}v(I_k)\\
    &\le \frac{1}{r}w(X)+5\sum_{k>n}v(I_k)\le \varepsilon.
\end{split}
\end{equation*}
    Since $\varepsilon$ was fixed but arbitrary this implies $v^*(E)=0$
    that is the set $E$ is a null set.

    Hence we have proved that the derivative $f'(x)$ exists and is finite
    for almost all points $x$ in the interval $(a,b).$
\ep

\bigskip


\bigskip

\section{Properties of monotone functions}

\bigskip

In this section we shall consider further properties of monotone functions. First notice
that since a function that is differentiable at a point is continuous at that point,
every monotone function on an interval is continuous almost everywhere on that interval.

\bigskip

\begin{thm}[Summability of monotone functions]
Let $X$ denote a closed bounded interval $[a,b]$ and $(X,V,v)$ the Riemann measure space
on the prering $V$ of all subintervals of $X.$

If $f$ is a real-valued monotone function on the interval $X$
then it is Lebesgue summable that is $f\in L(v,R).$
\end{thm}

\bigskip

\bp
    We may assume without loss of generality that the function $f$ is increasing.
    Split the interval $[a,b]$ into two disjoint intervals $I_{1,j}$
    of equal length and let $a_{1,j}$ denote the left end
    of the interval $I_{1,j}$ where   $j=1,2.$
    Put
\begin{equation*}
     s_1(x)=f(a_{1,1})c_{I_{1,1}}(x)+f(a_{1,2})c_{I_{1,2}}(x)\fa x\in X.
\end{equation*}
    Notice that
\begin{equation*}
     s_1(x)\le f(x)\fa x\in X.
\end{equation*}
    Proceeding in similar manner with each of the intervals $I_{1,j}$
    by dividing them into two intervals, and proceeding by induction
    we will obtain an increasing sequence of simple functions
\begin{equation*}
     s_n(x)=\sum_{j\le 2^n}f(a_{nj})c_{I_{nj}}(x)\le f(x)\fa x\in X.
\end{equation*}
    Notice that the sequence $s_n$ of simple functions is bounded from above
    by the summable function $f(b)c_X.$
\bigskip

    The sequence of values $s_n(x)$ converges at every point of continuity
    of the function $f$ to the value $f(x).$

    Indeed, take any positive $\varepsilon$ and select $\delta>0$ so that
\begin{equation*}
     |f(x)-f(y)|<\delta\qtext{if} y\in (x-\delta,x+\delta)\subset (a,b).
\end{equation*}
    Select $k$ so that
\begin{equation*}
     2^{-n}v(X)<\delta\fa n\ge k.
\end{equation*}
    Since each of the intervals $I_{nj}$ have the same length
\begin{equation*}
     v(I_{nj})=2^{-n}v(X)
\end{equation*}
    and they form a disjoint decomposition of the interval $X$
\begin{equation*}
     X=\bigcup_{j\le 2^n}I_{nj},
\end{equation*}
    we must have that
    for each $n>k$ there is exactly one interval $I_{nj}$ containing the point
    $x.$ Since for each such interval we have
\begin{equation*}
     I_{nj}\subset (x-\delta,x+\delta)\fa n>k
\end{equation*}
    we see that
\begin{equation*}
     |s_n(x)-f(x)|=|f(a_{nj})-f(x)|<\varepsilon\fa n>k.
\end{equation*}
    Thus $s_n(x)\to f(x)$ at every point of continuity of the function $f.$
    It follows from the Monotone Convergence Theorem that the function $f$
    belongs to the space $L(v,R)$ of Lebesgue summable functions.
\ep

\bigskip

\subsection{Differentiability of series of monotone functions}

For a reference to the following theorem see Fubini \cite{fubini-ser}.

\bigskip

\begin{thm}[Fubini]
Let $X$ denote a closed bounded interval $[a,b]$ and $(X,V,v)$ the Riemann measure space
on the
prering $V$ of all subintervals of $X.$

Assume that $h_n$ is a sequence of functions such
that all functions of the sequence are either
increasing or all are decreasing.

If the series
\begin{equation*}
 \sum_n h_n(x)=h(x)\fa x\in X
\end{equation*}
converges to a finite value for each $x,$
then there exists a null set $D$ such that the derivative
$h'(x)$ exists and
we have the equality
\begin{equation*}
 \sum_n h'_n(x)=h'(x)\fa x\in X\less D.
\end{equation*}
\end{thm}

\bigskip

\bp
    Without loss of generality we may assume that all functions of the sequence
    are increasing and $h_n(a)=0$ for all $n.$
    Assume that for each $n$ the derivative $h'_n(x)$ exists for $x\in X\less D_n$
    where $D_n$ is a null set.
    Since the function $h$ as a sum of a series of increasing functions is itself
    increasing the derivative $h'(x)$ exists for all $x\not\in D_0,$ where $D_0$
    is a null set. Let $D$ be the null set representing the union
\begin{equation*}
     D=\bigcup_{n\ge 0}D_n.
\end{equation*}
    Take any point $x\in X\less D$ and consider any point $y\in X$ such that $x<y.$
    Since the
    corresponding difference quotients of functions $h_n$ are all nonnegative,
     we must have the estimate
\begin{equation*}
     \sum_{n\le m}\frac{h_n(y)-h_n(x)}{y-x}\le \frac{h(y)-h(x)}{y-x}\fa m.
\end{equation*}
    Thus for a fixed index $m$ passing to the limit $y\to x$ we get
\begin{equation*}
     \sum_{n\le m}h'_n(x)\le h'(x)\fa m\qtext{and}x\not\in D.
\end{equation*}
    Since all terms $h'_n(x)$ are nonnegative and the above partial sums are
    bounded  we must have
\begin{equation*}
     \sum_{n=1}^\infty h'_n(x)\le h'(x)\fa x\not\in D.
\end{equation*}
    The above implies
\begin{equation*}
     h'_n(x)\to 0\fa x\not\in D.
\end{equation*}

    Introduce notation for the partial sums
\begin{equation*}
     s_n=h_1+h_2+\cdots+h_n\fa n
\end{equation*}
    and select a subsequence $s_{k_n}$ so that
\begin{equation*}
     |s_{k_n}(b)-h(b)|\le 2^{-n}\fa n.
\end{equation*}
    Notice that the series with terms $f_n=(h-s_{k_n})$
    consists of increasing functions and it converges for all $x\not\in D$
    to a finite function.
    Thus from the previous considerations we can conclude that
\begin{equation*}
     f'_n(x)=h'(x)-s'_{k_n}(x)\to 0\fa x\not\in D.
\end{equation*}
    Since the sequence of partial sums $s_n$ is increasing the above relation implies
    the convergence
\begin{equation*}
     s'_n(x)\to h'(x)\fa x\not\in D,
\end{equation*}
    which is equivalent to
\begin{equation*}
     \sum_{n=1}^\infty h'_n(x)=h'(x)\fa x\not\in D.
\end{equation*}
\ep

\bigskip

\begin{cor}[Fubini]
\label{cor-Fubini for increasing sequence}
Let $I\subset R$ be any interval of reals.
Assume that $f_n$ is a sequence of increasing functions on the interval $I.$

If the sequence $f_n$ increasingly converges to a finite-valued function
$f$ at every point of the interval $I$,
then the function $f$ is differentiable almost everywhere on the interval $I$
and  moreover the derivatives $f'_n(x)$ increasingly
converge to the derivative $f'(x)$  almost everywhere on $I.$
\end{cor}

\bigskip


\bigskip

\section{Invariant measures induce invariant integrals}

\bigskip

\begin{defin}[Invariant measure]
Let $X$ be an abstract space and $V$ a prering of subsets of $X.$
Assume that the triple $(X,V,v)$ forms a positive measure space
and $T$ represents a map from the space $X$ into $X.$

We shall say that the measure space $(X,V,v)$ is {\bf invariant} under
the map $T$ if
\begin{equation*}
 T^{-1}(A)\in V\qtext{and}v(T^{-1}(A))=v(A)\fa A\in V.
\end{equation*}
\end{defin}

\bigskip

The Riemann measure space $(R,V,v)$ over the reals $R$ is invariant under
any translation map
\begin{equation*}
 T(x)=x+h\fa x\in R,
\end{equation*}
since a translation operation maps intervals onto intervals and preserves the length
of the interval. The Riemann measure space
is also invariant under the reflection $x\mto -x.$
The same is true for Riemann measure space $(R^n,V,v).$

The counting measure space $(X,V,v)$ is invariant under any map $T$
of $X$ onto $X$ which is one-to-one.

\begin{thm}[Invariant integral]
Let $X$ be an abstract space and $V$ a prering of subsets of $X.$
Assume that the triple $(X,V,v)$ forms a positive measure space
and $Y$ a Banach space. Assume that $T$ is a map from the space $X$ into $X.$

If the measure space $(X,V,v)$ is invariant under the map $T,$
then the operator
\begin{equation*}
 f\mto f\circ T\fa f\in L(v,Y)
\end{equation*}
maps the space $L(v,Y)$ of Bochner summable functions into $L(v,Y)$
and preserves the integral, that is,
\begin{equation*}
 \int f\,dv=\int f\circ T\,dv\fa f\in L(v,Y).
\end{equation*}
\end{thm}

\bigskip

\bp
    Notice that the theorem is valid for simple functions.
    Indeed if $A\in V$ then for the characteristic function of $A$ we have
\begin{equation*}
     c_A\circ T(x)=1\iff T(x)\in A\iff x\in T^{-1}(A)\iff c_{T^{-1}(A)}(x)=1.
\end{equation*}
    Thus for a simple function
\begin{equation*}
     s=y_1c_{A_1}+\cdots y_kc_{A_k}
\end{equation*}
    we get
\begin{equation*}
     s\circ T=y_1c_{T^{-1}(A_1)}+\cdots y_kc_{T^{-1}(A_k)}
\end{equation*}
    and therefore
\begin{equation*}
\begin{split}
     \int s\circ T\,dv&=y_1\int c_{T^{-1}(A_1)}\,dv+\cdots +y_k\int c_{T^{-1}(A_k)}\,dv\\
    &=y_1v(T^{-1}(A_1))+\cdots +y_kv(T^{-1}(A_k))\\
    &=y_1v(A_1)+\cdots +y_kv(A_k)=\int s\, dv.\\
\end{split}
\end{equation*}

    As a consequence the transformation $f\mto f\circ T$ maps a basic sequence into
    a basic sequence. So take any Bochner summable function $f\in L(v,Y).$
    By definition of the space $L(v,Y)$ there exists a basic sequence $s_n\in S(V,Y)$
    and a null set $A$ such that
\begin{equation*}
     s_n(x)\into f(x)\qtext{when}x\not\in A.
\end{equation*}
    Thus
\begin{equation*}
     s_n\circ T(x)\into f\circ T(x)\qtext{when}T(x)\not\in A
\end{equation*}
    or equivalently
\begin{equation*}
     s_n\circ T(x)\into f\circ T(x)\qtext{when}x\not\in T^{-1}(A).
\end{equation*}
    Clearly the sequence $s_n\circ T$ is also basic.
    Since the map $A\mto T^{-1}(A)$ considered as a map from the power set $\Power(X)$
    into itself, is monotone and maps unions of sets into unions of their images,
    the set $T^{-1}(A)$ is a null set.
    Thus we must have that $f\circ T\in L(v,Y)$
    and
\begin{equation*}
     \int f\circ T\,dv=\lim_n\int s_n\circ T\,dv=\lim_n\int s_n\,dv=\int f\,dv.
\end{equation*}
\ep

\bigskip



\section{Integration over the space $R$ of reals}

\bigskip

If $(R,V,v)$ is the Riemann measure space and $X$ any subinterval
of $R$ then  let $V$ denote the
prering of all bounded subintervals of $X,$ and $v(A)$ the length of the
interval $A\subset X.$ Clearly the space $(X,V,v)$ is a measure
space as a subspace of the Riemann measure space. We shall call
this measure space $(X,V,v)$ the {\em Riemann measure space over the interval $X.$}


For the case of a Riemann measure space over an interval $X$
we shall use the customary
notation for the integral of a Bochner summable function $f\in
L(v,Y).$ We shall write
\begin{equation*}
    \begin{split}
    \int_{t_1}^{t_2}f(t)\,dt&=\quad\int c_{[t_1,t_2]}f\,dv
        \qtext{if}t_1\le t_2,\ t_1,t_2\in X,\\
    \int_{t_1}^{t_2}f(t)\,dt&=-\int c_{[t_2,t_1]}f\,dv
        \qtext{if}t_1> t_2,\ t_1,t_2\in X.\\
    \end{split}
\end{equation*}

\bigskip


Adopting the above notation yields a convenient formula for any
$f\in L(v,Y)$
\begin{equation*}
    \int_{t_1}^{t_2}f(t)\,dt+\int_{t_2}^{t_3}f(t)\,dt+\int_{t_3}^{t_1}f(t)\,dt=0\fa
    t_1,t_2,t_3\in X.
\end{equation*}

%
%
\bigskip

\begin{defin}[Indefinite integral]
Let $(X,V,v)$ denote a  Riemann measure space over an interval $X$ and $Y$ a Banach space.
Let $f\in L(v,Y)$ be a Bochner summable function.

By an {\bf indefinite integral} of the function $f$ we shall understand
any function of the form
\begin{equation*}
 F(x)=\int_a^xf(t)\,dt\fa x\in X\text{ and some }a\in X.
\end{equation*}
\end{defin}

\bigskip

We shall prove later a theorem due to Lebesgue that indefinite integral $F$ of a function
$f,$ summable with respect to Riemann measure, has a derivative almost everywhere and its
derivative $F'$ coincides with the function $f$ almost everywhere, but first notice
the following simple consequence of continuity of the integrand $f.$

\bigskip


\begin{pro}[Differentiability at continuity points]
Let $X$ be an interval and let $(X,V,v)$ denote the
Riemann measure space over $X$ and $Y$ any Banach space.

Let $f\in L(v,Y)$ be a Bochner summable function and let $F$ denote its indefinite integral
\begin{equation*}
 F(x)=\int_a^x f(t)\,dt\fa x\in I.
\end{equation*}

If the function $f$ is continuous at a point $p\in I,$
then the function $F$ is differentiable at $p$
and $F'(p)=f(p).$
\end{pro}

\bigskip
\bp
    Take any $\varepsilon>0 $ and select $\delta>0$ so that
\begin{equation*}
     |f(p+h)-f(p)|\le \varepsilon\qtext{if}|h|\le\delta\ \text{ and }\ x,\,x+h\in X.
\end{equation*}
Then we have
\begin{equation*}
     \left|\frac{F(p+h)-F(p)}{h}-f(p)\right|=
    \left|\frac{1}{h}\int_{p}^{p+h}(f(t)-f(p))\,dt\right|\le
    \varepsilon
\end{equation*}
if $0<|h|\le\delta.$    Thus $F'(p)=f(p).$
\ep

%
%

\bigskip

\begin{thm}[Lipschitzian property of indefinite integral]
Let $(I,V,v)$ denote the  Riemann measure space and $Y$ a Banach space.
Assume that $f$ is a Bochner summable function
on the interval $I$ and $F$ its indefinite integral.

If for some constant $M$ we have $|f(x)|\le M$ almost everywhere on $I$ then
\begin{equation*}
 |F(x)-F(y)|\le M|x-y|\fa x,\,y\in I,
\end{equation*}
that is the function $F$ is Lipschitzian on the interval $I.$
\end{thm}

\bigskip

\bp
    The proof is straightforward and we leave it to the reader.
\ep

\bigskip


\begin{defin}[Local summability]
Assume that $(X,V,v)$ is a positive measure space over a prering $V$ of
an abstract space $X.$

A function $f$ from  $X$ into a Banach space $Y$ is
said to be {\bf locally summable}  if for every set $A\in V$
the function $c_Af$ belongs to the space $L(v,Y)$
of Bochner summable functions.
\end{defin}

\bigskip

Clearly every summable function $f\in L(v,Y)$ is locally summable on $X$
but if the set $X$ is not in the prering $V$ then  function $c_X$ is locally
summable on $X$ but it does not have to be summable on $X.$

Consider
for example an infinite set $X$ and the counting measure, or the space $R$
of reals and the Riemann measure on the prering of all bounded intervals.

\bigskip


\begin{defin}[Newton's formula]
Assume that $F$ is a function from an interval $I\subset R$ into a Banach space $Y.$
We shall say that the function $F$ satisfies {\bf Newton's formula} on the
interval $I$ if
the derivative $F'(x)$ exists for almost all $x\in I,$ and it is
locally Bochner summable on $I,$ and
\begin{equation*}
 F(x_1)-F(x_2)=\int_{x_2}^{x_1}F'(t)\,dt\fa x_1,x_2\in I.
\end{equation*}
\end{defin}

\bigskip

\begin{thm}[Newton's formula and Lipschitzian functions]
Let $(I,V,v)$ be the Riemann measure space over an interval
$I.$ Assume that $F$ is a Lipschitzian function
from the interval $I$ into a Banach space $Y.$

If the derivative $F'(x)$ exists for almost all $x\in I$ then
the derivative is locally summable on $I$ and Newton's formula
holds
\begin{equation*}
 F(x_1)-F(x_2)=\int_{x_2}^{x_1}F'(t)\,dt\fa x_1,x_2\in I.
\end{equation*}
\end{thm}

\bigskip

\bp
    In the case when the interval $I$ is bounded on the right, notice that
    the Lipschitz condition implies the existence of the limit
    $$\lim_{x\to b}F(x)=y_0$$ where $b$ is the right end of the interval $I.$
    In such a case extend the function $F$ to the right by the formula
\begin{equation*}
     F(x)=y_0\fa x\ge b.
\end{equation*}
    This operation will not change the Lipschitzian property nor the
    differentiability almost everywhere of the function $F.$

    Take any sequence $h_n$ of positive numbers converging to zero
    and consider the sequence of functions
\begin{equation*}
     f_n(x)=\frac{1}{h_n}(F(x+h_n)-F(x))\fa x\in I.
\end{equation*}
    Notice that the functions $f_n$ are well defined and are continuous
    on the interval $I$ and as such they are locally Bochner
    summable on $I.$ If $M$ denotes the Lipschitz constant of $F$
    then
\begin{equation*}
     |f_n(x)|=\frac{1}{h_n}|F(x+h_n)-F(x)|\le M\fa x\in I.
\end{equation*}
    Moreover by assumption we have $\lim_n f_n(x)=F'(x)$ for almost all $x\in I.$
    Thus by Lebesgue's Dominated Convergence Theorem the function $F'$
    is locally summable on the interval $I$ and
\begin{equation*}
     \int_{x_1}^{x_2}F'(t)\,dt=
     \lim_n\int_{x_1}^{x_2}f_n(t)\,dt\fa x_1\le x_2,\,x_1,x_2\in I.
\end{equation*}

    From linearity of the integral and the fact that
    with respect to Riemann measure the integral is invariant
    under translations, we have
\begin{equation*}
 \begin{split}
    \int_{x_1}^{x_2}f_n(t)\,dt
    &=\frac{1}{h_n}\left[\int_{x_1}^{x_2}F(t+h_n)\,dt-\int_{x_1}^{x_2}F(t)\,dt\right]\\
    &=\frac{1}{h_n}\left[\int_{x_1+h_n}^{x_2+h_n}F(t)\,dt-\int_{x_1}^{x_2}F(t)\,dt\right]\\
    &=\frac{1}{h_n}[(\int_{x_1+h_n}^{x_2}F(t)\,dt+\int_{x_2}^{x_2+h_n}F(t)\,dt)\\
    &\qquad-(\int_{x_1}^{x_1+h_n}F(t)\,dt+\int_{x_1+h_n}^{x_2}F(t)\,dt)]\\
    &=\frac{1}{h_n}\left[\int_{x_2}^{x_2+h_n}F(t)\,dt
    -\int_{x_1}^{x_1+h_n}F(t)\,dt\right]\\
\end{split}
\end{equation*}

    By assumption the function $F$ is Lipschitzian, so it is continuous.
    Passing to the limit
     $n\to\infty$ in the above equality we get
\begin{equation*}
     \int_{x_1}^{x_2}F'(t)\,dt=F(x_2)-F(x_1)\fa x_1\le x_2,\,x_1,x_2\in I.
\end{equation*}
    For the case when the integration limits are in reversed order $x_2< x_1$,
    from the above formula and the definition of the integral with respect
    to Riemann measure follows
    that Newton's formula holds for all $x_1,x_2\in I.$
\ep

\bigskip

Newton's formula in the case when the derivative $F'$ is continuous is known
as the {\bf Fundamental Theorem of Calculus.} The existence of the derivative
$F'$ even at every point of the interval does not guarantee the validity
of Newton's formula.

Consider for instance the real-valued function defined by the formulas
\begin{equation}\label{non-absolutely continuous function}
 F(x)=x^2\sin(1/x^2)\fa x\not=0\qtext{and}F(0)=0.
\end{equation}

The derivative $F'$ exists at every point $x\in R$ but it is not summable on
any interval containing in its interior the point $x=0.$
One can prove this fact by showing that the function $F$ is not {\em absolutely
continuous,} which is a necessary condition for a function to be
an indefinite integral of a locally summable function. So let us
consider this notion.

\bigskip

\begin{defin}[Absolute continuity]
Assume that $(X,V,v)$ is a Riemann measure space
over an interval $X$ and $Y$ is a Banach space.

A function $F:X\mto Y$ is said to be {\bf absolutely continuous} on $X$
if for every $\varepsilon>0$ there exists a $\delta>0$ such that
for any finite system of disjoint intervals
\begin{equation*}
 I_k=(a_k,b_k)\subset X\qtext{where}k=1,2,\ldots,n
\end{equation*}
we have
\begin{equation*}
 \sum_k |F(b_k)-F(a_k)|<\varepsilon\qtext{if}\sum_k v(I_k)<\delta.
\end{equation*}
\end{defin}

\bigskip

Notice that if $F$ is a Lipschitzian function with constant $M$
then the formula for the number $\delta$ corresponding to $\varepsilon>0 $
in the above definition is
\begin{equation*}
 \delta=\varepsilon/M.
\end{equation*}
Thus every Lipschitzian function is absolutely continuous.
We have the following general theorem.

\bigskip

\begin{thm}[Absolute continuity of indefinite integral]
Assume that $(X,V,v)$ is a Riemann measure space
over an interval $X$ and $Y$ is a Banach space.

If $f\in L(v,Y)$ then its indefinite integral given by
\begin{equation*}
 F(x)=\int_a^x f(t)\, dt\fa x\in X\qtext{and some}a\in X
\end{equation*}
is absolutely continuous.
\end{thm}

\bigskip
\bp
    Take any $\varepsilon>0.$ By Lemma \ref{Density of simple functions in $L(v,Y)$}
    on density of simple functions in the space $L(v,Y)$ there exists a function
    $s\in S(V,Y)$ such that
\begin{equation*}
     \norm{f-s}<\varepsilon/2.
\end{equation*}
    Let $S$ denote the indefinite integral of the function $s$
\begin{equation*}
    S(x)=\int_a^x s(t)\, dt\fa x\in X
\end{equation*}
    and $G$ the indefinite integral of the function $g=f-s.$

    Let $M=\sup\set{|s(x)|:\ x\in X}.$ Since $M$ represents a Lipschitz
    constant of $S$ we have for any system of disjoint intervals
\begin{equation*}
    I_k=(a_k,b_k)\qtext{where}k=1,2,\ldots,n
\end{equation*}
    and $\delta=\varepsilon/(2M)$ that
\begin{equation*}
    \sum_k |S(b_k)-S(a_k)|<\varepsilon/2\qtext{if}\sum_k v(I_k)<\delta.
\end{equation*}
    Put $A=\bigcup_kI_k$ and notice the estimate
\begin{equation*}
    \sum_k |G(b_k)-G(a_k)|\le\sum_k \int_{a_k}^{b_k}|g(t)|\,dt
    \le\int_A|g(t)|\,dt\le \norm{g}<\varepsilon/2.
\end{equation*}
    Thus we have
\begin{equation*}
     \sum_k |F(b_k)-F(a_k)|\le \sum_k |G(b_k)-G(a_k)|+\sum_k |S(b_k)-S(a_k)|<\varepsilon.
\end{equation*}
\ep

\bigskip

We leave it to the reader to show that the function $F$ given by
the formula (\ref{non-absolutely continuous function}) is not absolutely continuous
in any interval containing the point $x=0$ and thus it is not
an indefinite integral of a  Bochner summable function in any such interval.

\bigskip

Our goal is to prove that any function $F$  from an interval $X$ into a Banach space $Y$
possessing almost everywhere on $X$ a locally summable derivative $F'$ satisfies
Newton's formula
\begin{equation*}
 F(y)-F(x)=\int_x^y F'(t)\,dt\fa x,y\in X
\end{equation*}
if and only if the function $F$ is absolutely continuous. To accomplish this we will
need several theorems of Lebesgue concerning this topic for real-valued functions.

\bigskip

\section{Lebesgue theory of absolutely continuous functions}

\bigskip

We start this section by proving several important lemmas and propositions
playing a key role in the development of the Lebesgue theory of
absolutely continuous functions.

\bigskip

\begin{pro}
Let $(X,V,v)$ be any measure space over a prering $V$ of an abstract space $X.$

If for some summable function $f\in L(v,R)$ we have
\begin{equation*}
 \int_A f(t)\,dt=0\fa A\in V
\end{equation*}
then $f(x)=0$ for almost all $x\in X.$
\end{pro}

\bigskip
\bp
    By assumption of the theorem for any set $A\in V$
    we have
    \begin{equation*}
         \int c_Af\,dv=0.
    \end{equation*}
    As a consequence for any simple function $s$ we must have
    \begin{equation*}
         \int sf\,dv=0\fa s\in S(V,R).
    \end{equation*}
    Since simple functions are dense in the Lebesgue space $L(v,R),$
    for any $\varepsilon>0$ there is a simple function $s_0\in S(V,R)$
    such that
    \begin{equation*}
         \norm{s_0-f}<\varepsilon.
    \end{equation*}
    If the simple function $s_0$ has representation
    \begin{equation*}
         s_0=\sum_k r_k c_{A_k},
    \end{equation*}
    where $\set{A_k}$ is a finite family of disjoint sets from the prering $V,$
    define function $s$ by the formula
    \begin{equation*}
         s=\sum_k \mathrm{sign}(r_k) c_{A_k},
    \end{equation*}
    where
    \begin{equation*}
         \mathrm{sign}(r)=r/|r|\qtext{ if }r\not=0\qtext{and}\mathrm{sign}(r)=0\qtext{ if }r=0.
    \end{equation*}
    By linearity and monotonicity of the integral we have
    \begin{equation*}
        \begin{split}
         \norm{s_0}&=\int |s_0|\,dv=\int ss_0\,dv=\int s(s_0-f)\,dv+\int sf\,dv\\
         &\le \int |s_0-f|\,dv+0=\norm{s_0-f}<\varepsilon.
        \end{split}
    \end{equation*}
    Thus from the triangle inequality for a norm we get
    \begin{equation*}
         \norm{f}\le\norm{s_0}+\norm{s_0-f}<2\varepsilon.
    \end{equation*}
    Since $\varepsilon$ was fixed but arbitrary we must have $f(x)=0$
    for almost all $x\in X.$
\ep


\bigskip

The above proposition can be immediately generalized to the case of Bochner summable
functions but first let us prove a lemma that is of interest by itself.

\bigskip

\subsection{Commutativity of a linear bounded operator with integral}

\bigskip

\begin{lemma}
\label{commutativity of int with lin op}
Let $(X,V,v)$ be any measure space over a prering $V$ of an abstract space $X.$
Assume that $Y,\,Z$ are  Banach spaces and $u$ a linear bounded operator from $Y$ to $Z.$

If $f\in L(v,Y)$ then $u\circ f\in L(v,Z)$
and
\begin{equation*}
 u(\int f\, dv)=\int u\circ f\, dv
\end{equation*}
\end{lemma}

\bigskip

\bp
    By definition of a Bochner summable function $f\in L(v,Y)$ there exists
    a basic sequence $s_n$ and a null set $D$ such that
\begin{equation*}
     s_n(x)\to f(x)\fa x\in X\less D.
\end{equation*}
    Since the operator $u$ is linear and bounded, the composition
    $S_n(x)=u(s_n(x))$ yields a basic sequence converging to the
    function $u(f(x))$ for every $x\in X\less D.$
    Thus we have $u\circ f\in L(v,Z)$ and by linearity and continuity of the operator $u$
    we get
\begin{equation*}
     \int u\circ f\,dv=\lim_n\int u\circ s_n\, dv=u(\,\lim_n\int s_n\, dv)=u(\int f\, dv).
\end{equation*}
\ep

\bigskip

\begin{pro}
Let $(X,V,v)$ be any measure space over a prering $V$ of an abstract space $X$
and let $Y$ be a Banach space.

If for some Bochner summable function $f\in L(v,Y)$ we have
\begin{equation*}
 \int_A f(t)\,dt=0\fa A\in V
\end{equation*}
then $f(x)=0$ for almost all $x\in X.$
\end{pro}

\bigskip

\bp
    Let $Y'$ denote the dual Banach space of the space $Y'.$
    It follows from the assumption of the theorem that for any set $A\in V$
    and any functional $y'\in Y'$
    we have
    \begin{equation*}
         0=\int y'\circ (c_Af)\,dv=\int c_Ay'\circ (f)\,dv.
    \end{equation*}
    As a consequence for any simple function $s\in S(V,Y')$ we must have
    \begin{equation*}
         \int u(s,f)\,dv=0\fa s\in S(V,Y'),
    \end{equation*}
    where $u$ denotes the bilinear bounded operator
    \begin{equation*}
         u(y',y)=y'(y)\fa y'\in Y'\qtext{and}y\in Y.
    \end{equation*}
    Since simple functions are dense in the Lebesgue space $L(v,Y),$
    for any $\varepsilon>0$ there is a simple function $s_0\in S(V,Y)$
    such that
    \begin{equation*}
         \norm{s_0-f}<\varepsilon.
    \end{equation*}
    If the simple function $s_0$ has representation
    \begin{equation*}
         s_0=\sum_k y_k c_{A_k},
    \end{equation*}
    where $\set{A_k}$ is a finite family of disjoint sets from the prering $V,$
    define function $s$ by the formula
    \begin{equation*}
         s=\sum_k y'_k c_{A_k},
    \end{equation*}
    where $y'_k\in Y'$ is linear functional such that
    \begin{equation*}
         y'_k(y_k)=|y_k|\qtext{and}|y'_k|=1..
    \end{equation*}
    By linearity and monotonicity of the integral we have
    \begin{equation*}
        \begin{split}
         \norm{s_0}&=\int |s_0|\,dv=\int u(s,s_0)\,dv=\int u(s,(s_0-f))\,dv+\int u(s,f)\,dv\\
         &\le \int |s_0-f|\,dv+0=\norm{s_0-f}<\varepsilon.
        \end{split}
    \end{equation*}
    Thus from the triangle inequality for a norm we get
    \begin{equation*}
         \norm{f}\le\norm{s_0}+\norm{s_0-f}<2\varepsilon.
    \end{equation*}
    Since $\varepsilon$ was fixed but arbitrary we must have $f(x)=0$
    for almost all $x\in X.$
\ep


\bigskip

\begin{lemma}
Let $(X,V,v)$ be the Riemann measure space over a closed bounded interval $X=[a,b].$

If for some Lebesgue summable function $f\in L(v,R)$ we have
\begin{equation*}
 \int_a^x f(t)\,dt=0\fa x\in X
\end{equation*}
then $f(x)=0$ for almost all $x\in X.$
\end{lemma}

\bigskip

\bp
    It follows from the assumption of the theorem
    and linearity of the integral, that for any interval $A\in V$
    we must have
    \begin{equation*}
         \int c_A(t)f(t)\,dt=0\fa A\in V.
    \end{equation*}
    Thus from the preceding proposition we get $f(x)=0$
    for almost all $x\in X.$
\ep

\bigskip

\begin{lemma}
Let $(X,V,v)$ be the Riemann measure space over a closed bounded interval $X=[a,b].$

If $f$ is a nonnegative bounded function such that $f\in L(v,R)$ then its
indefinite integral $F$ given by
\begin{equation*}
 F(x)=\int_a^x f(t)\,dt\fa x\in X
\end{equation*}
is differentiable almost everywhere on $X$ and
\begin{equation*}
 F'(x)=f(x)\qtext{for almost all}x\in X.
\end{equation*}
\end{lemma}

\bigskip

\bp
    Since the function $F$ is increasing on the interval $X$ it is differentiable
    almost everywhere. Denote its derivative by $g.$
    Thus $F'(x)=g(x)$ for almost all $x\in X.$

    Since the function $f$ by assumption is
    bounded its indefinite integral $F$ is Lipschitzian.
    Since we already established that for such a function $F$ the Newton formula
    is valid, we have
    \begin{equation*}
        F(x)-F(a)=\int_a^x F'(t)\,dt=\int_a^x g(t)\,dt\fa x\in X.
    \end{equation*}
        Since $F(a)=0$ we must have
    \begin{equation*}
        \int_a^x g(t)\,dt=F(x)=\int_a^x f(t)\,dt\fa x\in X.
    \end{equation*}
        That is
    \begin{equation*}
        \int_a^x (g(t)-f(t))\,dt=0\fa x\in X.
    \end{equation*}
    By preceding Lemma and the above we must have $(g(x)-f(x))=0$ almost
    everywhere and therefore
    \begin{equation*}
        F'(x)=g(x)=f(x)\qtext{for almost all}x\in X.
    \end{equation*}
\ep

\bigskip


The following proposition is just a stepping stone
to a more general result that will follow.

\bigskip

\begin{pro}[Lebesgue]
\label{pro-diff of real fun}
Let $(X,V,v)$ be the Riemann measure space over a closed bounded interval $X=[a,b].$

If $f$ is a nonnegative Lebesgue summable function $f\in L(v,R)$ then its
indefinite integral $F$ given by
\begin{equation*}
 F(x)=\int_a^x f(t)\,dt\fa x\in X
\end{equation*}
is differentiable almost everywhere on $X$ and
\begin{equation*}
 F'(x)=f(x)\qtext{for almost all}x\in X.
\end{equation*}
\end{pro}

\bigskip

\bp
    Since $X\in V$ the functions of the form $nc_X$ are simple
    and as such they are summable. Since the space $L(v,R)$
    forms a linear lattice we must have
    \begin{equation*}
        f_n=f\wedge nc_X\in L(v,R)\fa n=1,2,\ldots
    \end{equation*}
    Notice that the sequence $f_n$ increasingly converges everywhere on the set $X$ to the
    function~$f.$

    Introduce the indefinite integrals
    \begin{equation*}
        F_n(x)=\int_a^x f_n(t)\,dt\qtext{ and }F(x)=\int_a^x f(t)\,dt\fa x\in X.
    \end{equation*}
    Since the above indefinite integrals form increasing functions, their derivatives
    $F'_n(x)$ and $F'(x)$ exist for almost all $x\in X.$

    Since each function $f_n$ is summable
    and bounded, from the preceding lemma we have $F'_n(x)=f_n(x)$ for almost all $x\in X.$

    It follows from the Monotone Convergence Theorem,
    that the sequence $F_n$ increasingly converges
    everywhere on $X$ to the function $F$. Thus by Corollary to Fubini theorem
    \ref{cor-Fubini for increasing sequence} we have that the derivative
    $F'$ exists almost everywhere on $X$ and
    \begin{equation*}
        f_n(x)=F'_n(x)\to F'(x)\qtext{for almost all}x\in X.
    \end{equation*}

    Thus by Lebesgue's Dominated Convergence Theorem we have
    \begin{equation*}
        \int_a^x f_n(t)\, dt\to \int_a^x F'(t)\, dt\fa x\in X
    \end{equation*}
    and also
    \begin{equation*}
        \int_a^x f_n(t)\, dt\to \int_a^x f(t)\, dt\fa x\in X.
    \end{equation*}
    From the uniqueness of the limit we must have
    \begin{equation*}
        \int_a^x F'(t)\, dt=\int_a^x f(t)\, dt\fa x\in X,
    \end{equation*}
    or equivalently
    \begin{equation*}
        \int_a^x (F'(t)-f(t))\, dt=0\fa x\in X,
    \end{equation*}
    which implies $F'(x)=f(x)$ for almost all $x\in X.$
\ep

\bigskip

\begin{defin}[Lebesgue points of a summable function]
Let $(X,V,v)$ be any Riemann measure space over an interval $X$
and $Y$ a Banach space.
Assume that $f$ is a locally Bochner summable function from $X$ into $Y.$

A point $p\in X$ is called a {\bf Lebesgue point} of the function $f$ if
\begin{equation*}
 \lim_{h\to  0}\frac{1}{h}\int_p^{p+h}|f(t)-f(p\,)|\,dt=0.
\end{equation*}
\end{defin}

\bigskip

The proof of the following theorem is due to Lebesgue who proved it
for the case of real-valued functions. It carries over to Bochner summable
functions with only minor modifications.

\bigskip

\begin{thm}
Let $(X,V,v)$ be any Riemann measure space over an interval $X$
and $Y$ a Banach space.

If $f$ is a locally Bochner summable function from $X$ into $Y,$
then almost every point of $X$ is a Lebesgue point of the function $f.$
\end{thm}

\bigskip
\bp
    Without loss of generality we may assume that the set $X$ forms a closed
    bounded interval.
    Thus the function $f$ is summable on $X.$ Therefore there exists a basic
    sequence $s_n\in S(V,Y)$ of simple functions and a null set $D_0$ such that
    \begin{equation*}
         s_n(x)\to f(x)\fa x\in X\less D_0.
    \end{equation*}
    Notice that the image set $s_n(X)$ is finite for each $n$ and thus
    the set $B=\bigcup_n s_n(X)$ is at most countable.

    Represent the set $B$ as a sequence
    \begin{equation*}
         B=\set{y_1,\,y_2,\,\ldots}.
    \end{equation*}
    It follows from the triangle inequality for the norm that for any index $m$ the
    map $y\mto |y-y_m|$ is Lipschitzian and as such it maps
    a basic sequence into a basic sequence.
    Thus the sequence $|s_n(x)-y_m|$ for fixed $m$ is basic and it converges
    almost everywhere on $X$ to the function $|f(x)-y_m|.$ Therefore each such function
    belongs to the space $L(v,R).$

    By Lebesgue's Proposition \ref{pro-diff of real fun}
    there exist null sets $D_m$ such that
    \begin{equation*}
         |f(p)-y_m|=\lim_{h\to 0}\frac{1}{h}\int_p^{p+h}|f(t)-y_m|\,dt\fa p\in X\less D_m.
    \end{equation*}
    Put $D=\bigcup_{m\ge 0} D_m.$ For any $p\in X\less D$ and  any index $m$ we have
    \begin{equation*}
         \begin{split}
        &\limsup_{h\to0}\frac{1}{h}\int_p^{p+h}|f(t)-f(p)|\,dt\\
        &\le\limsup_{h\to0}\frac{1}{h}\int_p^{p+h}|f(t)-y_m|\,dt+
        \limsup_{h\to0}\frac{1}{h}\int_p^{p+h}|y_m-f(p)|\,dt\\
        &= |f(p)-y_m|+|y_m-f(p)|=2|f(p)-y_m|.
        \end{split}
    \end{equation*}
    Taking any $\varepsilon>0$ and choosing $m$ so that $|f(p)-y_m|<\varepsilon/2$
    we get
    \begin{equation*}
        \limsup_{h\to0}\frac{1}{h}\int_p^{p+h}|f(t)-f(p)|\,dt<\varepsilon
    \end{equation*}
    that is
    \begin{equation*}
        \lim_{h\to  0}\frac{1}{h}\int_p^{p+h}|f(t)-f(p)|\,dt=0.
    \end{equation*}
    Hence the proof is complete.
\ep

\bigskip

The following theorem is a simple consequence of the preceding one.
\bigskip

\begin{thm}[Differentiability a.e. of indefinite integral]
Let $(X,V,v)$ be any Riemann measure space over an interval $X$
and $Y$ a Banach space.

If $f$ is a locally Bochner summable function from $X$ into $Y,$
then its indefinite integral $F$ defined by the formula
\begin{equation*}
 F(x)=\int_a^x f(t)\,dt\fa x\in X\qtext{and some}a\in X
\end{equation*}
is differentiable almost everywhere on $X$ and we have the equality
\begin{equation*}
 F'(x)=f(x)\qtext{for almost all}x\in X.
\end{equation*}
\end{thm}

\bigskip
\bp
    Let $p\in X$ be a Lebesgue point of the function $f.$  Take any
    number $h\not=0$ such that $p+h\in X$ and form the difference
    quotient of the function $F.$ We have the following estimate
\begin{equation*}
 \begin{split}
    \left|\frac{F(p+h)-F(p)}{h}-f(p\,)\right|&=\left|\frac{1}{h}\int_p^{p+h}f(t)\,dt-f(p\,)\right|\\
    &=\left|\frac{1}{h}\int_p^{p+h}[f(t)-f(p\,)]\,dt\right|\\
    &\le \frac{1}{h}\int_p^{p+h}|f(t)-f(p\,)|\,dt.
\end{split}
\end{equation*}
    The last expression in the above
    estimates by definition of a Lebesgue point
    converges to zero when $h$ tends to zero.
    Thus $F'(p)=f(p)$ for almost all $p\in X.$
\ep

\bigskip

\subsection{Absolutely continuous function with zero derivative}

\bigskip

\begin{thm}
\label{absol cont fun with deriv zero}
Let $(X,V,v)$ be a Riemann measure space over a closed bounded interval $X=[a,b]$
on the prering $V$ of all subintervals of $X$ and $Y$ a Banach space.

If $F$ is an absolutely continuous function from $X$ into $Y$ and its
derivative $F'(x)$ exist for almost all $x\in X$ and is equal to zero
almost everywhere on $X$
then the function $F$ is constant on the interval $X.$
\end{thm}
\bigskip
\bp
    Define a set function $\eta$ on the prering $V$ by the formula
    \begin{equation*}
     \eta(A)=F(\beta)-F(\alpha)\fa A\in V,
    \end{equation*}
    where $\alpha\le \beta$ are the end points of the interval $A.$
    Notice that the set function $\eta$ is finitely additive on the prering $V$
    and as such it can be uniquely extended onto the ring $S(V)$ of simple
    sets generated by $V.$

    We shall prove that $F(a)=F(c)$ for any point $c\in [a,b].$
    To this end consider the set
    \begin{equation*}
    E=\set{x\in [a,c]:\ F'(x)=0}.
    \end{equation*}
    Take any $\varepsilon>0$ and introduce the family of closed intervals
    \begin{equation*}
        H=\set{I\subset [a,c]:\ I=\emptyset\text{ or }
        I\cap E\not=\emptyset,\, 0<v(I),\,|\eta(I)|<\varepsilon v(I)}
    \end{equation*}

    The family $H$ forms a Vitali covering of the set $E.$
    Thus there exists a sequence of disjoint intervals $I_j\in H$ forming
    an interlaced cover of the set $E$ that is
    \begin{equation*}
        E\subset [\bigcup_{j\le n}I_j]\cup[\bigcup_{j>n}S(I_j)]\fa n.
    \end{equation*}
    Since the set $[a,c]\less E$ is a null set, there exists a sequence
    of intervals $J_j$ such that $\sum_j v(J_j)<1$ and
    \begin{equation*}
        [a,c]\less E\subset \bigcup_{j>n}J_j\fa n.
    \end{equation*}
    Thus we have the relation
    \begin{equation*}
        [a,c]\subset [\bigcup_{j\le n}I_j]\cup[\bigcup_{j>n}S(I_j)]
        \cup[\bigcup_{j>n}J_j]\fa n.
    \end{equation*}
    The above yields
    \begin{equation*}
        A_n=[a,c]\less  [\bigcup_{j\le n}I_j]\subset[\bigcup_{j>n}S(I_j)]
        \cup[\bigcup_{j>n}J_j]\fa n.
    \end{equation*}
    Notice that each set $A_n$ belongs to the ring $S(V)$ of simple sets
    and from countable subadditivity of the measure $v_c$ we have
    \begin{equation*}
        \begin{split}
        v_c(A_n)&\le \sum_{j>n}v_c(S(I_j))+\sum_{j>n}v_c(J_j)\\
        &\le 5\sum_{j>n}v(I_j)+\sum_{j>n}v(J_j)\fa n
        \end{split}
    \end{equation*}
    implying $v_c(A_n)\to 0.$
    So for sufficiently large $n$ we must have $v_c(A_n)<\delta,$
    where $\delta$ is selected by absolute continuity of the function $F$
    so that for every disjoint finite collection
    of intervals $B_k\subset [a,c]$ we have
    \begin{equation*}
        \sum_k |\eta(B_k)|<\varepsilon\qtext{if}\sum_k v(B_k)<\delta.
    \end{equation*}
    Let $B_k\in V\, (k=1,2,\ldots,m)$ be a finite decomposition of the set $A_n.$
    Then the family of intervals
    \begin{equation*}
        \set{B_k:\ k\le m}\cup\set{I_j:\ j\le n}
    \end{equation*}
    represents a finite disjoint decomposition of the interval $[a,c].$
    Thus by finite additivity of the set function $\eta$ we have
    \begin{equation*}
        \eta([a,c])=\sum_{k\le m} \eta(B_k)+\sum_{j\le n} \eta(I_j)
    \end{equation*}
    yielding
    \begin{equation*}
        \begin{split}
        |F(c)-F(a)|&=|\eta([a,c])|\\
        &\le\sum_{k\le m} |\eta(B_k)|+\varepsilon\sum_{j\le n} v(I_j)\\
        &\le \varepsilon+\varepsilon v([a,c])=\varepsilon(1+c-a).
        \end{split}
    \end{equation*}
    Since $\varepsilon$ was fixed but arbitrary we must have
    \begin{equation*}
    F(c)=F(a)\fa c\in [a,b].
\end{equation*}
\ep

\bigskip

\subsection{Newton's formula and absolute continuity}

\bigskip

\begin{thm}
Let $(X,V,v)$ be a Riemann measure space over a closed bounded interval $X=[a,b]$
on the prering $V$ of all subintervals of $X$ and $Y$ a Banach space.

Assume that a function $F$ from $X$ into $Y$ has
a derivative $F'(x)$ for almost all $x\in X$ and
the derivative is Bochner summable that is $F'\in L(v,Y).$
Then the function $F$ satisfies the Newton formula
\begin{equation*}
 \int_x^y F'(t)\,dt=F(y)-F(x)\fa x,y\in X
\end{equation*}
if and only if
the function $F$ is absolutely continuous on $X.$
\end{thm}

\bigskip

\bp
    To prove the necessity of the condition take any function $F$ and assume
    that it satisfies the Newton formula. Thus we have
    \begin{equation*}
    F(x)=F(a)+\int_a^x F'(t)\,dt\fa x\in X.
    \end{equation*}
    Since we proved before that the indefinite integral of a Bochner summable
    function is absolutely continuous, from the above formula we can deduce
     function $F$ must be absolutely continuous.
\bigskip

    Conversely, assuming that the function $F$ is absolutely continuous on $X$
    introduce function
    \begin{equation*}
        H(x)=F(x)-\int_a^x F'(t)\,dt\fa x\in X.
    \end{equation*}
    The function $H$ as a difference of two absolutely continuous functions
    is absolutely continuous on $X.$
    Notice that $H'(x)=F'(x)-F'(x)=0$ for almost all $x\in X.$
    Thus by Theorem \ref{absol cont fun with deriv zero}
    on absolutely continuous function whose derivative is equal to zero almost
    everywhere, we get that $H$ is a constant function on the interval $X.$
    Therefore
    \begin{equation*}
        H(x)=H(a)=F(a)\fa x\in X
    \end{equation*}
    which implies
    \begin{equation*}
        F(x)-F(a)=\int_a^x F'(t)\, dt\fa x\in X
    \end{equation*}
    and the above implies the Newton formula
    \begin{equation*}
        F(y)-F(x)=\int_x^y F'(t)\,dt\fa x,y\in X.
    \end{equation*}
\ep

\bigskip

\subsection{A  consequences of Lebesgue-Bochner theory}

\bigskip

In Calculus courses we define an {\bf antiderivative}
$F$ of a continuous function $f$ as any function such that $F'=f$ and we
prove the Newton formula which is also called the {\em Fundamental Theorem of Calculus.}

Let $I$ be any interval and $Y$ a Banach space.
Denote by $L$ the family of all functions $$f:I\mto Y$$ that are locally Bochner summable
on $I$ with respect to Riemann measure space $(I,V,v).$
This family contains all continuous functions from $I$ into the space $Y.$

The above theorem shows that we can extend the notion of an antiderivative of
a function $f\in L$
by considering functions $F$ that are absolutely continuous on every bounded
subinterval of the interval $I$ and such that $F'=f$ almost everywhere on the interval~$I.$

Such antiderivatives satisfy Newton's formula. Moreover every function $f\in L$ possesses
such an antiderivative and any two antiderivatives of the same function differ by
a constant. This is one of the important consequences of Lebesgue-Bochner theory
of summable functions.

\bigskip


\bigskip

\section{Integration by parts for Lebesgue-Bochner summable functions}

\bigskip

\begin{thm}[Integration by parts]
Let $X$ denote a closed bounded interval $[a,b]$
and assume that $(X,V,v)$ represents the Riemann measure space.

Let $Y$ be a Banach space and $Z$
the field of reals $R$ or the field $C$ of complex numbers considered as a Banach
space depending on whether the Banach space $Y$ is over the field $R$ or $C.$

If $f\in L(v,Z)$ and $g\in L(v,Y)$ and $F$ and $G$ are absolutely continuous
functions on $X$ such that $F'=f$ and $G'=g$ almost everywhere on $X$ then
the following formula known as {\bf integration by parts} is true
\begin{equation*}
 \int_a^b f(t)G(t)\,dt + \int_a^b F(t)g(t)\,dt=F(b)G(b)-F(a)G(a).
\end{equation*}
\end{thm}

\bigskip

\bp
    Since the functions $F$ and $G$ are absolutely continuous on the closed
    bounded interval, they are bounded. Let $\norm{F}$ and  $\norm{G}$ denote
    the supremum norm on $X.$

    Put $H=FG.$
    Take any finite family of disjoint intervals $I_j=[a_j,b_j].$ We have the estimate
\begin{equation*}
 \begin{split}
    |H(b_j)-H(a_j)|&=|F(b_j)G(b_j)-F(b_j)G(a_j)+F(b_j)G(a_j)-F(a_j)G(a_j)|\\
    &\le \norm{F}|G(b_j)-G(a_j)|+\norm{G}|F(b_j)-F(a_j)|.
\end{split}
\end{equation*}
    Thus from the above estimate and from
    absolute continuity of the functions $F$ and $G$ follows the
    absolute continuity of their product $H.$

    Notice that the function $H$
    has derivative almost everywhere on the interval $X$ and
\begin{equation*}
     H'(x)=f(x)G(x)+F(x)g(x).
\end{equation*}
    Using Theorem \ref{Summability of $u(f,g)$} we get that the function
    $H'$ is in $L(v,Y)$ and applying Newton's formula we get
\begin{equation*}
    \int_a^b f(t)G(t)\,dt + \int_a^b F(t)g(t)\,dt=F(b)G(b)-F(a)G(a).
\end{equation*}
\ep

\bigskip



\bigskip

\section{Tensor product of measure spaces}

\bigskip

Assume now that we have two measure spaces $(X_i,V_i,v_i)$ each defined
on the prering $V_i$ of an abstract space $X_i$ for $i=1,2.$
Consider the Cartesian product
$X_1\times X_2.$
On the tensor product $V_1\otimes V_2$ of the
families $V_i$  of sets
\bb
    V_1\otimes V_2=\set{A_1\times A_2:\ A_1\in V_1\text{ and }A_2\in V_2}
\ee
define the set function by
\begin{equation*}
 v_1\otimes v_2(A_1\times A_2)=v_1(A_1)v_2(A_2)\fa A_i\in V_i\,(i=1,2).
\end{equation*}

\bigskip

\begin{thm}[Tensor product of measure spaces forms measure space]
Assume that $(X_i,V_i,v_i)$ are measure spaces for $i=1,2.$
Let the triple $(X,V,v)$
consist of $X=X_1\times X_2,$ $V=V_1\otimes V_2,$ and
$v_1\otimes v_2(A_1\times A_2)=v_1(A_1)v_2(A_2)$ for all $A=A_1\times A_2\in V.$

Then  the triple
\begin{equation*}
 (X,V,v)=(X_1\times X_2,V_1\otimes V_2,v_1\otimes v_2)
\end{equation*}
forms a positive measure space.
\end{thm}

\bigskip

\bp
    As we established before the tensor product of prerings forms a prering.
    To prove  that the set function $v=v_1\otimes v_2$
    is $\sigma$-additive take any set $A\times B$
    in $V$ and let $A_n\times B_n\in V$ denote a sequence of disjoint sets whose
    union is the set $A\times B.$ Notice the identity

    \begin{equation}\label{tensor decomposition}
        c_A(x_1)\,c_B(x_2)=\sum_n c_{A_n}(x_1)\,c_{B_n}(x_2)\fa x_1\in X_1,\ x_2\in X_2.
    \end{equation}
    Fixing $x_2$ and integrating with respect to $v_1$
    both sides of the equation (\ref{tensor decomposition} ) on the basis
    of Beppo Levi's theorem for a series \ref{levi for series}, we get
    \begin{equation*}
        v_1(A)\,c_B(x_2)=\sum_n v_1({A_n})\,c_{B_n}(x_2)\fa x_2\in X_2.
    \end{equation*}
    Integrating the above term by term with respect to $v_2$
    and applying again Beppo Levi's theorem for a series \ref{levi for series} we get
    \begin{equation*}
        v_1(A)\,v_2(B)=\sum_n v_1({A_n})\,v_2({B_n})
    \end{equation*}
    that is the set function
        $$v(A\times B)=v_1\otimes v_2(A\times B)
        =v_1(A)v_2(B)\fa A\times B\in V_1\otimes V_2$$
    is $\sigma$-additive.
    Hence the triple
        $$(X,V,v)=(X_1\times X_2,V_1\otimes V_2,v_1\otimes v_2)$$
    forms a measure space.
\ep

\bigskip

The above theorem has an immediate generalization to any finite
number of measure spaces.

\bigskip

\begin{thm}[Tensor product of $n$ measure spaces]
Assume that $(X_i,V_i,v_i)$ are measures over abstract spaces
$X_i$ and $V_i$ are prerings for $i=1,\ldots,n.$ Let the triple
$(X,V,v)$  consist of $X=X_1\times\cdots\times X_n,$
$V=V_1\otimes\cdots\otimes V_n,$ and
    $$v(A)=v_1\otimes\cdots\otimes v_n(A)=v_1(A_1)\cdots v_n(A_n)$$
for all $A=A_1\times\cdots\times A_n\in V.$ Then the triple
$(X,V,v)$ forms a measure space.
\end{thm}

\bigskip



\bigskip

\section{Fubini's theorems for Lebesgue and Bochner integrals}

\bigskip

Let $Y$ be a Banach space with norm $\abs{\ \,}.$
Assume that $(X_i,V_i,v_i)$ for
$i=1,2$ are measure spaces each on its prering $V_i$ of the abstract space $X_i.$
For the sake of brevity let $(X,V,v)$ denote
the product measure space $$(X_1\times X_2,V_1\otimes V_2,v_1\otimes v_2).$$

\bigskip

We shall adopt the following convention for a function $f$ of the point $(x_1,x_2)$
from the product space $X_1\times X_2.$ When $x_1$ is fixed, by the symbol $f(x_1,\cdot)$
we shall understand the function $x_2\mto f(x_1,x_2).$ To indicate the variable
of integration in expressions depending on the variable $x$ we
shall write
\begin{equation*}
 \int f\,dv=\int f(x)\,v(dx).
\end{equation*}
We shall also use abbreviation $v_1$-a.e. to indicate that the relation
preceding it holds almost everywhere with respect to measure $v_1.$

\begin{defin}[Space $Fub(Y)$] Denote by $Fub(Y)$ the set of
Bochner summable functions $f\in L(v,Y),$ for which the {\bf Fubini theorem}
is true, that is such that there exists a
$v_1$-null set $C$ and a summable function $h\in L(v_1,Y)$ such that
\begin{equation}\label{(a)}
    f(x_1,\cdot)\in L(v_2,Y)\qtext{ and }
    h(x_1)=\int f(x_1,x_2)\,v_2(dx_2)\qtext{ if }x_1\notin C,
\end{equation}
and moreover
\begin{equation}\label{(b)}
    \int f\, dv=\int h\  dv_1=\int\left(\int f(x_1,x_2)\,v_2(dx_2)\right)v_1(dx_1).
\end{equation}
\end{defin}

\bigskip

We shall prove that $$Fub(Y)=L(v_1\otimes v_2,Y)$$ using an argument
similar to the argument of S.~Saks \cite{saks} who used it to prove
the theorem for Lebesgue integrals generated by Lebesgue measures over sigma rings.
For a reference to the original theorem see Fubini \cite{fubini-int}.

\bigskip

From linearity of the integrals follows that the set $Fub(Y)$ is
linear.
Since for any fixed $y\in Y$ and a function of the form
\begin{equation*}
    s(x_1,x_2)=yc_{A_1\times A_2}(x_1,x_2)=yc_{A_1}(x_1)c_{A_2}(x_2)\fa
(x_1,x_2)\in X_1\times X_2
\end{equation*}
 we have
\begin{equation*}
    \int s\,dv=\int yc_{A_1}(x_1)c_{A_2}(x_2)\,dv=y\,v_1(A_1)v_2(A_2)=
    \int \left(\int yc_{A_2}\,dv_2\right) c_{A_1}\,dv_1,
\end{equation*}
every simple function $s$ is in the set $Fub(Y),$ that is we have
 $S(V,Y)\subset Fub(Y)$.

\bigskip

\begin{lemma}[$Fub(R)$ is closed under monotone pointwise convergence]
\label{fub-Lemma 1}
Assume that functions $f_n$ represent a sequence monotone
with respect to the relation less or equal everywhere on the Cartesian product
$X=X_1\times X_2.$

If the functions $f_n$ belong to the set $Fub(R),$ and
the sequence $f_n$ converges everywhere to a
finite-valued function $f,$ and the sequence of integrals $\int f_n\ dv$
is bounded, then the limit function $f$ also belongs to the set $Fub(R)$.
\end{lemma}

\bigskip

\bp
From the assumption of the lemma there exist $v_1$-null sets
$C_m$ and functions $h_m\in L(v_1,R)$ such that
$f_m(x_1,\cdot)\in L(v_2,R)$
 and $h_m(x_1)=\int f_m(x_1,\cdot)\,dv_2$ if $x_1\notin C_m$
and
\begin{equation*}
 \int h_m\ dv_1=\int f_m\ dv\qtext{for}m=1,2,\ldots
\end{equation*}

Put $C=\bigcup_mC_m$. Since $C$ is a $v_1$-null set, the sequence $h_m$ is monotone
with respect to the relation less or equal $v_1$-almost everywhere.

Therefore from the Monotone Convergence Theorem \ref{Monotone Convergence Theorem}
 we get that there exist a function $h\in L(v_1,R)$
and a $v_1$-null set $D$ such that
\begin{equation*}
    h_m(x_1)\to h(x_1)\qtext{ if } x_1\notin D
\end{equation*}
and $h_m$ converges to $h$ in the topology of the space $L(v_1,R)$. This implies
\begin{equation*}
    \int h_m\ dv_1\to\int h\ dv_1.
\end{equation*}

For any $x_1$ which is not
in the set $C\bigcup D$, the sequence of functions $f_m(x_1,\cdot)$,
converges monotonically at each point to
the function $f(x_1,\cdot)$. Therefore again, for a fixed $x_1\not\in C\bigcup D,$  from the
 Monotone Convergence Theorem \ref{Monotone Convergence Theorem},
 and the fact that the sequence
of integrals
\begin{equation*}
    \int f_m(x_1,\cdot)\,dv_2=h_m(x_1)
\end{equation*}
is bounded, since it is convergent,  we get
\begin{equation*}
    f(x_1,\cdot)\in L(v_2,R)
\end{equation*}
and the equality
\begin{equation*}
    h(x_1)=\lim_mh_m(x_1)=\lim_m\int f_m(x_1,\cdot)\,dv_2=
    \int f(x_1,\cdot)\,dv_2\qtext{ if }x_1\notin C\bigcup D.
\end{equation*}

Since by assumption the sequence $f_m$ converges monotonically everywhere to the function
$f$ and the sequence of integrals $\int  f_m\ dv$ is bounded,
therefore we have
\begin{equation*}
f\in L(v,R)\qtext{ and }\int f_m\,dv\to  \int f\,dv.
\end{equation*}
Thus we get
\begin{equation*}
\int h\ dv_1=\lim_m\int h_m\ dv_1=\lim_m\int f_m\ dv=\int f\ dv.
\end{equation*}
Hence $f\in Fub(R).$
\ep

\bigskip

\begin{defin}[Section of a set $A\subset X_1\times X_2$]
By a section of a set $A\subset X_1\times X_2$
corresponding to a fixed value of the variable $x_1\in X_1$ we shall
understand the set
\begin{equation*}
    A(x_1)=\set{x_2\in X_2:\ (x_1,x_2)\in A}.
\end{equation*}
Similarly we can define a section of the set for any fixed $x_2\in X_2.$
\end{defin}

\bigskip

\begin{lemma}[Almost all sections of a $v$-null set are $v_2$-null sets]
\label{fub-Lemma 2}
Assume that $(X_i,V_i,v_i)$ for
$i=1,2$ are measure spaces each over its prering $V_i$ of an abstract space $X_i.$
Let $(X,V,v)$ denote
the product measure space $$(X_1\times X_2,V_1\otimes V_2,v_1\otimes v_2).$$
If $A\subset X$ is a $v$-null set,
then there exists a $v_1$-null set $C\subset X_1$
such that for every $x_1\notin C$ the section
\begin{equation*}
    A(x_1)=\{x_2\in X_2\ :\ (x_1,x_2)\in A\}
\end{equation*}
of the set $A$ represents a $v_2$-null set.
\end{lemma}

\bigskip

\bp
Let $A$ be a $v$-null set. One can prove from the
definition of a null set, that there exists a sequence
of sets $A_n\in V$ such that
\begin{equation*}
    A\subset \bigcup_{n>m}A_n\qtext{ for }m=1,2,\ldots
\end{equation*}
and
\begin{equation*}
\sum^\infty_{n=1}v(A_n)\le 1.
\end{equation*}

Put
\begin{equation*}
B_{mn}=\bigcup_{m<j<m+n}A_j,\quad B_m=\bigcup_{n>m}A_n\qtext{ and }B=\bigcap_mB_m.
\end{equation*}
We see that $A\subset B$. Consider the characteristic function
$c_{B_m}$. Notice that since $V=V_1\otimes V_2$ is a prering
the family $S(V)$ of simple sets forms a ring. Since the set $B_{mn}$ is a finite
union of sets from the prering $V$ we must have $B_{mn}\in S(V).$
Hence the characteristic functions $c_{B_{mn}}$ are in the space $S(V,R)$ of simple functions.
Now when $n\to\infty$
the sequence converges increasingly everywhere to the function $c_{B_m}$.
Since for fixed $m$ the integrals of the functions $c_{B_{mn}}$ are bounded by the number
\begin{equation*}
\sum_{n>m}\ v(A_n)\qtext{for}m=1,2,\ldots
\end{equation*}
representing the remainder of a convergent series,
we get $c_{B_{m}}\in Fub(R)$ for all  $m$ and
\begin{equation*}
 \int c_{B_m}\ dv\le\sum_{n>m}\ v(A_n)\qtext{for}m=1,2,\ldots
\end{equation*}

When $m\to\infty$ the sequence $c_{B_m}$ converges monotonically everywhere to
the function $c_B$. Thus from Lemma \ref{fub-Lemma 1} we get $c_B\in Fub(R)$,
that is there
exists a function $h\in L(v_1,R)$ and a $v_1$-null set $C_1$ such that
$c_B(x_1,\cdot)\in L(v_2,R)$
\begin{equation*}
    h(x_1)=\int c_B(x_1,\cdot)\,dv_2\qtext{ if }x_1\notin C_1.
\end{equation*}
Moreover we have
\begin{equation*}
    \int h\ dv_1=\int c_B\, dv\leq\int c_{B_m}\ dv
    \leq\sum_{n>m}v(A_n)\qtext{ for all }m.
\end{equation*}

Since the function $c_B$ is non-negative, therefore we have
    $h(x_1)\geq0$ $v_1$-almost everywhere.
This yields
\begin{equation*}
\norm{h}=\int|h|\ dv_1=\int h\ dv_1.
\end{equation*}
Now from Theorem \ref{Basic properties of the space $L(v,Y)$},
concerning the basic properties of the space $L(v,Y)$,
we get that there exists a
$v_1$-null set $C_2$ such that $h(x_1)=0$ if $x_1\notin C_2$.

Put $C=C_1\bigcup C_2$. We see that
\begin{equation*}
0=\int c_B(x_1,\cdot)\ dv_2=\norm{c_B(x_1,\cdot)}_{v_2}\qtext{ if }x_1\notin C.
\end{equation*}
This yields
\begin{equation*}
c_B(x_1,x_2)=0\qtext{ if }x_2\notin C(x_1),
\end{equation*}
where $C(x_1)$ is a $v_2$-null set. Since we have the equality
\begin{equation*}
 c_B(x_1,x_2)=c_{B(x_1)}(x_2)\fa x_1\in X_1,\ x_2\in X_2,
\end{equation*}
where
\begin{equation*}
B(x_1)=\{x_2\in X_2\ :\ (x_1,x_2)\in B\},
\end{equation*}
we must have
\begin{equation*}
B(x_1)\subset C(x_1).
\end{equation*}
Thus the set $B(x_1)$ is a $v_2$-null set and so is the set
\begin{equation*}
A(x_1)\subset B(x_1).
\end{equation*}
Thus the Lemma is proved.
\ep

\bigskip

Using Lemma \ref{fub-Lemma 2} by a similar argument as in Lemma
\ref{fub-Lemma 1} we can prove the following:

\bigskip

\begin{lemma}[$Fub(R)$ is closed under monotone convergence $v$-a.e.]
\label{fub-Lemma 3}
Assume that functions $f_n$ represent a sequence monotone
with respect to the relation less or equal v-almost everywhere on the product space
$X=X_1\times X_2.$

If the functions $f_n$ belong to the set $Fub(R),$ and
the sequence $f_n$ converges v-almost everywhere to a
finite-valued function $f,$ and the sequence of integrals $\int f_n\ dv$
is bounded, then the limit function $f$ also belongs to the set $Fub(R)$.
\end{lemma}

\bigskip

\begin{thm}[Fubini-Bochner: $Fub(Y)=L(v,Y)$]\label{fub-Bochner}
Assume that $Y$ is a Banach space and $(X_i,V_i,v_i)$ for $i=1,2$ are measure spaces
over prerings $V_i.$
Let $(X,V,v)$ denote
the tensor product measure space $$(X_1\times X_2,V_1\otimes V_2,v_1\otimes v_2).$$
 Then every Bochner summable function $f\in L(v,Y)$ belongs to the family $Fub(Y),$
 that is
there exists a
$v_1$-null set $C$ and a summable function $h\in L(v_1,Y)$ such that
\begin{equation*}
    f(x_1,\cdot)\in L(v_2,Y)\qtext{ and }
    h(x_1)=\int f(x_1,x_2)\,v_2(dx_2)\qtext{ if }x_1\notin C,
\end{equation*}
and moreover
\begin{equation*}
    \int f\, dv=\int h\, dv_1=\int\left(\int f(x_1,x_2)\,v_2(dx_2)\right)\,v_1(dx_1).
\end{equation*}
\end{thm}

\bigskip

\bp
Take any Bochner summable function $f\in L(v,Y)$. It follows
from the definition of the space of summable functions that
there exists a basic sequence $s_n$ convergent $v$-almost everywhere
to the function $f$. By the definition of a basic sequence we have
\begin{equation*}
    s_n=h_1+\dots+h_n,\quad ||h_n||_v\leq M \ 4^{-n}\qtext{ for }n=1,2,\dots,
\end{equation*}
where $h_n\in S(V,Y)$ represents a sequence of simple functions.

Consider the sequence of real-valued functions
\begin{equation*}
    g_n(x)=|h_1(x)|+\cdots+|h_n(x)|\qtext{ for }x\in X.
\end{equation*}
We see that
\begin{equation*}
    g_n\in S(V,R)\qtext{ and }\int g_n\ dv\leq M\qtext{ for }n=1,2,\dots
\end{equation*}
Since the sequence $g_n$ is monotone and the sequence of integrals is bounded,
there exists a Lebesgue summable function
$g\in L(v,R)$ such that the sequence $g_n$ converges
$v$-almost everywhere to the function $g$. Since $g_n\in Fub(R)$
therefore $g\in Fub(R)$ in accord with Lemma \ref{fub-Lemma 3}.

Notice that
\begin{equation*}
    |s_n(x)|\leq g(x)\qtext{ for } n=1,2,\dots
\end{equation*}
$v$-almost everywhere. Therefore there exists a $v$-null set $A$
such that
\begin{equation*}
    |s_n(x)|\leq g(x)\qtext{for}x\notin A\qtext{and} n=1,2,\dots
\end{equation*}
and
\begin{equation*}
 s_n(x)\to f(x)\qtext{if}x\notin A.
\end{equation*}
Let $C_1$ be a $v_1$-null set such that the sections $A(x_1)$ of the set $A$ are
$v_2$-null sets if $x_1\notin C_1$.

Since $g\in Fub(R)$ there exists a $v_1$-null set $C_2$ and a function
$h\in L(v_1,R)$ such that
\begin{equation*}
   g(x_1,\cdot)\in L(v_2,R)\qtext{ and }h(x_1)
   =\int g(x_1,\cdot)\,dv_2\qtext{ if }x_1\not\in C_2
\end{equation*}
and
\begin{equation*}
       \int h\,dv_1=\int g\,dv.
\end{equation*}

Take any point
\begin{equation*}
   x_1\notin C_1\bigcup C_2=C.
\end{equation*}
We have
\begin{equation*}
   |s_n(x_1,x_2)| \leq g(x_1,x_2)\qtext{for}x_2\notin A(x_1)\qtext{and}n=1,2,\dots
\end{equation*}
Thus the sequence $s_n(x_1,\cdot)$ is dominated by a summable function
and converges $v_2$-almost everywhere to the function $f(x_1,\cdot)$.
Therefore from the
Dominated Convergence Theorem \ref{Dominated Convergence Theorem} we get
\begin{equation*}
   f(x_1,\cdot)\in L(v_2,R)
\end{equation*}
and
\begin{equation*}
   \tilde{s}_n(x_1)=\int s_n(x_1,\cdot)\,dv_2\to \int f(x_1,\cdot)\,dv_2.
\end{equation*}
Notice the estimate
\begin{equation*}
   |\tilde{s}_n(x_1)|\leq\int |s_n(x_1,\cdot)|\,dv_2
   \leq\int g(x_1,\cdot)\,dv_2=h(x_1)\qtext{ if }x_1\notin C.
\end{equation*}
The function defined by the formula
\begin{equation*}
   \tilde{f}(x_1)=\int f(x_1,\cdot)\,dv_2\qtext{if}x_1\notin C,
   \quad\tilde{f}(x_1)=0\qtext{if} x_1\in C,
\end{equation*}
being the limit almost everywhere of the sequence of simple
functions $\tilde{s}_n$, dominated by the summable function $h$,
is summable, that is $\tilde{f}\in L(v_1,Y)$. Moreover
\begin{equation*}
   \int\tilde{f}\,dv_1=\lim_n\int\tilde{s}_n\,dv_1=\lim_n\int s_n\,dv=\int f\,dv.
\end{equation*}
Thus we have proved that $L(v,Y)=Fub(Y).$
\ep

\bigskip

It is convenient to formulate the Fubini theorem in a more concise way.
Introduce first the following notion.

\bigskip

\begin{defin}[Meaningful iterated integral]
Assume that we are given two positive measure spaces $(X_i,V_i,v_i)$ on prerings $V_i $
of some abstract spaces $X_i$ and a Banach space  $Y.$
Let $(X,V,v)$ denote the tensor product of these measure spaces.

We shall say that the iterated integral for a function $f:X\mto Y$ is {\bf meaningful}
\begin{equation*}
 \int \left(\int f(x_1,x_2)\ v_2(dx_2)\right) v_1(dx_1)
\end{equation*}
if the inner integral exists as Bochner integral and
its value $h(x_1)$ for $v_1$-almost
all points $x_1$  yields a Bochner summable function $h\in L(v_1,Y).$
\end{defin}

\bigskip

\begin{thm}[Fubini theorem]
Assume that we are given two positive measure spaces $(X_i,V_i,v_i)$ on prerings $V_i $
of some abstract spaces $X_i$ and a Banach space  $Y.$
Let $(X,V,v)$ denote the tensor product of these measure spaces.

If $f$ belongs to the space $L(v,Y)$ of Bochner summable functions then the
following iterated integrals are meaningful and we have the equalities
\begin{equation*}
 \int f\ dv= \int \left(\int f(x_1,x_2)\ v_2(dx_2)\right) v_1(dx_1)=
 \int \left(\int f(x_1,x_2)\ v_1(dx_1)\right) v_2(dx_2)
\end{equation*}
\end{thm}

\bigskip


\bigskip

\section{Fubini theorem in terms of equivalence classes}

\bigskip

 As before assume that $X$ is an abstract space and $Y$  a Banach
space. We shall use the notation $\abs{\ \,}$ to denote the norm
in $Y.$ Assume that $(X,V,v)$ is a measure space over a prering
$V$ of sets of $X.$ The functional $f\mto \norm{f}$ considered on the
space $L(v,Y)$ represents a semi-norm. By identifying functions equal almost everywhere
we may obtain a Banach space.

To perform this identification formally proceed as follows.
Let $F(X,Y)$ denote the space of all functions from the
space $X$ into the space $Y.$

\bigskip

\begin{defin}[Equivalence relation generated by equality almost everywhere]
Define the relation $\equiv$ between functions in $F(X,Y)$ by the
condition
\begin{equation*}
    (f\equiv g)\iff (f(x)=g(x)\quad v-\text{almost everywhere}).
\end{equation*}

This relation is an {\bf equivalence relation} and thus it splits
the space $F(X,Y)$ into disjoint classes of functions that are
equal almost everywhere. Denote by $[f]$ the class of functions
containing the function $f.$ Now if $L\subset F(X,Y)$ is any
sub-collection of functions, let $L_0$ denote the collection of
classes generated by this equivalence relation, that is
\begin{equation*}
    L_0=\set{[f]:\  f\in L}.
\end{equation*}
\end{defin}

\bigskip

In particular $L_0(v,Y)$ will denote the space of classes that are
generated by the space $L(v,Y)$ of Bochner summable functions, and
$M_0(v,Y)$ the collection generated by the space $M(v,Y)$ of
Bochner measurable functions as defined in Bogdanowicz \cite{bogdan14}.

For $p>0$ let $L_0^p(v,Y)$ denote the collection generated by the
space $L^p(v,Y),$ where
\begin{equation*}
    L^p(v,Y)=\set{f\in M(v,Y):\ \int |f(x)|^p\,v(dx)<\infty}.
\end{equation*}

\bigskip

It follows from Theorem \ref{Theorem 1} that the functional
\begin{equation*}
    \norm{g}_v=\int |f|\,dv\qtext{if}f\in g\in L_0(v,Y).
\end{equation*} is
well defined and represents a norm on the space $L_0(v,Y)$.
Moreover the pair
\begin{equation*}
   (L_0(v,Y),\ \norm{\ }_v)
\end{equation*}
forms a Banach space.

Define the integral operator on the space $L_0(v,Y)$ by
\begin{equation*}
      \int g\ dv=\int f\ dv\qtext{ if }g\in L_0(v,Y)\qtext{ and }[f]=g.
\end{equation*}
From the linearity of the integral operator $\int f\ dv$ and from
the estimate
\begin{equation*}
      \left|\int f\, dv\right|\leq \norm{f}_v\qtext{ for all } f\in L(v,Y)
\end{equation*}
we get that the integral operator on the space $L_0(v,Y)$
is well defined. Indeed, if $[f_1]=[f_2]$ that is $
      f_1(x)=f_2(x)\quad v-a.e.
$
we have
\begin{equation*}
      \int f_1dv=\int f_2dv
\end{equation*}
according to Theorem \ref{Theorem 1}. The operator $\int g\, dv$
is linear on the space $
      L_0(v,Y)
$ and
\begin{equation*}
      \left|\int g\, dv\right|\leq \norm{g}_v\qtext{ for all } g\in L_0(v,Y).
\end{equation*}

It will be convenient to use also the following notation when it
is important to indicate the variable of integration
\begin{equation*}
      \int g(x)\,v(dx)=\int g\ dv.
\end{equation*}

\bigskip

Now consider the space $L_0(v,Y)$ for the product measure
$v=v_1\otimes v_2$. Take a function
\begin{equation*}
      f\in L_0(v,Y).
\end{equation*}
We shall say that the double integral \begin{equation*}
      I=\int\left(\int f(x_1,x_2)\ v_2(dx_2)\right)\ v_2(dx_2)
\end{equation*}
has a meaning if for every function
  $g\in L(v,Y)$  being a representative of the class $f$ there exists a $v_1$-null
set $C$ and a function
\begin{equation*}
      h\in L(v_1,Y)
\end{equation*}
such that the function $g(x_1,x_2)$ considered as the function of
the variable $x_2$ is $v_2$-summable if $x_1\notin C$ and
\begin{equation*}
      h(x_1)=\int g(x_1,x_2)\ v_2(dx_2)\quad (x_1\notin C),
\end{equation*}
and by the value of the double integral we shall understand the
value
\begin{equation*}
  I=\int h\ dv_1.
\end{equation*}

It follows from the next theorem that this definition is correct.
It will be convenient to use the following notation
\begin{equation*}
  \int f\ dv=\int f(x_1,x_2)\, v_1(dx_1)\, v_2(dx_2)
\end{equation*}
for an integral generated by the product measure.

\bigskip


\begin{thm}
Assume that $Y$ is a Banach space and $X_i$ are abstract spaces with
prerings $V_i.$
Assume that $(X_i,V_i,v_i)$ for $i=1,2$ form positive measure spaces.
Let $(X,V,v)$ denote the tensor product measure
space $$(X_1\times X_2,V_1\otimes V_2,v_1\otimes v_2).$$

Then for every function
\begin{equation*}
  f\in L_0(v,Y)
\end{equation*}
the iterated integrals in the following formula have a meaning and
they are equal
\begin{equation*}
\begin{split}
    \int f(x_1,x_2)\ v_1(dx_1)\, v_2(dx_2)
    &=\int\left(\int f(x_1,x_2)\, v_1(dx_1)\right)\, v_2(dx_2)\\
    &=\int\left(\int f(x_1,x_2)\, v_2(dx_2)\right)\, v_1(dx_1).
\end{split}
\end{equation*}
\end{thm}

\bigskip

\bp
    The proof of the theorem follows immediately from the above
    definitions and from Fubini-Bochner Theorem \ref{fub-Bochner}.
\ep

\bigskip

Now let us find the relations between the completions
$v_c,(v_1)_c,(v_2)_c$ of the measures $v=v_1\otimes v_2,\ v_1,\
v_2$. We remind the reader that if $v$ is a measure on a prering
$V$ of an abstract space $X$ then the completion $v_c$ is defined
on the family of summable sets
\begin{equation*}
   V_c=\{A\subset X :\ c_A\in L(v,R)\}
\end{equation*}
by the formula
\begin{equation*}
   v_c(A)=\int c_A\, dv=\norm{c_A}_v.
\end{equation*}

\bigskip

\begin{thm}
Assume that $X_i$ are abstract spaces with
prerings $V_i.$
Assume that $(X_i,V_i,v_i)$ for $i=1,2$ form positive measure spaces.
Let $(X,V,v)$ denote the tensor product measure
space $$(X_1\times X_2,V_1\otimes V_2,v_1\otimes v_2).$$

Then for every summable set $A\subset V_c$ there
exists $v_1$-null set $C$ such that
\begin{equation*}
       A(x_1)\in(V_2)_c\qtext{ if }x_1\notin C
\end{equation*}
and the function $h$ given by the formula
\begin{equation*}
       h(x_1)=(v_2)_c(A(x_1))\qtext{ if }x_1\notin C
\end{equation*}
belongs to the space $L_0(v_1,R)$ and
\begin{equation*}
       v_c(A)=\int (v_2)_c(A(x_1))\ dv_1.
\end{equation*}
\end{thm}

\bigskip

\bp
    Consider the function
    \begin{equation*}
        f(x_1,x_2)=c_A(x_1,x_2)\qtext{ for all }(x_1,x_2)\in X_1\times X_2.
    \end{equation*}
    Notice that
    \begin{equation*}
        f(x_1,x_2)=c_{A(x_1)}(x_2)\qtext{ for all }(x_1,x_2)\in X_1\times X_2.
    \end{equation*}
    Using Fubini-Bochner Theorem \ref{fub-Bochner}
    and the definition of the completion of a measure on a prering,
    we can conclude the proof.
\ep

\bigskip

If $V$ is a family of sets denote by $V^\sigma$ the family of all
sets of the form
\begin{equation*}
  A=\bigcup_nA_n
\end{equation*}
where $A_n$ is a sequence of sets from the family $V$.

\begin{thm}
\label{fub-Theorem 4}
[\{$A$ is $v_1$-null \& $B\in (V_2)^\sigma$\}$\impl$ $A\times B$ is $v$-null]
Assume that $X_i$ are abstract spaces with
prerings $V_i.$
Assume that $(X_i,V_i,v_i)$ for $i=1,2$ form positive measure spaces.
Let $(X,V,v)$ denote the tensor product measure
space $$(X_1\times X_2,V_1\otimes V_2,v_1\otimes v_2).$$

If the set $A$ is a
$v_1$-null set and $B\in(V_2)^\sigma$ then the set $A\times B$ is
a $v$-null set.
\end{thm}

\bigskip

\bp
    Take any sequence of sets $B_n\in V_2$ and let
    $B=\bigcup_nB_n$ denote their union. The set $A$ is a $v_1$-null
    set if and only if there exists a sequence of sets $A_n\in V_1$
    such that
    \begin{equation*}
        \sum^\infty_{n=1} v_1(A_n)<\infty
    \end{equation*}
    and
    \begin{equation*}
        A\subset\bigcup_{n>m}A_n\qtext{ for }m=1,2,\dots
    \end{equation*}

    Notice the equality
    \begin{equation*}
        A\times B=\bigcup_n A\times B_n
    \end{equation*}
    Consider a fixed set $A\times B_n$. We have
    \begin{equation*}
        A\times B_n\subset \bigcup_{k>m}A_k\times B_n\qtext{ for }m=1,2,\dots.
    \end{equation*}
    and
    \begin{equation*}
        \sum_{k>m}v(A_k\times B_n)
        =\sum_{k>m}v_1(A_k)v_2(B_n)
        =v_2(B_n)\sum_{k>m}v_1(A_k)\to 0
    \end{equation*}
    as $m\to\infty.$ This implies that the set $A\times B_n$ is a $v$-null set. Thus
    the set $A\times B$ as the union of a countable number of $v$-null
    sets is a $v$-null set.
\ep

\bigskip

\begin{defin}[Family $N(v,Y)$ of null functions]
If $Y$ is a Banach space and $(X,V,v)$ is a measure space over a
prering $V$ let
\begin{equation*}
    N(v,Y)=\set{f\in L(v,Y):\ \int |f(x)|\,v(dx)=0}.
\end{equation*}
Functions belonging to this family will be called {\bf null
functions.}
\end{defin}

\bigskip

\begin{cor}
\label{fub-Corolary 1}
Assume that $X_i$ are abstract spaces with
prerings $V_i.$
Assume that $(X_i,V_i,v_i)$ for $i=1,2$ form positive measure spaces.
Let $(X,V,v)$ denote the tensor product measure
space $$(X_1\times X_2,V_1\otimes V_2,v_1\otimes v_2).$$

If $f_1\in L(v_1,R)$ and
$f_2\in N(v_2,R)$ then the function $f$ defined by the formula
\begin{equation*}
    f(x_1,x_2)=f_1(x_1)\ f_2(x_2)\qtext{ for }(x_1,x_2)\in X_1\times X_2
\end{equation*}
belongs to the set $N(v,R)$.
\end{cor}

\bigskip

\bp
    It follows from the definition of the space $L(v,R)$ of
    summable functions that for every function
    \begin{equation*}
        f_1\in L(v_1,R)
    \end{equation*}
    there exists a set $A\in (V_1)^\sigma$ such that
    \begin{equation*}
        f_1(x)=0\qtext{ if }x\notin A.
    \end{equation*}

    From the definition of a null function we have
    \begin{equation*}
        f_2(x)=0\ \qtext{ if }x\notin B
    \end{equation*}
    for some $v_2$-null set $B$. This implies
    \begin{equation*}
        f(x_1,x_2)=f_1(x_1)f_2(x_2)=0\qtext{ if }(x_1,x_2)\notin A\times B.
    \end{equation*}
    According to Theorem \ref{fub-Theorem 4} the set $A\times B$ is
    $v$-null set. Thus we have
    \begin{equation*}
        f\in N(v,R).
    \end{equation*}
    This concludes the proof of the corollary.
\ep

\bigskip

Let $v$ be a measure on a prering $V$. Denote by
\begin{equation*}
S=S(V)
\end{equation*}
the family of simple sets generated by the prering. Let $S_\delta$
be the family of all sets of the form
\begin{equation*}
A=\bigcap_nA_n
\end{equation*}
where $A_n$ is a sequence of simple sets. It is easy to see that
\begin{equation*}
S_\delta\subset V_c
\end{equation*}

Denote by $V_v$ the family of all sets of the form
\begin{equation*}
A=\bigcup_nA_n
\end{equation*}
where $A_n$ is an increasing sequence of sets from the family
$S_\delta$ such that the sequence of numbers $v_c(A_n)$ is
bounded. Let $N_v$ denote the family of all $v$-null sets.

\bigskip

It follows from a result in the paper of
Bogdanowicz~\cite{bogdan14}, page~259, Theorem~4, Part~6,
that a set $A$ is summable
if and only if there
exist a set $B\in V_v$ and a set $C\in N_v$ such that the set $A$
can be represented as the symmetric difference
\begin{equation*}
       A=B\div C=(B\setminus C)\cup(C\setminus B).
\end{equation*}

\bigskip

\begin{thm}
Assume that the triples $(X,W,w)$ and $(Y,U,u)$ represent measure
spaces over prerings $W,U,$ respectively. Let $(Z,V,v)$ denote the
product measure space $$(X\times Y,W\otimes U,w\otimes u).$$
Then the measure $v_c$ is an extension of the product measure
$\rho=w_c\otimes u_c$ and moreover $\rho_c=v_c$.
\end{thm}

\bigskip

\bp
    Denote by $V,W,U,$ respectively, the domains of the measures
    $v,w,u.$ Take two sets $A\in W_c$ and $B\in U_c$. These sets can
    be represented in the form
    \begin{equation*}
    A=C\div D\qtext{ and }B=E\div F
    \end{equation*}
    where
    \begin{equation*}
    C\in W_w,\ D\in N_w,\ E\in U_u,\ F\in N_u.
    \end{equation*}
    It is plain that
    \begin{equation*}
    G\in V_v\subset V_c,\qtext{ where }G=C\times E.
    \end{equation*}
    This implies
    \begin{equation*}
    c_G\in L(v,R).
    \end{equation*}

    \bigskip

    Notice that the condition $A=C\div D$ is equivalent to $A\div C=D$
    which in turn is equivalent to the condition $c_A-c_C=h\in
    N(w,R)$. Similarly one can prove that $c_B-c_E=g\in N(u,R)$.
    \\
    Using the identity
    \begin{equation*}
    c_{A\times B}(x_1,x_2)=c_A(x_1)c_B(x_2)\qtext{ for all }x_1\in X_1,x_2\in X_2
    \end{equation*}
    we easily get the representation
    \begin{equation*}
    c_{A\times B}=c_{C\times E}+k
    \end{equation*}
    where \begin{equation*}
    k=hg+hc_E+gc_C.
    \end{equation*}
    It follows from Corollary 1 that the function $k$ is a $v$-null
    function. This means that the functions
    \begin{equation*}
    c_{A\times B}\ ,\ c_{C\times E}=c_G
    \end{equation*}
    are equal $v$-almost everywhere. Since the function $c_G$ is a
    $v$-summable function, therefore according to Theorem \ref{Theorem
    1} the function $c_{A\times B}$ is also summable, that is we have
    \begin{equation*}
    A\times B\in V_c.
    \end{equation*}

    Now from Fubini-Bochner Theorem \ref{fub-Bochner} we get
    \begin{equation*}
    w_c(A)u_c(B)=\int c_{A\times B}dv=v_c(A\times B)\qtext{for all }A\times B\in V.
    \end{equation*}

    The relation $\eta\subset\mu$ between two functions
    will mean that the graph of the function $\eta$
    is a subset of the graph of the function $\mu.$
    Thus we have proved
    \begin{equation*}
    \rho=w_c\cdot u_c\subset v_c.
    \end{equation*}
    Since the measure $\rho_c$ is the smallest complete measure being
    an extension of the measure $\rho$ and the measure $v_c$ is
    complete we get from the previous relation the relation
    $\rho_c\subset v_c$. Similarly from the relation $v\subset \rho
    \subset \rho_c$ we get $v_c\subset\rho_c$. This proves $\rho=v_c$.
\ep

\bigskip

Since the completions of the measures $v,\ \rho,\ v_c$ coincide
therefore according to Bogdanowicz~\cite{bogdan14}, page~267, Section~7, Theorem~8,
they generate the same
Lebesgue integral as the product measure.

\bigskip

\begin{cor}The measures:
$v=w\otimes u,$ $\rho=w_c\otimes u_c,$ and $v_c,$ all generate the
same Lebesgue-Bochner integration.
\end{cor}

\bigskip


\bigskip

\section{Tonelli's theorem for the extended Lebesgue integral.}

\bigskip

Let $X$ be an abstract space with a prering $V$ of sets and
$Y$  a Banach space. Assume that the triple $(X,V,v)$ forms a positive measure space.

We define the space $M(v,Y)$ of {\bf Bochner measurable functions}
as in Bogdanowicz \cite{bogdan14}.
A functions $f:X\to Y$ belongs to the space $M(v,Y)$
if it has its support in  a set $B\in V^\sigma,$ that is
\begin{equation*}\label{ton-(1)}
    \quad f(x)=0\qtext{ if }x_1\notin B,
\end{equation*}
and
\begin{equation*}\label{ton-(2)}
    \quad c_Aa(f(\cdot))\in L(v,Y)\qtext{ for all } A\in V,
\end{equation*}
where
\begin{equation*}
    a(y)=(1+|y|)^{-1}y\qtext{ for all }y\in Y.
\end{equation*}

This definition allows us, using the properties of the space of
 Lebesgue-Bochner summable functions $L(v,Y)$ to get the
 classical properties of, and relations between, the spaces
\begin{equation*}
L(v,Y),\ M(v,Y),\ L(v,R),\ M(v,R).
\end{equation*}
For details see Bogdanowicz \cite{bogdan14}.

Let $M^+(v)$ denote the space of all functions $f$
from the set $X$ into the closed extended interval $[ 0,\infty]$
for which there exist a set $B\in V_\sigma$ such that
\begin{equation*}
     f(x)=0\qtext{ if }x\notin B,
\end{equation*}
and
\begin{equation*}
c_Ab(f(\cdot))\in L(v,R)\qtext{ for all } A\in V,
\end{equation*}
where
\begin{equation*}
b(y)=(1+y)^{-1}y\qtext{ if }y\in [ 0,\infty),\qtext{ and }b(\infty)=1.
\end{equation*}

Functions belonging to the space $M^+(v)$ will be called Lebesgue measurable
extended functions.

It was proven in \cite{bogdan14} that $f\in M^+(v)$ if and
only if there exist non-negative Lebesgue summable functions
    $f_n\in L(v,R)$
such that the sequence $f_n$ increasingly converges almost everywhere to the
function $f$. This allows us to extend the integral onto $M^+(v)$ by the formula
\begin{equation*}
\int f\,dv=\lim_n \int f_n\,dv.
\end{equation*}
It follows from the Monotone Convergence Theorem that this
extension of the integral onto the space $M^+(v)$ does not depend
on the choice of the sequence $f_n.$ Thus the definition is correct.

The following implications have been proved in \cite{bogdan14}
\begin{equation*}
    f\in M(v,Y)\impl |f(\cdot)|\in M^+(v)
\end{equation*}
and
\begin{equation*}
    \{f\in M(v,Y)\text{ and }\int |f(\cdot)|\,dv<\infty\}\impl f\in L(v,Y).
\end{equation*}

\bigskip

\begin{thm}[Tonelli]\label{ton-Tonelli}
Assume that $X_i$ are abstract spaces with
prerings $V_i.$
Assume that $(X_i,V_i,v_i)$ for $i=1,2$ form positive measure spaces.
Let $(X,V,v)$ denote the tensor product measure
space $$(X_1\times X_2,V_1\otimes V_2,v_1\otimes v_2).$$

 Then for every extended measurable function $f\in M^+(v)$ there exists a
$v_1$-null set $C$ and an extended measurable function $h\in M^+(v_1)$ such that
\begin{equation*}
    f(x_1,\cdot)\in M^+(v_2)\qtext{ and }
    h(x_1)=\int f(x_1,x_2)\,v_2(dx_2)\qtext{ if }x_1\notin C,
\end{equation*}
and moreover
\begin{equation*}
    \int f\, dv=\int h\, dv_1=\int\left(\int f(x_1,x_2)\,v_2(dx_2)\right)\,v_1(dx_1).
\end{equation*}
\end{thm}

\bigskip

\bp
Take any function $f\in M^+(v)$.
There exists a sequence of non-negative Lebesgue summable functions $f_n\in L(v,R)$
that increasingly converges to the function $f$ except on a $v$-null set $A$.
Denote by $C_0$ a $v_1$-null set such that for each $x_1\notin C_0$, the section $A(x_1)$
of the set $A$ forms a $v_2$-null set.

It follows from Fubini's theorem that there exist functions
$h_n\in L(v_1,R)$ and a sequence $C_n$ of $v_1$-null sets such that
\begin{equation*}
f_n(x_1,\ )\in L(v_2,R)\qtext{ and }h_n(x_1)=
\int f_n(x_1,x_2)\,dv_2\qtext{ if } x_1\notin C_n.
\end{equation*}
It is easy to see that one may assume that the functions $h_n$ are non-negative
and
\begin{equation*}
    \int h_n\,dv_1=\int f_n\,dv\fa n=1,2,\ldots
\end{equation*}
Put
\begin{equation*}
C=\bigcup^{\,\infty}_{n=0}C_n.
\end{equation*}
If $x_1\notin C$ then the sequence of functions $f_n(x_1,\cdot )$
increasingly converges to the function $f_n(x_1,\cdot )$ for $x_2\notin A(x_1),$
where $A(x_1)$ represents a $v_2$-null set. This implies
\begin{equation*}
    f(x_1,\cdot)\in M^+(v_2)
\end{equation*}
and
\begin{equation*}
    \int f(x_1,\cdot)\,dv_2=\lim_n\int f_n(x_1,\cdot)\,dv_2.
\end{equation*}

If $x_1\notin C$ then the sequence
\begin{equation*}
    h_n(x_1)=\int f_n(x_1,\cdot)\,dv_2
\end{equation*}
increasingly converges to a value $h(x_1)\in [ 0,\infty].$
Put $h(x_1)=0$ for $x_1\in C$.
Notice that
\begin{equation*}
    h(x)=\lim_n h_n(x_1)=\lim_n\int f_n(x_1,\cdot)\,dv_2
    =\int f(x_1,\cdot)\,dv_2\qtext{ if } x_1\notin C.
\end{equation*}
Since the sequence $h_n\in L(v_1, R)$ increasingly
converges almost everywhere to the
function $h$ therefore we have $h\in M^+(v_1)$ and
\begin{equation*}
    \int h\,dv_1=\lim_n\int h_n\,dv_1=\lim_n\int f_n\,dv=\int f\, dv.
\end{equation*}
Thus the theorem is proved.
\ep

\bigskip


\bigskip

\section{Extension of vector measure}

\bigskip


Now let us consider the case when the space $Y$ is the space $R$
of reals and $Z$ any Banach space and the bilinear operator $u$ is
the multiplication operator $u(r,z)=rz.$

\begin{pro}[Isomorphism and isometry of $K(v,Z)$ and $K(v_c,Z)$]
Assume that $(X,V,v)$ is  a measure space on a prering $V$ of
subsets of an abstract space $X.$

Every vector measure $\mu\in K(v,Z)$ can be extended from the
prering $V$ onto the delta ring $V_c$ of summable sets by the formula
\begin{equation*}
    \mu_c(A)=\int u(c_A,d\mu)\fa A\in V_c.
\end{equation*}
This extension establishes isometry and isomorphism between the
Banach spaces $K(v,Z)$ and $K(v_c,Z).$
\end{pro}

\bigskip


\bigskip

\section{Absolute continuity of vector measures}

\bigskip


\begin{defin}[Absolute continuity of a vector measure]
\label{Absolute continuity of a vector measure} Given a Lebesgue
measure space $(X,V,v)$ on a $\sigma$-algebra $V$ and
$\sigma$-additive vector measure $\mu$  from $V$ into a Banach
space $Y.$ We say that the vector measure $\mu$ is absolutely
continuous with respect to the measure $v$ if
\bb
    \mu(A)=0\qtext{whenever} v(A)=0.
\ee
\end{defin}

\bigskip

\begin{thm}[Phillips]
\label{Phillips} Assume that $(X,V,v)$ is a Lebesgue measure space
on a $\sigma$-algebra $V$ such that $v(X)<\infty$ and $Y$ is a
reflexive Banach space.

Assume that $\mu$ is a $\sigma$-additive vector measure of finite
variation from the $\sigma$-algebra $V$ into $Y.$

If $\mu$ is absolutely continuous with respect to the Lebesgue
measure $v,$ then there exists a Bochner summable function $g\in
L(v,Y)$ such that \bb
    \mu(A)=\int_Ag\,dv\fa A\in V.
\ee
\end{thm}

\bigskip

For the proof of this remarkable theorem see Diestel and Uhl
\cite[page 76]{diestel-uhl}. This result can be found in the
original paper of Phillips \cite{phillips}.

\bigskip

\end{document}